\pdfoutput=1

\documentclass[final]{siamart171218}
\usepackage[T1]{fontenc}
\usepackage[latin2]{inputenc}
\usepackage{amsmath,amssymb,amsbsy,amsfonts}
\usepackage{mathtools}
\usepackage{tikz}
\usepackage{enumerate}
\usepackage{hyperref}
\usepackage{siunitx}
\usepackage{ragged2e}
\usepackage{array}
\usepackage{cprotect}
\usepackage{fnpct}
\usepackage{stackengine}
 
\usepackage[textsize=tiny]{todonotes}

\usepackage{etoolbox}
\makeatletter
\patchcmd{\@addmarginpar}{\ifodd\c@page}{\ifodd\c@page\@tempcnta\m@ne}{}{}
\makeatother

\usepackage{marginnote}
\makeatletter
\renewcommand{\@todonotes@drawMarginNoteWithLine}{%
  \begin{tikzpicture}[remember picture, overlay, baseline=-0.75ex]%
      \node [coordinate] (inText) {};%
  \end{tikzpicture}%
  \marginnote[{
      \@todonotes@drawMarginNote%
      \@todonotes@drawLineToLeftMargin%
  }]{
      \@todonotes@drawMarginNote%
      \@todonotes@drawLineToRightMargin%
  }%
}
\makeatother

\usepackage{etoolbox}
\makeatletter
\patchcmd{\@mn@@@marginnote}{\begingroup}{\begingroup\@twosidefalse}{}{\fail}
\makeatother
 
\def\tens#1{\pmb{\mathbb{#1}}}
\def\vec#1{\boldsymbol{#1}}

\def\supp{\mathop{\operatorname{supp}}\nolimits}
\def\esssup{\mathop{\operatorname{ess\,sup}}}
\def\diver{\mathop{\operatorname{div}}\nolimits}

\def\llparen{\left(\kern-.25em\left(}
\def\rrparen{\right)\kern-.25em\right)}

\def\tr{\mathop{\operatorname{tr}}\nolimits}

\def\lin{\mathop{\operatorname{span}}\nolimits}
\def\d{\mathrm{d}}
\def\pp#1#2{\frac{\partial #1}{\partial #2}}
\def\dd#1#2{\frac{\d #1}{\d #2}}

\def\RR{\mathbb{R}}
\def\bx{\vec{x}}

\def\bb{\vec{b}}

\def\bs{\vec{s}}
\def\ba{\vec{a}}
\def\bn{\vec{n}}
\def\bm{\vec{m}}

\def\bz{\vec{z}}

\def\btau{\vec{\tau}}
\def\bc{\vec{c}}

\def\bh{\vec{h}}

\def\bF{\vec{F}}

\def\bA{\tens{A}}
\def\bT{\tens{T}}
\def\bQ{\tens{Q}}
\def\bD{\tens{D}}

\def\bO{\tens{O}}

\def\bS{\tens{S}}
\def\bI{\tens{I}}
\def\bP{\tens{P}}
\def\bH{\tens{H}}

\def\bv{\vec{v}}
\def\bfi{\vec{\varphi}}

\def\b0{\vec{0}}
\def\bom{\vec{\omega}}
\def\bw{\vec{w}}

\def\bu{\vec{u}}
\def\bU{\vec{U}}

\def\nug{\nu_\mathrm{g}} 
\def\alg{\alpha_\mathrm{g}} 
\def\cG{\mathcal{G}} 
\def\cB{\mathcal{B}} 

\newcommand{\Cinf}{\mathcal{C}^\infty}
\newcommand{\D}{\mathcal{C}^\infty_{0}}

\newcommand{\bC}{\boldsymbol{\mathcal{C}}}
\newcommand{\bL}{\mathbf{L}}
\newcommand{\bW}{\mathbf{W}}

\newcommand{\Lndiv}[1]{\bL^{#1}_{\bn,\diver}}

\newcommand{\Cinfz}{\bC^\infty_{0}}
\newcommand{\Cinfzdiv}{\bC^\infty_{0,\diver}}
\newcommand{\Vz}{\bW^{3,2}_{0}}
\newcommand{\Vzdiv}{\bW^{3,2}_{0,\diver}}
\newcommand{\Wz}[1]{\bW^{1,#1}_{0}}
\newcommand{\Wzdiv}[1]{\bW^{1,#1}_{0,\diver}}

\newcommand{\Cinfn}{\bC^\infty_{\bn}}
\newcommand{\Cinfndiv}{\bC^\infty_{\bn,\diver}}
\newcommand{\Vn}{\bW^{3,2}_{\bn}}
\newcommand{\Vndiv}{\bW^{3,2}_{\bn,\diver}}
\newcommand{\Wn}[1]{\bW^{1,#1}_{\bn}}
\newcommand{\Wndiv}[1]{\bW^{1,#1}_{\bn,\diver}}

\newcommand{\freeslip}{free-slip}

\newcommand{\navierslip}{Navier-slip}
\newcommand{\navierslipalt}{slip}
\newcommand{\noslip}{no-slip}
\newcommand{\noslipalt}{stick}
\newcommand{\slip}{slip}

\newcommand{\fixsupscript}{\smash[t]{\vphantom{)}}}
\newcommand{\fixsupscriptmore}{\smash[t]{\vphantom{\bigr)}}}
 
\newsiamthm{Theorem}{Theorem}
\newsiamthm{Lem}{Lemma}
\newsiamremark{Rem}{Remark}
\newsiamthm{Col}{Corollary}

\def\figtau{0.382}

\tikzstyle{dot}=[draw,fill=white,circle,inner sep=0pt,minimum size=4pt]

\newcommand{\smallfigscale}{0.75}
\newcommand{\smallfig}[2]{%
  \begin{tikzpicture}[scale=\smallfigscale,baseline=15.0]%
    \clip (-0.10, -0.25) rectangle (1.10, 1.30);%
    \draw[ultra thin] (-0.1,  0.0) -- (1.1, 0.0);%
    \draw[ultra thin] ( 0.0, -0.1) -- (0.0, 1.1);%
    \draw[variable=\t,domain=0:1,smooth,thick] #1;%
    \draw[variable=\t,domain=0:1.1,dashed,thin] #2;%
  \end{tikzpicture}%
}

\newcommand{\bigfigscale}{2.5}
\newcommand{\bigfig}[1]{%
  \begin{tikzpicture}[scale=\bigfigscale]%
    \footnotesize
    \draw[ultra thin] (-0.1,  0.0) -- (1.1, 0.0) node[below]{$|\bD|$};%
    \draw[ultra thin] ( 0.0, -0.1) -- (0.0, 1.1) node[left ]{$|\bS|$};%
    \draw[variable=\t,domain=0:1,smooth,ultra thick] #1;%
  \end{tikzpicture}%
}

\def\shearlim {\smallfig{plot ({tanh(4*\t)},{\t        })}{plot ({1}, {\t})}}
\def\thick    {\smallfig{plot ({\t^0.5    },{\t        })}{}}
\def\ns       {\smallfig{plot ({\t        },{\t        })}{}}
\def\thin     {\smallfig{plot ({\t        },{\t^0.5    })}{}}
\def\stresslim{\smallfig{plot ({\t        },{tanh(4*\t)})}{plot ({\t}, {1})}}

\def\rigid     {\smallfig{plot[domain=0:\figtau] ({0},{\t}) -- plot[domain=0:1-\figtau] ({0     },{\figtau+\t})}{}}
\def\rigidthick{\smallfig{plot[domain=0:\figtau] ({0},{\t}) -- plot[domain=0:1-\figtau] ({\t^0.5},{\figtau+\t})}{}}
\def\rigidlin  {\smallfig{plot[domain=0:\figtau] ({0},{\t}) -- plot[domain=0:1-\figtau] ({\t    },{\figtau+\t})}{}}
\def\rigidthin {\smallfig{plot[domain=0:\figtau] ({0},{\t}) -- plot[domain=0:1-\figtau] ({\t^2  },{\figtau+\t})}{}}
\def\rigidlim  {\smallfig{plot[domain=0:\figtau] ({0},{\t}) -- plot[domain=0:1        ] ({\t    },{\figtau   })}{}}

\def\eulerlim  {\smallfig{plot[domain=0:\figtau] ({\t},{0}) -- plot[domain=0:1        ] ({\figtau   },{\t    })}{}}
\def\eulerthick{\smallfig{plot[domain=0:\figtau] ({\t},{0}) -- plot[domain=0:1-\figtau] ({\figtau+\t},{\t^2  })}{}}
\def\eulerlin  {\smallfig{plot[domain=0:\figtau] ({\t},{0}) -- plot[domain=0:1-\figtau] ({\figtau+\t},{\t    })}{}}
\def\eulerthin {\smallfig{plot[domain=0:\figtau] ({\t},{0}) -- plot[domain=0:1-\figtau] ({\figtau+\t},{\t^0.5})}{}}
\def\euler     {\smallfig{plot[domain=0:\figtau] ({\t},{0}) -- plot[domain=0:1-\figtau] ({\figtau+\t},{0     })}{}}

\newcommand{\smallfigcolorscale}{0.75}
\newcommand{\smallfigcolor}[2]{%
  \begin{tikzpicture}[scale=\smallfigcolorscale,baseline=15.0]%
    \clip (-0.15, -0.22) rectangle (1.15, 1.15);%
    \draw[ultra thin] (-0.1,  0.0) -- (1.1, 0.0);%
    \draw[ultra thin] ( 0.0, -0.1) -- (0.0, 1.1);%
    \draw[variable=\t,domain=0:1.1,smooth,very thick,dotted] #1;%
    \draw[variable=\t,domain=0:1,smooth,ultra thick,red] #2;%
  \end{tikzpicture}%
}

\def\nsrest       {\smallfigcolor  {plot ({\t        },{\t        })} {node[dot] at (0,0) {}}}
\def\nsmove       {\smallfigcolor{}{plot ({\t        },{\t        })}                        }
\def\rigidrest    {\smallfigcolor  {plot ({0         },{\t        })} {node[dot] at (0,0) {}}}
\def\rigidmove    {\smallfigcolor{}{plot ({0         },{\t        })}                        }

\def\eulermove    {\smallfigcolor{}{plot ({\t        },{0         })}                        }
\def\rigidlinrest {\smallfigcolor{plot[domain=0:\figtau] ({0},{\t}) -- plot[domain=0:1-\figtau] ({\t},{\figtau+\t})} {node[dot] at (0,0) {}}}
\def\rigidlinone  {\smallfigcolor{plot[domain=0:1-\figtau] ({\t},{\figtau+\t})}{plot[domain=0:\figtau] ({0},{\t})}}
\def\rigidlintwo  {\smallfigcolor{plot[domain=0:\figtau] ({0},{\t})}{plot[domain=0:1-\figtau] ({\t},{\figtau+\t})}}

\def\eulerlinone  {\smallfigcolor{plot[domain=0:1-\figtau] ({\figtau+\t},{\t})}{plot[domain=0:\figtau] ({\t},{0})}}
\def\eulerlintwo  {\smallfigcolor{plot[domain=0:\figtau] ({\t},{0})}{plot[domain=0:1-\figtau] ({\figtau+\t},{\t})}}
\def\eulerlinthree{\smallfigcolor{plot[domain=0:\figtau] ({\t},{0})  plot[domain=0.5*(1-\figtau):1-\figtau]  ({\figtau+\t},{\t})}{plot[domain=0:0.5*(1-\figtau)] ({\figtau+\t},{\t})}}
\def\eulerlinfour {\smallfigcolor{plot[domain=0:\figtau] ({\t},{0})--plot[domain=0:0.5*(1-\figtau)]  ({\figtau+\t},{\t})}{plot[domain=0.5*(1-\figtau):1-\figtau] ({\figtau+\t},{\t})}}

\def\bigeuler  {\bigfig{plot ({\t},{0 })}}
\def\bigns     {\bigfig{plot ({\t},{\t})}}
\def\bigrigid  {\bigfig{plot ({0 },{\t})}}
\def\bigbingham{\bigfig{plot[domain=0:\figtau] ({0},{\t}) -- plot[domain=\figtau:1] ({\t-\figtau},{\t});
                        \node[dot,label={left:$\sigma_*$}] at (0, \figtau) {}}}
\def\bigeulerns{\bigfig{plot[domain=0:\figtau] ({\t},{0}) -- plot[domain=\figtau:1] ({\t},{\t-\figtau})
                        node[inner sep=1,outer sep=1,draw,ultra thin,right,xshift=3,yshift=0] {$r=2$}
                        plot[domain=\figtau:1] ({\t}, {(\t^0.5/\figtau^0.5)*(\t-\figtau)})
                        node[inner sep=1,outer sep=1,draw,ultra thin,right,xshift=3,yshift=0] {$r=\frac52$}
                        ;
                        \node[dot,label={below:$\delta_*$}] at (\figtau, 0) {}}}
\def\bigeulerlad{\bigfig{plot[domain=0:\figtau] ({\t},{0}) --
                         plot[domain=\figtau:1] ({\t}, {(0.5+0.5*(\t^0.5/\figtau^0.5))*(\t-\figtau)})
                         ;
                         \node[dot,label={below:$\delta_*$}] at (\figtau, 0) {}}}

\def\bigpowerlawscale{1.5}
\def\bigpowerlaw{
  \begin{tikzpicture}[scale=\bigpowerlawscale]
    \footnotesize
    \draw[ultra thin] (-0.2,  0.0) -- (2.2, 0.0) node[below]{$|\bD|$};
    \draw[ultra thin] ( 0.0, -0.2) -- (0.0, 2.2) node[left ]{$|\bS|$};
    \node[dot,label={below:$d_*$}]      at (1, 0) {};
    \node[dot,label={left:$2\nu_*d_*$}] at (0, 1) {};
    \draw[variable=\t,domain=0:2,thick,samples=400]
      plot ({\t}         ,{(\t)^(1/19)}) node[draw,ultra thin,right,xshift=3,yshift=-3] {$r=\frac{20}{19}$};
    \draw[variable=\t,domain=0:2,thick,samples=400]
      plot ({\t}         ,{(\t)^(1/2)} ) node[draw,ultra thin,right,xshift=3,yshift=3] {$r=\frac{3}{2}$};
    \draw[variable=\t,domain=0:2,thick,samples=400]
      plot ({\t}         ,{\t}         ) node[draw,ultra thin,above right,xshift=1.6,yshift=1.6] {$r=2$};
    \draw[variable=\t,domain=0:2,thick,samples=400]
      plot ({(\t)^(1/2)} ,{\t}         ) node[draw,ultra thin,above,yshift=3,xshift=8] {$r=3$};
    \draw[variable=\t,domain=0:2,thick,samples=400]
      plot ({(\t)^(1/19)},{\t}         ) node[draw,ultra thin,above,yshift=3,xshift=-8] {$r=20$};
  \end{tikzpicture}
}

\def\biglimitingstressscale{1.5}
\def\biglimitingstress{
  \begin{tikzpicture}[scale=\biglimitingstressscale]
    \footnotesize
    \draw[ultra thin] (-0.2,  0.0) -- (2.2, 0.0) node[below]{$|\bD|$};
    \draw[ultra thin] ( 0.0, -0.2) -- (0.0, 2.2) node[left ]{$|\bS|$};
    \node[dot,label={below:$d_*$}]      at (1, 0) {};
    \node[dot,label={left:$2\nu_*d_*$}] at (0, 1) {};
    \draw[variable=\t,domain=0:2,thick,samples=400,every node/.style={inner sep=1,outer sep=1}]
      plot ({\t},{(1+(\t)^(1/2))^(-2)*\t}) node[draw,ultra thin,right,xshift=3,yshift= 0.0] {$a=\frac12$}
      plot ({\t},{(1+(\t)^(1/1))^(-1)*\t}) node[draw,ultra thin,right,xshift=3,yshift= 0.0] {$a=1$}
      plot ({\t},{(1+(\t)^(2))^(-1/2)*\t}) node[draw,ultra thin,right,xshift=3,yshift=-1.5] {$a=2$}
      plot ({\t},{(1+(\t)^(4))^(-1/4)*\t}) node[draw,ultra thin,right,xshift=3,yshift= 2.5] {$a=4$}
      ;
    \draw[variable=\t,dashed,every node/.style={inner sep=1,outer sep=1}]
      plot[domain=0:1] ({\t},{\t}) node[draw,ultra thin,solid,above,xshift=0,yshift=3] {$a=\infty$}
      --plot[domain=1:2] ({\t},{1});
  \end{tikzpicture}
}
\def\biglimitingstrainscale{1.5}
\def\biglimitingstrain{
  \begin{tikzpicture}[scale=\biglimitingstrainscale]
    \footnotesize
    \draw[ultra thin] (-0.2,  0.0) -- (2.2, 0.0) node[below]{$|\bD|$};
    \draw[ultra thin] ( 0.0, -0.2) -- (0.0, 2.2) node[left ]{$|\bS|$};
    \node[dot,label={below:$d_*$}]      at (1, 0) {};
    \node[dot,label={left:$2\nu_*d_*$}] at (0, 1) {};
    \draw[variable=\t,domain=0:2,thick,samples=400,every node/.style={inner sep=1,outer sep=1}]
      plot ({(1+(\t)^(1/2))^(-2)*\t},{\t}) node[draw,ultra thin,above,yshift=3,xshift= 0.0] {\rotatebox{90}{$b=\frac12$}}
      plot ({(1+(\t)^(1/1))^(-1)*\t},{\t}) node[draw,ultra thin,above,yshift=3,xshift=-1.0] {\rotatebox{90}{$b=1$}}
      plot ({(1+(\t)^(2))^(-1/2)*\t},{\t}) node[draw,ultra thin,above,yshift=3,xshift=-2.5] {\rotatebox{90}{$b=2$}}
      plot ({(1+(\t)^(4))^(-1/4)*\t},{\t}) node[draw,ultra thin,above,yshift=3,xshift= 2.5] {\rotatebox{90}{$b=4$}}
      ;
    \draw[variable=\t,dashed,every node/.style={inner sep=1,outer sep=1}]
      plot[domain=0:1] ({\t},{\t}) node[draw,ultra thin,solid,right,yshift=0,xshift=3] {$b=\infty$}
      --plot[domain=1:2] ({1},{\t});
  \end{tikzpicture}
}

\pgfkeys{
  /pgf/declare function={Si1(\r,\nu,\d,\D)=2*\nu*\d*(\D/\d)^(\r-1);},
  /pgf/declare function={Si2(\r,\nu,\d,\D)=2*\nu*(1+\D^2/\d^2)^(\r/2-1)*\D*2^(1-\r/2);},
  /pgf/declare function={Si(\r,\q,\nuratio,\D)=Si1(\r,1/(1+\nuratio),1,\D)+Si2(\q,\nuratio/(1+\nuratio),1,\D);},
  /pgf/declare function={Sii1(\nu,\D)=2*\nu*\D;},
  /pgf/declare function={Dii(\rp,\nu,\d,\Stwo)=1/(2*\nu)*(1+(\Stwo)^2/(2*\nu*\d)^2)^(\rp/2-1)*\Stwo*2^(1-\rp/2);},
}
\def\bigaddstressixscale{1.5}
\def\bigaddstressiyscale{0.75}
\def\bigaddstressi{
  \begin{tikzpicture}[xscale=\bigaddstressixscale,yscale=\bigaddstressiyscale]
    \footnotesize
    \draw[ultra thin] (-0.2,  0.0) -- (2.2, 0.0) node[below]{$|\bD|$};
    \draw[ultra thin] ( 0.0, -0.2) -- (0.0, 4.6) node[left ]{$|\bS|$};
    \node[dot,label={below:$d_*$}]      at (1, 0) {};
    \node[dot,label={left:$2(\nu_*+\tilde\nu_*)d_*$}] at (0, 2) {};
    \draw[variable=\t,domain=0:2,thick,samples=400,every node/.style={inner sep=1,outer sep=1}]
      plot ({\t},{Si(1.5,2.5,0.1250,\t)}) node[draw,ultra thin,right,xshift=  3,yshift=  0] {$\frac{\tilde\nu_*}{\nu_*}=\frac18$}
      plot ({\t},{Si(1.5,2.5,1.0000,\t)}) node[draw,ultra thin,right,xshift=  3,yshift=  0] {$\frac{\tilde\nu_*}{\nu_*}=1$}
      plot ({\t},{Si(1.5,2.5,8.0000,\t)}) node[draw,ultra thin,right,xshift=  3,yshift=  0] {$\frac{\tilde\nu_*}{\nu_*}=8$}
      ;
  \end{tikzpicture}
}
\def\bigaddstressiixscale{1.5}
\def\bigaddstressiiyscale{0.375}
\def\bigaddstressii{
  \begin{tikzpicture}[xscale=\bigaddstressiixscale,yscale=\bigaddstressiiyscale]
    \footnotesize
    \draw[ultra thin] (-0.2,  0.0) -- (2.2, 0.0) node[below]{$|\bD|$};
    \draw[ultra thin] ( 0.0, -0.2) -- (0.0, 9.2) node[left ]{$|\bS|$};
    \node[dot,label={below:$d_*$}]      at (1, 0) {};
    \node[dot,label={left:$2(\nu_*+\tilde\nu_*)d_*$}] at (0, 4) {};
    \draw[variable=\t,domain=0:4,thick,samples=100,every node/.style={inner sep=1,outer sep=1}]
      plot[domain=0:4.935] ({Dii(10/9,1.0,1.0,\t)},{\t+Sii1(1.0,Dii(1.25,1.0,1.0,\t))}) node[draw,ultra thin,above,xshift=  0,yshift= 3] {$r'=\frac{10}9$}
      plot[domain=0:3.990] ({Dii(2.0 ,1.0,1.0,\t)},{\t+Sii1(1.0,Dii(2.0 ,1.0,1.0,\t))}) node[draw,ultra thin,right,xshift=  3,yshift= 0] {$r'=2$}
      plot[domain=0:2.282] ({Dii(10. ,1.0,1.0,\t)},{\t+Sii1(1.0,Dii(5.0 ,1.0,1.0,\t))}) node[draw,ultra thin,right,xshift=  3,yshift= 0] {$r'=10$}
      ;
  \end{tikzpicture}
}

\def\simpleflows{
  \def\simpleflowsxscale{3.0}
  \def\simpleflowsyscale{1.5}
  \begin{tikzpicture}[xscale=\simpleflowsxscale,yscale=\simpleflowsyscale,declare function={
    ypoint(\eps,\C,\delta)=2^0.5*\delta*\eps/(2*\C);
    ueps(\eps,\y,\C,\delta)=
      + (abs(\y) < ypoint(\eps,\C,\delta))
          * -\C/(\eps)
          * \y
          * \y
      + (abs(\y) >= ypoint(\eps,\C,\delta))
          * -\C/(1+\eps)
          * (\y+sign(\y)*2^0.5/(2*\C)*(1+(1+\eps)^0.5))
          * (\y+sign(\y)*2^0.5/(2*\C)*(1-(1+\eps)^0.5))
    ;
    u(\y,\C,\delta)= -\C*(\y)^2 - 2^0.5*\delta*abs(\y);
  }]
    \footnotesize
    \draw[ultra thin] (-1.1,-2.6) -- ( 1.1,-2.6) node[below]{$y$};
    \draw[ultra thin] (-1.1,-2.6) -- (-1.1, 0.1) node[left ]{$u$};
    \node[dot,label={below:$y_0$}] at (0,-2.6) {};
    \node[dot,xshift={ ypoint(1.0 ,1,1)*\simpleflowsxscale*1cm},yshift={ueps(1.0 , ypoint(1.0 ,1,1),1,1)*\simpleflowsyscale*1cm}] at (0,0) {};
    \node[dot,xshift={-ypoint(1.0 ,1,1)*\simpleflowsxscale*1cm},yshift={ueps(1.0 ,-ypoint(1.0 ,1,1),1,1)*\simpleflowsyscale*1cm}] at (0,0) {};
    \node[dot,xshift={ ypoint(0.5 ,1,1)*\simpleflowsxscale*1cm},yshift={ueps(0.5 , ypoint(0.5 ,1,1),1,1)*\simpleflowsyscale*1cm}] at (0,0) {};
    \node[dot,xshift={-ypoint(0.5 ,1,1)*\simpleflowsxscale*1cm},yshift={ueps(0.5 ,-ypoint(0.5 ,1,1),1,1)*\simpleflowsyscale*1cm}] at (0,0) {};
    \node[dot,xshift={ ypoint(0.25,1,1)*\simpleflowsxscale*1cm},yshift={ueps(0.25, ypoint(0.25,1,1),1,1)*\simpleflowsyscale*1cm}] at (0,0) {};
    \node[dot,xshift={-ypoint(0.25,1,1)*\simpleflowsxscale*1cm},yshift={ueps(0.25,-ypoint(0.25,1,1),1,1)*\simpleflowsyscale*1cm}] at (0,0) {};
    \node[dot] at (0,0) {};
    \draw[variable=\t,domain=-1:1,thick,samples=400,every node/.style={inner sep=1,outer sep=1}]
      plot ({\t},{ueps(1.0 ,\t,1,1)}) node[draw,ultra thin,right,xshift=3,yshift=0] {$\epsilon_*=1$}
      plot ({\t},{ueps(0.5 ,\t,1,1)}) node[draw,ultra thin,right,xshift=3,yshift=0] {$\epsilon_*=\frac12$}
      plot ({\t},{ueps(0.25,\t,1,1)}) node[draw,ultra thin,right,xshift=3,yshift=0] {$\epsilon_*=\frac14$}
      plot ({\t},{u   (     \t,1,1)}) node[draw,ultra thin,right,xshift=3,yshift=0] {$\epsilon_*=0$}
      ;
  \end{tikzpicture}
}
 \usepackage{tikz}
\usetikzlibrary{calc,decorations.pathmorphing,decorations.markings}

\tikzstyle{spring}=[thick,decorate,decoration={zigzag,pre length=0.3cm,post
  length=0.3cm,segment length=6}]

\tikzstyle{damper}=[thick,decoration={markings,
  mark connection node=dmp,
  mark=at position 0.5*\pgfdecoratedpathlength-10pt with
  {
    \node (dmp) [thick,inner sep=0pt,transform shape,rotate=-90,minimum
width=15pt,minimum height=8pt,draw=none] {};
    \draw [thick] ($(dmp.north east)+(20pt,0)$) -- (dmp.south east) -- (dmp.south
west) -- ($(dmp.north west)+(20pt,0)$);
    \draw [thick] ($(dmp.north)+(0,-5pt)$) -- ($(dmp.north)+(0,5pt)$);
  }
}, decorate]
 
\usetikzlibrary{external}
\makeatletter
\if@review
\renewcommand{\todo}[2][]{\tikzexternaldisable\@todo[#1]{#2}\tikzexternalenable}
\renewcommand{\missingfigure}[2][]{\tikzexternaldisable\@missingfigure[#1]{#2}\tikzexternalenable}
\else
\renewcommand{\todo}[2][]{}
\renewcommand{\missingfigure}[2][]{}
\fi
\makeatother

\title{On the classification of incompressible fluids
       and a~mathematical analysis of the equations that
       govern their motion\thanks{Submitted \today.
       \funding{Development of sections \ref{Ch0}, \ref{Ch1} was supported by the
       project LL1202 in the programme ERC-CZ funded by the Ministry of
       Education, Youth and Sports of the Czech Republic.
       Section~\ref{Ch2} was supported by the project GA \v{C}R 18-12719S
       funded by the Grant Agency of the Czech Republic.}}}

\author{%
  Jan Blechta\thanks{%
    Chemnitz University of Technology, Faculty of Mathematics
    (\email{jan.blechta@math.tu-chemnitz.de}).}
  \and
  Josef M\'{a}lek\thanks{%
    Charles University, Faculty of Mathematics and Physics
    (\email{malek@karlin.mff.cuni.cz}).}
  \and
  K.R. Rajagopal\thanks{%
    Texas A\&M University, Department of Mechanical Engineering
    (\email{krajagopal@tamu.edu}).}
}

\headers{On the classification of incompressible fluids}
        {J.~Blechta, J.~M\'{a}lek, and K.R.~Rajagopal}

\begin{document}

\makeatletter
\if@review\else\tikzexternaldisable\fi
\makeatother

\maketitle

\begin{abstract}
  In the first part of the paper we provide a~new classification of
incompressible fluids characterized by a~continuous monotone
relation between the velocity gradient and the Cauchy stress. The considered
class includes Euler fluids, Navier-Stokes fluids, classical power-law fluids
as well as stress power-law fluids, and their
various generalizations including the fluids that we refer to as activated fluids,
namely fluids that behave as an Euler fluid prior
activation and behave as a~viscous fluid once activation takes place. We
also present a~classification concerning boundary conditions that are
viewed as the constitutive relations on the boundary. In the second part of the
paper, we develop a~robust mathematical theory for activated Euler fluids
associated with different types of the boundary conditions ranging from no-slip to free-slip
and include Navier's slip as well as stick-slip. Both steady and
unsteady flows of such fluids in three-dimensional domains are analyzed.
 \end{abstract}



\section{Introduction} \label{Ch0}
The concept of a fluid defies precise definition as one can always come up with
a counter-example to that definition that seems to fit in with our
understanding of what constitutes a fluid. As Goodstein \cite{goodstein1985}
appropriately remarks ``Precisely what do we mean by the term liquid? Asking
what is a liquid is like asking what is life; we usually know when we see it,
but the existence of some doubtful cases make it hard to define precisely.'' The
concept of a fluid is treated as a primitive concept in mechanics, but
unfortunately it does not meet the fundamental requirement of a primitive, that
of being amenable to intuitive understanding. This makes the study under
consideration that much more difficult as it is our intent to classify fluid
bodies. In this study we shall consider a~subclass of the idealization of a~fluid,
namely that of incompressible fluid bodies. While no material is truly incompressible,
in many bodies the change of volume is sufficiently small to be ignorable. Our
ambit will include at one extreme materials that could be viewed as
incompressible Euler fluids and at the other extreme materials that offer so
much resistance to flow that they are ``rigid-like'' in their response, with a
whole host of ``fluid-like'' behavior exhibited by bodies whose response lie in
between these two extremes, such as fluids exhibiting shear thinning/shear
thickening, stress thinning/stress thickening, etc.

Before discussing the constitutive classification of fluid bodies, it would be
useful to consider another type of classification that is used, namely that of
flow classification with regard to the flows of a specific fluid, so that we do
not confuse these two types of classifications. One of the most useful
approximations and an integral part of fluid dynamics is the boundary layer
approximation for the flow of a Navier-Stokes fluid (see Prandtl
\cite{prandtl1905}, Schlichting \cite{schlichting1960}). The main tenet of the
approximation is the notion that for flows of a Navier-Stokes fluid past a
solid boundary, at sufficiently high Reynolds number, the vorticity is confined
to a thin region adjacent to the solid boundary. In this region, referred to as
the boundary layer, the flow is dominated by the effects of viscosity while
these effects fade away as one moves further away from the solid boundary.
Sufficiently far from the boundary, the effects of viscosity are negligible and
the equations governing the flow are identical to those for an Euler (ideal)
fluid and one solves the problem by melding together the solution for the Euler
fluid far from the boundary and the boundary layer approximation in a thin
layer adjacent to the boundary. The reason for developing such an approach is
the fact that in the ``boundary layer region'' an approximation is obtained for
the Navier-Stokes equations that is more amenable to analysis than the fully
non-linear equations. It is important to bear in mind that boundary layer
theory is an approximation of the Navier-Stokes equations in different parts of
the flow domain, that which is immediately adjacent to the solid boundary and
that which is away from the solid boundary. Great achievements in the field of
aerodynamics are a testimony to the efficacy and usefulness of such an
approximation with regard to solving analytically or computationally relevant
problems in a~particular geometrical setting. On the other hand, rigorous
analysis of the Prandtl  boundary layer equations is, despite significant
effort, far from being satisfactory (see \cite{Alexandre}, \cite{GerardVaret1,
GerardVaret2, GerardVaret3}, \cite{Kukavica1, Kukavica2}, \cite{LiuWangTang},
\cite{Lombardo}, \cite{Masmoudi1}, \cite{Oleinik2, Oleinik1}, \cite{SaCaf1,
SaCaf2}). An~alternative viewpoint for modeling the boundary layer phenomena
might thus bring some new insight on this issue.

The boundary layer approximation is not a constitutive approximation based on
different flow regimes though it seems to resemble such an approximation. That
is, one does not assume different constitutive assumptions for different
regions in the flow domain, based on some kinematical or other criterion, but
based on the value for the Reynolds number one merely carries out an
approximation of the Navier-Stokes equation in the flow domain. It is possible,
for instance, to assign different constitutive relations, based on the shear
rate, namely the fluid being an Euler fluid below a~certain shear rate and
a~Navier-Stokes or a~non-Newtonian fluid above the critical shear rate (such
a~classification is considered in Section~\ref{Sub1.5}), or as another
possibility a~non-Newtonian fluid if the shear rate is below a certain value
and a Navier-Stokes fluid above that shear rate, or any such assumption for the
constitutive response of the material, and to solve the corresponding equations
for the balance of linear momentum in the different flow domains. Such distinct
constitutive responses below and above a~certain kinematical criterion is akin
to models for the inelastic response of bodies wherein below a certain value
of the strain or stress, the body behaves as an elastic body while for values
above the critical value the body responds in an inelastic manner, which in turn
might lead to certain parts of a body to respond like an elastic body while
other parts could be exhibiting inelastic response. To make matters clear, in
a~solid cylindrical body that is undergoing torsion, a yield condition based on
the strain would lead to the body beyond a certain radius to respond
inelastically while below that threshold for the radius it responds as an
elastic body. In such an approach different constitutive relations are used in
different domains while in the classical boundary layer theory one uses
approximation of the equations of motion of a~particular fluid.

In this paper, we adopt the approach of assuming different constitutive
response relationships in different flow domains of the fluid, based on the
value of the shear rate or the value of the shear stress. We consider the
possibility that the character of the fluid changes when the certain
``activation'' criterion is reached. Here we consider an ``activation'' criterion
that is based on the shear rate or the shear stress, but it could be any other
criterion, say for instance the level of the electrical field in an
electrorheological fluid, or the temperature which changes the character of the
material from a fluid to a gas or a fluid to a solid, etc. Such an approach also
provides an alternative way to viewing the classical boundary layer
approximation in that it allows the fluid to behave like an Euler fluid in a
certain flow domain and a Navier-Stokes fluid elsewhere. Furthermore, based on other
criteria such as the Reynolds number we can carry out further approximations
with regard to governing equations in the different flow domains.

Within the context of the Navier-Stokes theory, boundary layers occur at flows
at sufficiently high Reynolds numbers. However, in the case of some
non-Newtonian fluids it is possible to have regions that are juxtaposed to a
solid boundary where the vorticity is concentrated even in the case of creeping
flow, i.e., flows wherein the inertial effect is neglected when the Reynolds
number is zero (see Mansutti and Rajagopal \cite{mansuttiKRR1991},
\cite{KRR1995}, \cite{harley-momoniat-rajagopal-2018}).
Thus, boundary layers are connected with the nonlinearities in
the governing equation and are not a consequence of just high Reynolds numbers.
Boundary layers can also occur at high Reynolds number in non-Newtonian fluids
of the differential type (see Mansutti et al. \cite{mansuttipontrelliKRR1993})
and of the integral type (see Rajagopal and Wineman \cite{KRRwineman1983}). It
is also possible that in non-Newtonian fluids one can have multiple decks with
dominance of different physical mechanisms in the different decks, and in these
different layers one can have the effects of viscosity, elasticity, etc.,
being significant, the delineation once again being determined within the context of a
specific governing equation (see Rajagopal et al.
\cite{KRRguptawineman1980,RajGuptaNa1983}).
On the other hand, we could have a more complicated situation wherein the flow
is characterized by different constitutive equations in different domains, and
in these different domains it might happen that one can further delineate
different subregions.

In the first part of this study, we provide a systematic classification of
the response of incompressible fluid-like materials ranging from the ideal Euler fluid
to non-Newtonian fluids that exhibit shear thinning/shear thickening,
stress thinning/stress thickening, as well as those responses where the
constitutive character of the material changes due to a threshold based on a
kinematical, thermal, stress or some other quantity (an example of the same is the
Bingham fluid which does not flow below a certain value of the shear stress
and starts to flow once the threshold is overcome) based on some criterion
concerning the level of shear rate or shear stress. We also provide a
systematic study of both activated and non-activated boundary conditions
ranging from \freeslip{} to \noslip. In carrying out our classification, we come
across the delineation of a class of fluids that, to our best knowledge, seems
to have not been studied by fluid dynamicists. This class of fluids is
characterized by the following intriguing dichotomy: (i) when the shear rate is below a
certain critical value the fluid behaves as the Euler fluid (i.e., there is no
effect of the viscosity, the shear stress vanishes), on the other hand (ii) if
the shear rate exceeds the critical value, dissipation starts to take place and fluid
can respond as a~shear (or stress) thinning or thickening fluid or as a
Navier-Stokes fluid. Implicit constitutive theory, cf. \cite{RajImpl05} and
also \cite{RajImpl03,rs08}, provides an elegant framework
to express such responses involving the activation criterion in a compact and
elegant manner that is also more suitable for further mathematical and
computational analysis.

In the second part of the paper, we study the mathematical properties of
three-dimensional internal flows in bounded smooth domains for fluids belonging
to this new class. We subject such flows to different types of boundary
conditions including \noslip, \navierslip, \freeslip{} and activated boundary
conditions like \noslipalt-\navierslipalt.
For this class of fluids and boundary conditions we prove the
global-in-time existence of a~weak solution in the sense of Leray to initial
and boundary value problems.

\section{Classification of incompressible fluids} \label{Ch1}
We are interested in studying isothermal flows
of \emph{incompressible fluids}, i.e., fluids which are assumed
to undergo only ischoric motions. For simplicity we consider
only fluids with uniform, constant density~$\rho_*>0$.
The balance of mass, the balance of linear momentum,
and the balance of angular momentum for such fluids
take, respectively, the following form
\begin{subequations}
  \label{eq:bal}
  \begin{align}
    \label{eq:balmass}
    \diver\bv &= 0,
    \\
    \label{eq:balmom}
    \rho_* \left( \pp{\bv}{t} + \diver(\bv\otimes\bv) \right)
      &= \diver\bT + \rho_*\bb,
    \\
    \label{eq:balang}
    \bT &= \bT^\top,
  \end{align}
\end{subequations}
where $\bv$ is the velocity field,
$\bT$ is the Cauchy stress,
and $\bb$ is a~given external body force.
Here, $\bm\otimes\bn$ denotes the second order tensor with the components
$\left(\bm\otimes\bn\right)_{ij} = m_i n_j$,
$\bA^\top$ denotes the transpose of a~tensor $\bA$,
and divergence of~$\bA$ is the vector with the
components $\left(\diver\bA\right)_i = \sum_{j=1}^3 \pp{\bA_{ij}}{x_j}$.

The incompressibility constraint~\eqref{eq:balmass}
can be written alternatively as
\begin{equation}
\tr\bD = \bD \mathbin{:} \bI = 0, \label{1.0a}
\end{equation} 
where $\bD$ (sometimes denoted $\bD\bv$) stands for the symmetric part of the velocity gradient, i.e., $\bD = \tfrac12(\nabla \bv + (\nabla \bv)^\top)$, $\bv$ being the velocity.

Due to this restriction, it is convenient to split the Cauchy stress tensor $\bT$ into its traceless (deviatoric) part $\bS$ and the mean normal stress, denoted by $m$ (more frequently expressed as $-p$), i.e.,
\begin{equation}
\bS = \bT - \frac13 (\tr \bT) \bI \quad \text{ and } \quad -p=m= \frac13 \tr\bT. \label{1.1}
\end{equation} 
Hence
\begin{equation}
\label{eq:cauchy}
\bT = m\bI + \bS = -p\bI + \bS,
\end{equation}
and, by virtue of \eqref{1.0a}, the stress power $\bT\mathbin{:} \bD$ satisfies
\begin{equation*}
\bT\mathbin{:}\bD  = \bS \mathbin{:} \bD.
\end{equation*} 

The main result of this section will be the classification of fluids using a simple framework that is characterized by a relation between $\bS$ and $\bD$, i.e., we are interested in materials whose response can be incorporated into the setting given by the \emph{implicit} constitutive equation
\begin{equation}
\cG(\bS,\bD) = \bO. \label{1.21}
\end{equation} 
The only restriction that we place is the requirement that response has to be \emph{monotone.}
For relevant discussion concerning non-monotone responses, we refer the reader to~\cite{malek_prusa_krr_2010}, \cite{leroux_krr_2013},
and \cite{JMPT2019}.

The incompressible Navier-Stokes fluid is a special sub-class of \eqref{1.21} where the relation between $\bS$ and $\bD$ is linear. This can be written either as
\begin{equation}
\bS = 2\nu_{*} \bD \quad \text{ with } \nu_{*}>0, \label{1.22}
\end{equation} 
where $\nu_{*}$ is called the (shear) \emph{viscosity}, or as
\begin{equation}
\bD = \alpha_{*} \bS \coloneqq \frac{1}{2\nu_{*}} \bS \quad \text{ with } \alpha_{*}>0, \label{1.22b}
\end{equation} 
where the coefficient $\alpha_{*}$ is called the \emph{fluidity}. Note that the stress power takes then the form
\begin{equation*}
\bS \mathbin{:} \bD = 2\nu_{*} |\bD|^2 = \alpha_{*} |\bS|^2.
\end{equation*} 

\begin{figure}
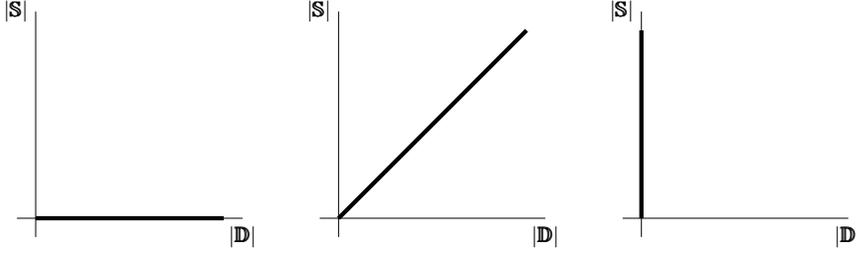

  \hfill\bigeuler\hfill\bigns\hfill\bigrigid\hfill\hfill
  \caption{From left to right, response of
           the Euler fluid~\eqref{1.24},
           the Navier-Stokes fluid~\eqref{1.22} (or~\eqref{1.22b}),
           and fluid allowing only motions fulfilling~\eqref{1.25}}
  \label{fig:euler-ns-rigid}
\end{figure}

There are two limiting cases when the stress power $\bS\mathbin{:}\bD$ vanishes. Either
\begin{equation}
\bS = \bO, \label{1.24}
\end{equation}
which implies that the fluid under consideration is the incompressible Euler fluid ($\bT = -p\bI$, see \eqref{1.1}), or
\begin{equation}
\bD = \bO \quad \text{ for all admissible flows.} \label{1.25}
\end{equation}
The latter corresponds to the situation where the body admits merely rigid body motions. More precisely, the flows fulfilling~\eqref{1.25} can be characterized through
\begin{equation*}
\bv (t,x) = \ba(t) \times x + \bb(t) \quad \text{ in all admissible flows.}
\end{equation*}
Response of models~\eqref{1.22} (or~\eqref{1.22b}), \eqref{1.24}, and \eqref{1.25} is shown in the Figure~\ref{fig:euler-ns-rigid}.

\subsection{Classical power-law fluids} \label{Sub1.1}
Classical power law fluids are described by
\begin{equation}
   \bS = 2\tilde\nu_{*} |\bD|^{r-2} \bD, \label{1.1.1}
\end{equation} 
which leads to
\begin{equation*}
\bS \mathbin{:} \bD = 2\tilde\nu_{*} |\bD|^{r}.
\end{equation*} 
Since $|\bD|^{r-2}\bD$ should have meaning for $|\bD|\to 0$, we require a~lower bound on $r$, namely
\begin{equation} 
r>1. \label{1.1.4}
\end{equation} 
Otherwise, if $r=1$ then $\lim_{|\bD|\to 0} \bD/|\bD|$ does not exist, and if $r<1$ then $|\bS|\to +\infty$ and stress concentration occurs at points where $\bD$ vanishes (i.e., $\bS$ plays the role of penalty for point where $\bD$ could vanish).

In what follows we study power-law fluids with the power-law index satisfying \eqref{1.1.4} and we shall investigate the responses of these fluids for $r\to 1$ and $r\to \infty.$ The latter corresponds to the case when the dual exponent $r'\coloneqq r/(r-1)$ tends to $1$.

We introduce the \emph{generalized viscosity} through
\begin{equation}
\nug(|\bD|) = \tilde\nu_{*} |\bD|^{r-2}.\label{1.1.5}
\end{equation}
In order to have the same units for $\nug$ as for the viscosity $\nu_{*}$ that appears in the formula for the Navier-Stokes fluid (see \eqref{1.22}), $\bD$ should scale as $d_{*}$ that has the unit \si{s^{-1}}. Thus, we replace
\eqref{1.1.5} by
\begin{equation*}
  \nug(|\bD|) = \nu_{*} \left(\frac{|\bD|}{d_{*}}\right)^{r-2}
  \quad\text{ where } [d_{*}]=\si{s^{-1}} \text{ and } [\nu_{*}]= \si{kg.m^{-1}.s^{-1}}
\end{equation*}
and we replace \eqref{1.1.1} by
\begin{equation}
  \bS = 2\nu_{*} \left(\frac{|\bD|}{d_{*}}\right) ^{r-2} \bD
  \quad\text{ with } [d_{*}]=\si{s^{-1}} \text{ and } [\nu_{*}]= \si{kg.m^{-1}.s^{-1}}.
  \label{1.1.7}
\end{equation} 
Of course, $\nu_{*}$ in \eqref{1.1.7} and $\nu_{*}$ in \eqref{1.22} are different in general; they however have the same units.

On considering \eqref{1.1.7}, we notice the following equivalence%
\footnote{%
  Indeed, starting for example from the formula on the left-hand side of~\eqref{1.1.8},
  we conclude that
  \begin{equation*}
    |\bS| = \frac{2\nu_*}{d_*^{r-2}} |\bD|^{r-1},
  \end{equation*}
  which implies
  \begin{equation*}
    |\bD|^{2-r} = \left( \frac{d_*^{r-2}}{2\nu_*} |\bS| \right)^\frac{2-r}{r-1}.
  \end{equation*}
  Hence
  \begin{equation*}
    \bD
    = \frac{1}{2\nu_*} \left( \frac{|\bD|}{d_*} \right)^{2-r} \bS
    = \left( \frac{1}{2\nu_*} \right)^{1+\frac{2-r}{r-1}} \left(d_*^{r-2}\right)^{\frac{2-r}{r-1}+1} |\bS|^\frac{2-r}{r-1} \bS
    = \frac{1}{2\nu_*} \left(\frac{|\bS|}{2\nu_*d_*} \right)^\frac{2-r}{r-1} \bS,
  \end{equation*}
  which leads to the formula on the right-hand side of~\eqref{1.1.8}.
}
\begin{equation}
\bS = 2\nu_{*} \left(\frac{|\bD|}{d_{*}}\right) ^{r-2} \bD  \quad\iff\quad \bD = \frac{1}{2\nu_{*}}
\left( \frac{|\bS|}{2\nu_*d_{*}}\right)^{\frac{2-r}{r-1}} \bS ,\label{1.1.8}
\end{equation} 
which gives rise to the following expressions for the generalized viscosity and \emph{generalized fluidity}
\begin{equation*}
\nug(|\bD|) = \nu_{*} \left(\frac{|\bD|}{d_{*}}\right) ^{r-2} \quad \text{ and } \quad \alg(|\bS|) = \frac{1}{2\nu_{*}}
\left( \frac{|\bS|}{2\nu_*d_{*}}\right)^{\frac{2-r}{r-1}}.
\end{equation*} 
It also allows us to express the stress power in the form ($r' \coloneqq r/(r-1)$)
\begin{equation*}
\begin{split}
\bS \mathbin{:} \bD &= \left(\frac{1}{r} + \frac{1}{r'}\right) \bS \mathbin{:} \bD = \frac{1}{r} \bS \mathbin{:} \bD + \frac{1}{r'} \bS \mathbin{:} \bD \\ 
& = 2\nu_*d_*^2 \left( \frac{1}{r} \left(\frac{|\bD|}{d_*}\right)^r + \frac{1}{r'} \left(\frac{|\bS|}{2\nu_*d_*}\right)^{r'} \right).
\end{split}
\end{equation*} 
Summarizing,
\begin{equation}
\bS = 2 \nug(|\bD|^2) \bD = 2\nu_* \left( \frac{|\bD|}{d_*}\right)^{r-2} \bD \,\iff \, \bD = \alg(|\bS|^2) \bS = \frac{1}{2\nu_*}\left( \frac{|\bS|}{2\nu_*d_*} \right)^{r'-2} \bS,
\label{1.1.11}
\end{equation}
emphasizing that the equivalence in \eqref{1.1.11} holds only if $r\in (1,+\infty)$ (which is equivalent to $r' \in (1, +\infty)$).

Generalizing the approach used in \cite{mrr95}, we will investigate the limits of $\bS$ and $\nug$ as $\bD$ tends to zero or infinity,  or vice versa, study limits of $\bD$ and $\alg$ as $\bS$ vanishes or tends to infinity.

Letting $|\bD|\to 0{+}$ we obtain, starting from the formula on the left-hand side of~\eqref{1.1.11},
\begin{equation*}
\begin{aligned} 
|\bS| &\to 0 &\text{if } r>1, & \\
|\bS| &\le 2\nu_* d_* &\text{if } r=1, &
  \qquad \text{(as } |\bD|\to 0{+} \text{)} \\
|\bS| &\to +\infty &\text{if } r<1, & \\
\end{aligned}
\end{equation*}
and
\begin{equation} \label{1.1.12a}
\begin{aligned} 
|\nug(|\bD|)| &\to 0 &\text{if } r>2, & \\
\nug(|\bD|) &= \nu_* &\text{if } r=2, &
  \qquad \text{(as } |\bD|\to 0{+} \text{)} \\
|\nug(|\bD|)| &\to +\infty &\text{if } r<2. & \\
\end{aligned}
\end{equation}
Thus, we note that the Cauchy stress $\bT$ in the fluid tends to a purely spherical stress when $r$ is greater than $1$ and the norm of $\bD$ tends to $0+$, or put differently, the constitutive relation for the fluid reduces to that for an Euler fluid.

Similarly, letting $|\bD|\to +\infty$ we have
\begin{equation*}
\begin{aligned} 
|\bS| &\to +\infty &\text{if } r>1, & \\
|\bS| &\le 2\nu_* d_* &\text{if } r=1, &
  \qquad \text{(as } |\bD|\to +\infty \text{)} \\
|\bS| &\to 0 &\text{if } r<1, & \\
\end{aligned}
\end{equation*}
and
\begin{equation*}
\begin{aligned} 
|\nug(|\bD|)| &\to +\infty &\text{if } r>2, & \\
\nug(|\bD|) &= \nu_* &\text{if } r=2, &
  \qquad \text{(as } |\bD|\to +\infty \text{)} \\
|\nug(|\bD|)| &\to 0 &\text{if } r<2. & \\
\end{aligned}
\end{equation*}
In order to investigate the behavior of $\bD$ and fluidity in the limiting case, it is useful to employ the expression on the right-hand side of \eqref{1.1.11}. Thus, for $|\bS| \to 0{+}$, we get
\begin{equation*}
\begin{aligned} 
|\bD| &\to 0 &\text{if } r'>1, & \\
|\bD| &\le d_* &\text{if } r'=1, &
  \qquad \text{(as } |\bS|\to 0{+} \text{)} \\
|\bD| &\to +\infty &\text{if } r'<1, & \\
\end{aligned}
\end{equation*}
and
\begin{equation} \label{1.1.16}
\begin{aligned} 
|\alg(|\bS|)| &\to 0 &\text{if } r'>2, & \\
\alg(|\bS|) &= \frac{1}{2\nu_*} &\text{if } r'=2, &
  \qquad \text{(as } |\bS|\to 0{+} \text{)} \\
|\alg(|\bS|)| &\to +\infty &\text{if } r'<2. & \\
\end{aligned}
\end{equation}
Similarly, letting $|\bS| \to +\infty$ we have
\begin{equation*}
\begin{aligned} 
|\bD| &\to +\infty &\text{if } r'>1, & \\
|\bD| &\le d_* &\text{if } r'=1, &
  \qquad \text{(as } |\bS|\to +\infty \text{)} \\
|\bD| &\to 0 &\text{if } r'<1, & \\
\end{aligned}
\end{equation*}
and
\begin{equation*}
\begin{aligned} 
|\alg(|\bS|)| &\to +\infty &\text{if } r'>2, & \\
\alg(|\bS|) &= \frac{1}{2\nu_*} &\text{if } r'=2, &
  \qquad \text{(as } |\bS|\to +\infty \text{)} \\
|\alg(|\bS|)| &\to 0 &\text{if } r'<2. & \\
\end{aligned}
\end{equation*}

\begin{figure}
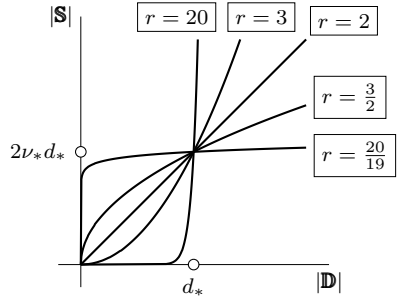

  \centering\bigpowerlaw
  \caption{Response of the power-law model~\eqref{1.1.7} for various values of $r$}
  \label{fig:powerlaw}
\end{figure}

Next, we study the response of the classical power-law fluid with regard to its dependence on the value of power-law index ($r\to 1+$ and $r' \to 1+$). Figure \ref{fig:powerlaw} illustrates behavior for both large $r$ and $r'$ approaching $1$.

Letting $r\to 1+$ in \eqref{1.1.11}, we observe that, for $\bD\neq \bO$,
$\bS = 2\nu_* d_* \frac{\bD}{|\bD|}$ holds and thus $|\bS|\le 2\nu_* d_*$. On
the other hand for $\bD = \bO$ and for any $\bS\in \RR^{3\times
3}_\mathrm{sym}$ such that $|\bS|\le 2\nu_* d_*$ we can find sequences
$\{\bD_n\}_{n=1}^{\infty} \subset\RR^{3\times 3}_\mathrm{sym}$ and
$\{r_n\}_{n=1}^{\infty} \subset(1,\infty)$ such that
\begin{align*}
  \bD_n &\rightarrow \bD=\bO,
  &
  r_n &\rightarrow 1,
  &
  2\nu_* d_* |\bD_n|^{r_n-2}\, \bD_n &\rightarrow \bS
  &
  \text{as } n &\rightarrow \infty.
\end{align*}
Consequently, \eqref{1.1.11} for $r\to 1+$ approximates the response
that could be referred to as
\emph{rigid/free-flow like behavior}:
\begin{equation} \label{1.1.17}
\begin{aligned}
&\text{if } |\bD| \neq 0 &&\text{ then } \bS = 2\nu_*d_* \frac{\bD}{|\bD|}, \\
&\text{if } |\bD| = 0 &&\text{ then } |\bS| \le 2\nu_*d_*.
\end{aligned}
\end{equation}
Instead of viewing \eqref{1.1.17} as multivalued response (both in the variables $\bD$ and $\bS$), it is possible to write \eqref{1.1.17} as a continuous graph over the Cartesian product $\RR^{3\times 3} \times \RR^{3\times 3}$ (see the framework \eqref{1.21}) defined through a~(scalar) equation
\begin{equation}
(|\bS| - 2\nu_* d_*)^{+} + \bigl| 2\nu_* d_* \bD - |\bD| \bS \bigr| = 0. \label{1.1.18}
\end{equation}

For determining the behavior of \eqref{1.1.11} as $r\to +\infty$ we prefer to study \eqref{1.1.11}$_2$ for $r' \to 1+$ and analogous to the above consideration
we observe that
\begin{equation} \label{1.1.19}
\begin{aligned} 
&\text{if } |\bS| \neq 0 &&\text{ then } \bD = \frac{d_*}{2\nu_*} \frac{\bS}{|\bS|}, \\
&\text{if } |\bS| = 0 &&\text{ then } |\bD| \le \frac{d_*}{2\nu_*}.
\end{aligned}
\end{equation}
We can call this response \emph{Euler/rigid like response}.
We can again rewrite \eqref{1.1.19} as
\begin{equation}
(2\nu_* |\bD| - d_*)^{+} + \bigl| 2\nu_* |\bS| \bD - d_* \bS \bigr| = 0. \label{1.1.20}
\end{equation}

The models~\eqref{1.1.17} and \eqref{1.1.19} are examples of fluids described within the context of an~activation criterion. More examples will be discussed in Subsections~\ref{Sub1.4} and~\ref{Sub1.5}. The slash formalism {\tt name1}/{\tt name2} (that we will use also below) means that material behaves as {\tt name1} before activation and as {\tt name2} after the activation criterion is met.
 
\subsection{Generalized power-law fluids and stress power-law fluids} \label{Sub1.2}
The formula~\eqref{1.1.11} suggests the introduction of \emph{generalized
power-law fluids} and \emph{generalized stress power-law fluids}
by requiring that, for the former,
\[
  \bS = 2\nug \left( |\bD|^2 \right) \bD
\]
and, for the latter,
\[
  \bD = \alg \left( |\bS|^2 \right) \bS
\]
where $\nug$ and $\alg$ are non-negative continuous functions referred to as
\emph{the generalized viscosity} and \emph{the generalized fluidity}.
The quantities $|\bD|^2 = \tr\bD^2$ and $|\bS|^2 = \tr\bS^2$ representing
the second invariants of $\bD$ and $\bS$, respectively, can be
viewed as natural higher dimensional generalizations of the shear-rate
and the shear stress, respectively.

We further introduce \emph{zero shear rate viscosity} as
\[
  \nu_0 \coloneqq \lim_{|\bD| \to 0} \nug\left( |\bD|^2 \right)
\]
and \emph{zero shear stress fluidity} through
\[
  \alpha_0 \coloneqq \lim_{|\bS| \to 0} \alg\left( |\bS|^2 \right).
\]

It follows from~\eqref{1.1.12a} that for classical power-law fluids
the zero shear rate viscosity vanishes for $r>2$, is finite for $r=2$,
and blows up if $r\in(1,2).$ A~similar behavior can be inferred from~\eqref{1.1.16}
for the zero shear stress fluidity: for $r'>2$ (i.e., for $r\in(1,2)$)
$\alpha_0$ is zero, for $r=2$ it is positive, and for $r'\in(1,2)$
(it means for $r>2$) the generalized fluidity becomes singular
in the vicinity of the origin.

Such behavior is not experimentally observed in any fluid; more frequently
both the zero shear rate viscosity and zero shear stress fluidity are
finite. The most popular generalizations of the classical power-law fluid
that exhibit these features as $|\bD| \to 0{+}$ (resp. $|\bS| \to 0{+}$)
take the form
\begin{equation} \label{A1}
  \bS = 2\nu_* \left( \frac12 + \frac12\frac{|\bD|^2}{d_*^2} \right)^\frac{r-2}{2} \bD
\end{equation}
and
\begin{equation} \label{A2}
  \bD = \alpha_* \left( \frac12 + \frac12\frac{|\bS|^2}{(2\nu_*d_*)^2} \right)^\frac{r'-2}{2} \bS.
\end{equation}

We refer the reader to~\cite{malek_prusa_krr_2010} for further details;
here we emphasize two observations. First, although both~\eqref{A1} and
\eqref{A2} are invertible for $r\in(1,+\infty)$, \eqref{A2} is not an~inverse
of~\eqref{A1} and vice versa (compare it to~\eqref{1.1.11}). Second, both
formulas are defined for all $r\in(-\infty,+\infty)$ and
$r'\in(-\infty,+\infty)$ respectively. The relationship $r'=\frac{r}{r-1}$
is however understood as a relation valid for $r>1.$

If $r<1$, then $\bS$ considered as a function of $\bD$ given by~\eqref{A1}
is not monotone;
the same holds for~\eqref{A2} if $r'<1.$ We refer the reader
to~\cite{malek_prusa_krr_2010} for more details and to~\cite{leroux_krr_2013}
for further nontrivial extensions.
 
\subsection{Fluids that can be viewed as a~mixture of power-law fluids} \label{Sub1.3}
It is natural to consider the possibility that the total response of the
fluid-like material is given as the sum of particular responses (a~simplified
scenario for the mixtures) of the individual (in our case two) contributors,
i.e.,
\begin{equation} \label{B1}
  \bS = \bS_1 + \bS_2.
\end{equation}

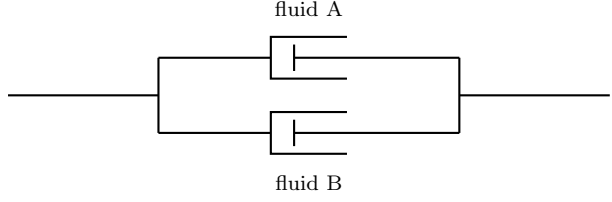
\begin{figure}
  \footnotesize
  \centering
  \begin{tikzpicture}
    \draw [thick]  (0, 0) -- (2, 0);
    \draw [thick]  (2,-0.5) -- (2, 0.5);
    \draw [damper] (2,+0.5) -- (6,+0.5) node[midway,above=1.5em] {fluid A};
    \draw [damper] (2,-0.5) -- (6,-0.5) node[midway,below=1.5em] {fluid B};
    \draw [thick]  (6,-0.5) -- (6, 0.5);
    \draw [thick]  (6, 0) -- (8, 0);
  \end{tikzpicture}
  \caption{1D mechanical analogue of an additive stress model~\eqref{B1}}
  \label{fig:two_dashpots}
\end{figure}

\begin{figure}
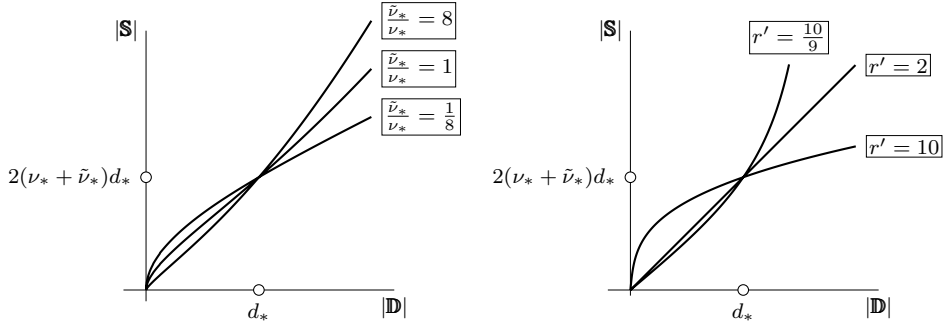

  \hfill\bigaddstressi\hfill\bigaddstressii\hfill\hfill
  \caption{Response of models represented by~\eqref{B1}.
           On the left, the case~\eqref{Si} with $r=\frac32$ and $q=\frac52.$
           On the right, the case~\eqref{Sii} with $\nu_*=\tilde\nu_*.$}
  \label{fig:addstress}
\end{figure}

One may think of putting two dashpots into a parallel arrangement; the
response of one dashpot captures behavior of a fluid A, and the response
of another dashpot corresponds to a fluid B; see Figure~\ref{fig:two_dashpots}.
To illustrate the potential of this setting, we consider for illustration
three examples:
\begin{enumerate}[(i)]
  \item \label{Si}
    $\bS_1 = 2\nu_* \left( \frac{|\bD|}{d_*} \right)^{r-2} \bD$
    and
    $\bS_2 = 2\tilde\nu_* \left( \frac12 + \frac12\frac{|\bD|^2}{d_*^2} \right)^\frac{q-2}{2} \bD$
    with
    $r\in(1,2)$,
    $q>2$;
  \item \label{Sii}
    $\bS_1 = 2\nu_* \bD$ and $\bS_2$ fulfills
    $\bD = \frac{1}{2\tilde\nu_*}
      \left( \frac12 + \frac12\frac{|\bS_2|^2}{(2\tilde\nu_*d_*)^2} \right)^\frac{r'-2}{2} \bS_2$
    with $r'\in(1,+\infty)\setminus\{2\}$;
  \item \label{Siii}
    $\bS_1$ responds as in~\eqref{1.1.18} and $\bS_2$ responds as
    in~\eqref{1.1.20}.
\end{enumerate}
The responses of fluids modeled by the constitutive
expressions~\eqref{Si} and \eqref{Sii} are shown
in Figure~\ref{fig:addstress} for some choice of parameters.
Further examples will be provided in Subsection~\ref{Sub1.5}.

\subsection{Fluids with bounded shear rate or bounded shear stress} \label{Sub1.4}
We have seen in Subsection~\ref{Sub1.1} that the models~\eqref{1.1.17}
and \eqref{1.1.19} exhibit an interesting feature, namely, the stress~$\bS$ is bounded as the symmetric part of the velocity gradient is varied, and vice versa.
While there are several mathematical advantages to the stress $\bS$ expressed as a~function of $\bD$, as this would greatly simplify the number of equations that one needs to consider when one substitutes it into the balance of linear momentum (even an explicit expression of the symmetric part of the velocity gradient as a function of $\bS$ might lead to a simplified structure for the equations) experiments on colloidal fluids clearly show that a fully implicit theory is necessary (see for example \cite{Boltenhagen,Holmes} and further references in \cite{prusaperlacova}).
When considering the models~\eqref{A1} and \eqref{A2} with $r=1$ and $r'=1$ (and adjusting a scaling by factor $\sqrt{2}$), we obtain
\begin{equation} \label{C1}
  \bS = 2\nu_* \frac{\bD}{\sqrt{1+\frac{|\bD|^2}{d_*^2}}}
\end{equation}
and
\begin{equation} \label{C2}
  \bD = \alpha_* \frac{\bS}{\sqrt{1+\frac{|\bS|^2}{(2\nu_*d_*)^2}}},
\end{equation}
and we observe that, in the case of~\eqref{C1},
\begin{equation*}
  |\bS| \leq 2\nu_*d_* \quad \text{ for all } \bD,
\end{equation*}
and, in the case of~\eqref{C2},
\begin{equation*}
  |\bD| \leq d_* \quad \text{ for all } \bS.
\end{equation*}

It is convenient to generalize~\eqref{C1} and \eqref{C2} in the following
manner: for parameters $a,b \in(0,+\infty)$ consider
\begin{equation} \label{C3}
  \bS = 2\nu_* \frac{\bD}{\left(1+\left(\frac{|\bD|}{d_*}\right)^a\right)^\frac1a}
\end{equation}
and
\begin{equation} \label{C4}
  \bD = \alpha_* \frac{\bS}{\left(1+\left(\frac{|\bS|}{2\nu_*d_*}\right)^b\right)^\frac1b}.
\end{equation}

\begin{figure}
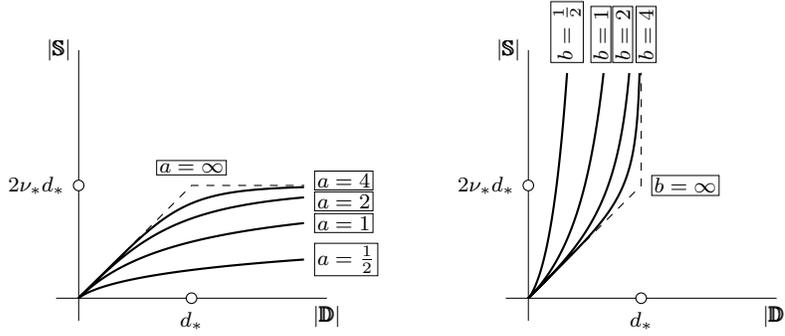

  \hfill\biglimitingstress\hfill\biglimitingstrain\hfill\hfill
  \caption{Response of stress-limiting model~\eqref{C3} on the left and
           shear-limiting model~\eqref{C4} on the right}
  \label{fig:limiting}
\end{figure}

In both cases, it is worth studying the behavior of the fluids for large
$a$ and $b$. When $a\to +\infty$ in~\eqref{C3}, the constitutive relation
approximates the response of the activated fluid which behaves as the Navier-Stokes
fluid prior to the activation and the magnitude of the stress remains bounded;
analogously, when $b\to +\infty$ in~\eqref{C4}, the constitutive relation
approximates the response such that the magnitude of $\bD$ remains bounded;
see Figure~\ref{fig:limiting}.

\subsection{Activated fluids} \label{Sub1.5}
In this section we study two classes of fluids: the first class is activation
based on the value of the stress (similar in character to a~Bingham
fluid) while the second class is activation based on the value of
the shear rate.

\begin{figure}
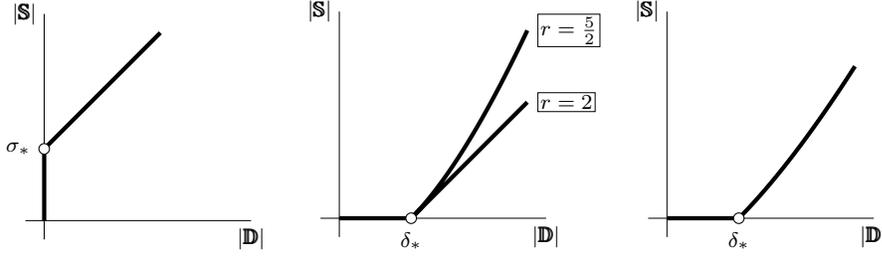

  \hfill\bigbingham\hfill\bigeulerns\hfill\bigeulerlad\hfill\hfill
  \caption{Response of the Bingham fluid (on the left),
           activated Euler/Navier-Stokes fluid (in the middle, $r=2$),
           Euler/power-law fluid (in the middle, $r=\frac52$, $d_*=\delta_*$), and
           Euler/Ladyzhenskaya fluid (on the right, $r=\frac52$, $d_*=\delta_*$, $\nu_*=\tilde\nu_*$)}
  \label{fig:activated}
\end{figure}

The first class of fluids that are studied flow only if the generalized shear stress
$|\bS| = \left( \tr\bS^2 \right)^\frac12$ exceeds a~certain critical
value $\sigma_*,$ referred to as \emph{the yield stress.} Once the fluid flows,  we
assume the fluid behavior is described by the constitutive expression for a~generalized power-law or
a~generalized stress power-law fluid. In the parts of the subdomain where
$|\bS|$ is below $\sigma_*$ the fluid can only
translate or rotate as a rigid body. Such responses are traditionally
described (see~\cite{DuvantLions}) through the dichotomy
\begin{equation} \label{D1}
\begin{aligned}
  |\bS| &\leq \sigma_* &&\iff& \bD &= \bO, \\
  |\bS| &>    \sigma_* &&\iff& \bS &= \sigma_* \frac{\bD}{|\bD|} + \bS_2
  \quad\text{ with }\left\{\begin{array}{lll}
    \text{either } & \bS_2 = 2\nug\left( |\bD|^2 \right) \bD, \\
    \text{or }     & \bD   = \alg\left( |\bS_2|^2 \right) \bS_2.
  \end{array}\right.
\end{aligned}
\end{equation}
In the case of the stress  $\bS_2 = 2\nu_*\bD$ we obtain the constitutive representation for the \emph{Bingham fluid} (see
Figure~\ref{fig:activated} on the left) and if $\bS_2 = 2\nug\left( |\bD|^2
\right) \bD$ then we obtain the constitutive representation for the \emph{Herschel-Bulkley} fluid. It is
worth of observing that~\eqref{D1} can be equivalently written within the context of 
the framework for implicit constitutive equations~\eqref{1.21}. Specifically, considering \eqref{D1} with the expression 
$\bS_2 = 2\nug\left( |\bD|^2 \right) \bD$
the equivalent formulation can be expressed as
\begin{equation*}
  2\nug\left( |\bD|^2 \right) \bD
  = \frac{\left( |\bS| - \sigma_* \right)^+}{|\bS|} \bS
\end{equation*}
where $(t)^+=\max\{t,0\}$ for $t\in\RR$.

On the other hand, considering~\eqref{D1} with the expression $\bD = \alg\left( |\bS_2|^2 \right) \bS_2$, the equivalent representation reads 
\begin{equation*}
  \bD =  \alg\left( \left| \bS - \sigma_* \frac{\bD}{|\bD|}\right|^2 \right) \frac{\left( |\bS| - \sigma_* \right)^+}{|\bS|} \bS,
\end{equation*}
which is in our opinion worthy of detailed investigation.

The next class is a~dual to \eqref{D1} in the following sense. If the generalized shear rate
$|\bD|$ is below a critical value $\delta_*,$ the flow is frictionless.
Inner friction between the fluid layers becomes important only when $|\bD|$
exceeds $\delta_*.$ Then the fluid can flow as a~Navier-Stokes fluid, or
a~power-law fluid (see Figure~\ref{fig:activated} on the right), or a~generalized
power-law fluid, or a~generalized stress power-law fluid. To summarize,
analogous to~\eqref{D1}, we can describe such a~response through the relation
\begin{equation*}
\begin{aligned}
  |\bD| &\leq \delta_* &&\iff& \bS &= \bO, \\
  |\bD| &>    \delta_* &&\iff& \bD &= \delta_* \frac{\bS}{|\bS|} + \alg\left( |\bS|^2 \right) \bS.
\end{aligned}
\end{equation*}
It is not surprising that this relation can be written as an explicit relation for the stress $\bS$
in terms of the symmetric part of the velocity gradient $\bD$, namely
\begin{equation} \label{D5}
  \alg\left( |\bS|^2 \right) \bS = \frac{\left( |\bD|-\delta_* \right)^+}{|\bD|} \bD.
\end{equation}
If there is $\nug$ such that
\begin{equation} \label{D6}
  \bD = \alg\left( |\bS|^2 \right) \bS
  \quad\iff\quad
  \bS = 2\nug\left( |\bD|^2 \right) \bD
\end{equation}
then~\eqref{D5} can be written in the form (\eqref{D6} describes the behavior
of the fluid after activation)
\begin{equation*}
  \bS = 2\nug\left( |\bD|^2 \right) \frac{\left( |\bD|-\delta_* \right)^+}{|\bD|} \bD.
\end{equation*}
Explicit equivalence holds for the Navier-Stokes
fluid (see~\eqref{1.22b}) and for the standard power-law fluids
(see~\eqref{1.1.11}). Then we obtain the response
\begin{equation} \label{D7}
  \bS = 2\nu_* \frac{\left( |\bD|-\delta_* \right)^+}{|\bD|} \bD
\end{equation}
and
\begin{equation} \label{D8}
  \bS = 2\nu_* \left(\frac{|\bD|}{d_*}\right)^{r-2} \frac{\left( |\bD|-\delta_* \right)^+}{|\bD|} \bD,
\end{equation}
\begingroup\setlength\emergencystretch{\hsize}\hbadness=10000
respectively (see Figure~\ref{fig:activated} on the right).
We call the response~\eqref{D7} \emph{the Euler/Navier-Stokes fluid}
and the response~\eqref{D8} \emph{the Euler/power-law fluid}. If
\begin{equation*}
  \bS = \left( 2\nu_* + 2\tilde\nu_* \left(\frac{|\bD|}{d_*}\right)^{r-2} \right)
        \frac{\left( |\bD|-\delta_* \right)^+}{|\bD|} \bD,
\end{equation*}
\endgroup
with $r>2$ we call this
response \emph{Euler/Ladyzhenskaya fluid}, since O.~A.~Ladyzhenskaya was the first
to consider the generalization of the Navier-Stokes constitutive equation to the form
\begin{equation} \label{D9}
  \bS = \left( 2\nu_* + 2\tilde\nu_* |\bD|^{r-2} \right) \bD
\end{equation}
and showed that unsteady internal flows of such
fluids in a bounded smooth container admit \emph{unique} weak solutions
if $r>\frac52$ (or $\frac{d+2}{2}$ in general dimension $d$);
see~\cite{ladyzhenskaya_1967,ladyzhenskaya_1968,ladyzhenskaya_1969}.
Improvement of the uniqueness result to the range $r\geq\frac{11}{5}$
is established in~\cite{bulicek-ettwein-kaplicky-prazak-2010,bulicek-kaplicky-prazak-2018}.

\bigskip
\paragraph{Simple shear flows of the Euler/Navier-Stokes fluid}

\begin{figure}
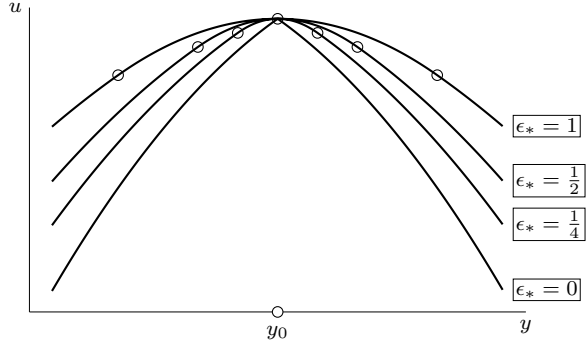

  \centering\simpleflows
  \caption{Simple shear flows of the regularized Euler/Navier-Stokes fluid for
           various values of the added viscosity $\epsilon_*\nu_*$ and fixed
           $C>0$. Circles mark the point of activation
           $y=y_0\pm\frac{\sqrt{2}\delta_*\epsilon_*}{2|C|}$
           where $|u'|=\sqrt{2}\delta_*$. The degenerate
           case $\epsilon_*=0$ is activated everywhere, i.e.,
           $|u'|\geq\sqrt{2}\delta_*$.
           Note that the velocity profiles are determined up to
           an~additive constant (because no boundary conditions are enforced);
           here we take $u_0=0$.}
  \label{fig:simpleflows}
\end{figure}

For the sake of illustration of response of the fluid~\eqref{D7}
we consider a~simple shear flow of such a~fluid. In order to characterize
such fluids it is also useful to consider a~modified model which we call
\emph{the regularized Euler/Navier-Stokes fluid} given by
\begin{equation} \label{D10}
  \bS = 2\nu_* \left( \epsilon_* + \frac{\left( |\bD|-\delta_* \right)^+}{|\bD|} \right) \bD
\end{equation}
with an~extra parameter $\epsilon_*\geq0$. We will consider solutions of
the balance equations~\eqref{eq:bal}, \eqref{eq:cauchy} for
the response of the~fluids described by~\eqref{D10}
(and the degenerate case~\eqref{D7}) in
$\RR^2$, with the velocity taking the form
\begin{equation*}
  \bv(x,y) = \left( u(y), 0 \right) \qquad x,y\in\RR
\end{equation*}
for some $u:\RR\rightarrow\RR$ absolutely continuous
on every compact interval in $\RR$. We note that now
we deal with solutions of the governing equations in $\RR^2$
so there are no boundary conditions involved. It is easy to check
that if $\epsilon_*>0$ then all such solutions
of~\eqref{eq:bal}, \eqref{eq:cauchy}, and~\eqref{D10} fulfill
\begin{subequations}
\label{D1213}
\begin{gather}
  \left( \epsilon_* + \mathcal{H}(|u'|-\sqrt{2}\delta_*) \right) u'' = -2 C
  \quad\text{ a.e. in }\RR, \\
  p(x) = -2\nu_* C x + p_0,
\end{gather}
\end{subequations}
with some $C\in\RR$, $\mathcal{H}(t)=1$ if $t>0$
and $\mathcal{H}(t)=0$ otherwise.  The formulas~\eqref{D1213}
represent a~generalization of the well-known equations for simple shear
flows of the Navier-Stokes fluid, which is a~special case with
$\delta_*=0$.
In fact all the solutions of~\eqref{D1213} take the form
\begin{subequations}
\begin{align}
  \label{D14}
  u(y) &=
  \left\{\begin{array}{ll}
    \! \! \!
    -\frac{C}{\epsilon_*} (y\!-\!y_0)^2
    +u_0
    &
    |y\!-\!y_0|\!\leq\!\frac{\sqrt{2}\delta_*\epsilon_*}{2|C|},
    \\
    \! \! \!
    -\frac{C}{1+\epsilon_*}\!
    \left(\! (y\!-\!y_0)^2
           \!+\! \sqrt{2}\delta_*\!\left|\frac{y-y_0}{C}\right|
           \!-\! \epsilon_*\!\left(\frac{\sqrt{2}\delta_*}{2C}\right)^{\!2}
    \right)
    +u_0
    &
    |y\!-\!y_0|\!\geq\!\frac{\sqrt{2}\delta_*\epsilon_*}{2|C|},
  \end{array}\right.
  \! \! \! \! \! \!
  \\
  p(x) &= -2\nu_* C x + p_0,
\end{align}
\end{subequations}
with any $C,y_0,u_0,p_0\in\RR$. In the interval
$\bigl\{y\in\RR,\, |y-y_0|\leq\frac{\sqrt{2}\delta_*\epsilon_*}{2|C|} \bigr\}$
the fluid is in the regime below the ``activation'' threshold (where
$|u'|\leq\sqrt{2}\delta_*$) with the viscosity $\epsilon_*\nu_*$ while
outside this interval the threshold is exceeded and the generalized viscosity has the value
$\nu_*\bigl(\epsilon_*+1-\frac{\sqrt{2}\delta_*}{|u'|}\bigr)$.
Taking the limit $\epsilon_*\rightarrow0+$ one obtains
\begin{subequations}
\label{D16D17}
\begin{align}
  \label{D16}
  u(y) &= -C \left( (y\!-\!y_0)^2 + \sqrt{2}\delta_*\!\left|\tfrac{y-y_0}{C}\right| \right)
  +u_0, \\
  p(x) &= -2\nu_* C x + p_0,
\end{align}
\end{subequations}
which indeed is a~solution of the balance equations for the Euler/Navier-Stokes
fluid \eqref{D7} with any $C,y_0,u_0,p_0\in\RR$. Now the flow exceeds
the activation threshold everywhere except of $y=y_0$ where
$u'(y_0\pm) = \mp\sqrt{2}\delta_*\frac{C}{|C|}$ and the shear rate jumps there
by virtue of the vanishing viscosity. In Figure~\ref{fig:simpleflows} we display
a~family of solutions~\eqref{D14} for varying $\epsilon_*$ and \eqref{D16}
(matching $\epsilon_*=0$) for fixed values of $C,y_0$.

Apart of the family of solutions~\eqref{D16D17} the Euler/Navier-Stokes
fluid~\eqref{D7} admits also the simple shear flows which do not exceed the threshold
$|u'|=\sqrt{2}\delta_*$, the viscosity is zero and therefore any admissible velocity
profile is a~solution. Such solutions are characterized by
\begin{subequations}
\label{D18D19}
\begin{gather}
  |u'|\leq\sqrt{2}\delta_*
  \quad\text{ a.e. in }\RR, \\
  p = p_0,
\end{gather}
\end{subequations}
with some locally Lipschitz continuous $u:\RR\rightarrow\RR$
and $p_0\in\RR$. In fact the families \eqref{D16D17} and
\eqref{D18D19} are all possible weak solutions for a~simple shear flow
of the Euler/Navier-Stokes fluid.
 
\subsection{Classification of incompressible fluids} \label{Sub1.6}
The previous exposition should indicate the broad spectrum of fluid
responses that can be described within the setting
\begin{equation} \label{E1}
  b\left( |\bD|^2 \right) \bD = a\left( |\bS|^2 \right) \bS,
\end{equation}
where $a$ and $b$ are continuous (not necessarily always differentiable)
functions.

\begin{table}
\centering
\caption{Summary of systematic classification of fluid-like response
         with the corresponding $|\bD|$-$|\bS|$ diagrams}
\renewcommand{\arraystretch}{1.5}
\setlength\arrayrulewidth{.1pt}
\footnotesize
\begin{tabular}{p{0.195\linewidth}r|p{0.195\linewidth}r|p{0.195\linewidth}r}
  Euler/\allowbreak{}rigid
&
  \eulerlim
&
  \mbox{Navier-Stokes}/\allowbreak{}limiting \mbox{shear-rate}
&
  \shearlim
&
  fluid body allowed to move only rigidly
&
  \rigid
\\\hline
  Euler/\allowbreak\mbox{shear-thickening}\!
&
  \eulerthick
&
  \mbox{shear-thickening}\!
&
  \thick
&
  rigid/\allowbreak\mbox{shear-thickening}\!
&
  \rigidthick
\\\hline
  Euler/\allowbreak\mbox{Navier-Stokes}
&
  \eulerlin
&
  Navier-Stokes
&
  \ns
&
  rigid/\allowbreak\mbox{Navier-Stokes}
&
  \rigidlin
\\\hline
  Euler/\allowbreak\mbox{shear-thinning}
&
  \eulerthin
&
  \mbox{shear-thinning}
&
  \thin
&
  rigid/\allowbreak\mbox{shear-thinning}
&
  \rigidthin
\\\hline
  Euler
&
  \euler
&
  limiting
&
  \stresslim
&
  rigid/\allowbreak\mbox{free-flow}
&
  \rigidlim
\\\hline
  \multicolumn{2}{c|}{
    $|\bD| \le \delta_* \iff \bS = \bO$
  }
&
  \multicolumn{2}{c|}{
    no activation
  }
&
  \multicolumn{2}{c}{
    $|\bS| \le \sigma_* \iff \bD = \bO$
  }
\end{tabular}
\label{tab_summary}
\end{table}
 
The following Table~\ref{tab_summary} summarizes these observations in
a~different way, paying attention to the broad range of models covered by~\eqref{E1}.
It includes the Euler (frictionless) fluid at one extreme and a fluid that only moves rigidly
at the other extreme and contains the responses ranging from the fluids
enforcing the activation criterion $|\bD|\le\delta_* \iff
\bS=\bO$ through non-activated fluids to the fluids that are governed by the activation
criterion $|\bS|\le\sigma_* \iff \bD=\bO.$ Vertically, the range of~$r$
is iterated from top to bottom:
$r=+\infty$, $r\in(2,+\infty)$, $r=2$, $r\in(1,2)$, and $r=1.$ Thus, at the
bottom left corner we have perfectly frictionless Euler fluid and at the top
right corner we have a fluid that can only undergo rigid motions. In the middle of the
table the Navier-Stokes model is placed.
 
\subsection{Activated boundary conditions} \label{Sub1.7}
Boundary conditions can have as much impact on the nature of the flow as
the constitutive equation for the fluid in the bulk. We illustrate it explicitly
in this subsection. Here, for the sake of
clarity, we restrict ourselves to internal flows, i.e.,
we consider (sufficiently smooth) domain $\Omega\subset\RR^3$ and
we assume that
\begin{equation} \label{F1}
  \bv\cdot\bn = 0 \quad \text{on }\partial\Omega,
\end{equation}
where $\bn: \partial\Omega\to\RR^d$ denotes the mapping that assigns the outward
unit normal vector to any $\bx\in\partial\Omega.$ The behavior of the
fluid at the tangential direction near the boundary is described by the
equations that reflect mutual interaction between the solid boundary and the fluid
flowing adjacent to the boundary. These are constitutive equations which we shall concern ourselves with.
In order to specify them, in the spirit of previous parts of the paper, we first recall
some energy estimates.

We consider the balance equations~\eqref{eq:bal} together with~\eqref{eq:cauchy}.
For simplicity, we assume that $\bb=\b0$.
Forming the scalar product of the equation~\eqref{eq:balmom} with $\bv$ and
integrating the result over~$\Omega,$ we arrive at
\begin{equation*}
  \dd{}{t} \int_\Omega \rho_*\frac{|\bv|^2}{2} \d\bx
    + \int_\Omega \diver\left( \rho_*\frac{|\bv|^2}{2} \bv \right) \d\bx
    = \int_\Omega \diver\left( \bS\bv \right) \d\bx
    - \int_\Omega \bS\mathbin{:}\bD\, \d\bx
    - \int_\Omega \diver\left( p\bv \right) \d\bx,
\end{equation*}
where we have used \eqref{eq:balmass} twice. Gauss's theorem and the
requirement~\eqref{F1} then lead to
\begin{equation} \label{F4}
  \dd{}{t} \int_\Omega \rho_*\frac{|\bv|^2}{2} \d\bx
    + \int_\Omega \bS\mathbin{:}\bD\, \d\bx
    + \int_{\partial\Omega} \left(-\bS\right)\mathbin{:}\left(\bv\otimes\bn\right) \d S
    = 0.
\end{equation}
By virtue of the symmetry of $\bS$, see~\eqref{eq:balang} and~\eqref{eq:cauchy},
we obtain
\begin{equation*}
  \left(-\bS\right) \mathbin{:} \left(\bv\otimes\bn\right)
  = \left(-\bS\right) \mathbin{:} \left(\bn\otimes\bv\right)
  = \left(-\bS\bn\right)_{\btau} \cdot \bv_{\btau},
\end{equation*}
where $\bz_{\btau}$ denotes the projection of $\bz: \partial\Omega \to \RR^d$
to the plane tangent to $\partial\Omega$ (at the point of
$\partial\Omega$ under consideration). Finally, introducing the notation
\begin{equation*}
  \bs \coloneqq \left(-\bS\bn\right)_{\btau}
  \quad\text{(projection of the normal traction to the tangent plane)},
\end{equation*}
we can rewrite~\eqref{F4} in the form
\begin{equation*}
  \dd{}{t} \int_\Omega \rho_*\frac{|\bv|^2}{2} \d\bx
    + \int_\Omega \bS\mathbin{:}\bD \d\bx
    + \int_{\partial\Omega} \bs\cdot\bv_{\btau} \d S
    = 0.
\end{equation*}

The discussion in Section~\ref{Ch1} thus far has been focused on discussing models
within the context of the framework $\cG\left(\bS,\bD\right)=\bO$ (see~\eqref{1.21}).
In fact, the discussion concerned the restricted class of models of the
form
\begin{equation} \label{F7}
  a\left( |\bD|^2 \right) \bD = b\left( |\bS|^2 \right) \bS,
\end{equation}
where $a$ and $b$ were non-negative continuous (not necessarily everywhere
differentiable) functions.

In a manner similar to the class of models defined through~\eqref{F7}, we could develop analogously the
identical framework of relations linking $\bs$ and $\bv_{\btau}$, i.e., to consider
various classes of boundary conditions that fit the form
\begin{equation*}
  \bh\left( \bs, \bv_{\btau} \right) = \b0,
\end{equation*}
where we deal with a (monotone) continuous function $\bh:\RR^d\times\RR^d\to\RR^d,$
or a more restrictive class
\begin{equation} \label{F9}
  \tilde a\left( |\bv_{\btau}|^2 \right) \bv_{\btau} = \tilde b\left( |\bs|^2 \right) \bs,
\end{equation}
where $\tilde a,\tilde b:[0,+\infty)\to[0,+\infty)$ are non-negative
continuous functions.

We shall not consider the problem within the context of such generality for the following two
reasons: (i) we do not have enough experimental data that would support
nonlinear relations between $\bs$ and $\bv_{\btau}$, and (ii) the extension
of the framework developed to take into account such nonlinear relations is
straightforward and follows in a~manner similar to that discussed
in Subsections~\ref{Sub1.1}--\ref{Sub1.5}.

In what follows, we restrict ourselves to both activated and non-activated
responses that are (after activation) linear. In a manner similar to that in the introductory
part of Section~\ref{Ch1} where we considered the product $\bS\mathbin{:}\bD,$ we notice here that the
product $\bs\cdot\bv_{\btau}$ vanishes if either
\begin{equation} \label{F10}
  \bs = \b0 \quad\text{on }\partial\Omega,
\end{equation}
or
\begin{equation} \label{F11}
  \bv_{\btau} = \b0 \quad\text{on }\partial\Omega\text{ for all admissible flows}.
\end{equation}
The condition~\eqref{F10} is referred to as the \emph{\freeslip} condition
condition, expressing the fact that the boundary
exhibits no friction to the flow, in the sense that the shear stress vanishes,
and fluid flows tangentially to the boundary.
On the other hand, the condition~\eqref{F11}, referred to as the \emph{\noslip}
boundary condition, requires that flows adhere to the boundary. It has the
character of a boundary constraint.

A linear relation between $\bs$ and $\bv_{\btau}$ is known as the \emph{\navierslip}:
\begin{equation} \label{F12}
  \bs = \gamma_* \bv_{\btau} \quad
  \text{ with }\gamma_*>0.
\end{equation}

We consider two types of activated boundary conditions. First, the relation
\begin{equation} \label{F13}
\begin{aligned}
  |\bs| &\leq s_* &&\iff& \bv_{\btau} &= \b0, \\
  |\bs| &>    s_* &&\iff& \bs      &= s_* \frac{\bv_{\btau}}{|\bv_{\btau}|} + \gamma_*\bv_{\btau},
\end{aligned}
\end{equation}
which has been coined as the {\noslipalt/\navierslipalt} condition
in the literature, but which we will
refer to as the \emph{\noslip/\navierslip} condition for consistency;
$s_*$ is the yield stress that is positive. It can
be rewritten into the form~\eqref{F9} through its equivalent characterization
\begin{equation} \label{F14}
  \gamma_* \bv_{\btau} = \frac{\left( |\bs|-s_* \right)^+}{|\bs|} \bs.
\end{equation}

In analogy, the second type of activated condition is given through the
description
\begin{equation} \label{F15}
\begin{aligned}
  |\bv_{\btau}| &\leq v_* &&\iff& \bs      &= \b0, \\
  |\bv_{\btau}| &>    v_* &&\iff& \bv_{\btau} &= v_* \frac{\bs}{|\bs|} + \frac{1}{\gamma_*}\bs,
\end{aligned}
\end{equation}
where $v_*>0.$ The equivalent description of~\eqref{F15}, that can be referred
to as the \emph{\freeslip/\navierslip}
condition, takes the form
\begin{equation} \label{F16}
  \bs = \gamma_* \frac{\left( |\bv_{\btau}|-v_* \right)^+}{|\bv_{\btau}|} \bv_{\btau}.
\end{equation}
These responses are summarized in Table~\ref{tab_summary_bc}.

\begin{table}
\centering
\caption{Classification of boundary activation of fluid response
         with the corresponding $|\bv_{\btau}|$-$|\bs|$ diagrams.
         The last row reflects the usage of the term \emph{slip}
         in Section~\ref{Ch2} to describe a~broad class of boundary
         conditions of the slip type.}
\renewcommand{\arraystretch}{1.5}
\setlength\arrayrulewidth{.1pt}
\footnotesize
\begin{tabular}{p{0.145\linewidth}|p{0.195\linewidth}|p{0.145\linewidth}|p{0.195\linewidth}|p{0.145\linewidth}}
  \mbox\freeslip
&
  \mbox\freeslip/\allowbreak\mbox\navierslip
&
  \mbox\navierslip
&
  \mbox\noslip/\allowbreak\mbox\navierslip
&
  \mbox\noslip
\\
  \euler
&
  \eulerlin
&
  \ns
&
  \rigidlin
&
  \rigid
\\[2em]
  $\bs = \b0$
&
  $|\bv_{\btau}| \leq v_* \iff \bs = \b0$
&
  $\bs \sim \bv_{\btau}$
&
  $|\bs| \le s_* \iff \bv_{\btau} = \b0$
&
  $\bv_{\btau} = 0$
\\
  \eqref{F10}
&
  \eqref{F15} or \eqref{F16}
&
  \eqref{F12}
&
  \eqref{F13} or \eqref{F14}
&
  \eqref{F11}
\\\hline
  \multicolumn{4}{c|}{\slip} & \noslip
\end{tabular}
\label{tab_summary_bc}
\end{table}

\bigskip
\paragraph{Simple shear flows of the Navier-Stokes fluid and the Euler/Navier-Stokes
           fluid subject to activated boundary conditions}

Finally, in order to emphasize the role of boundary conditions in
determining the nature of the flow, we consider
Poiseuille flow between two parallel plates located at $y=\pm L$.
All types of boundary conditions listed in Table~\ref{tab_summary_bc} will be
considered. We can write boundary
conditions~\eqref{F14} and \eqref{F16} together as
\begin{equation} \label{F32}
  \gamma_* \frac{\left( |\bv_{\btau}|-v_* \right)^+}{|\bv_{\btau}|} \bv_{\btau}
  = \frac{\left( |\bs|-s_* \right)^+}{|\bs|} \bs
  \quad\text{ on }\partial(\RR\times(-L,L))
\end{equation}
requiring that at least one of $v_*$ and $s_*$ is zero. Let us consider a~simple
shear flow of the Euler/Navier-Stokes fluid~\eqref{D16D17} in domain
$\RR\times(-L,L)$ for given $L>0$. As shown
in Section~\ref{Sub1.5} simple shear flows of the Euler/Navier-Stokes
fluid are in the form~\eqref{D16D17} or \eqref{D18D19}.
Let us assume symmetry $y_0=0$ in~\eqref{D16}; it will be obvious later
that the converse is not possible.
Normalizing~\eqref{D16} to a given flow rate $Q\in\RR$ such that
\begin{equation} \label{F18}
  \int_{-L}^L u(u) \d y = Q
\end{equation}
we obtain
\begin{subequations}
\label{F3031}
\begin{align}
  u &= - C \left( y^2 + \sqrt{2}\delta_* \left|\frac{y}{C}\right| \right)
       + \frac{Q}{2L}
       + C \left( \frac{L^2}{3} + \frac{\sqrt{2}\delta_*L}{2|C|} \right)
  ,\\
  p &= -2\nu_* C x + p_0
\end{align}
\end{subequations}
being defined for any $C\in\RR\setminus\{0\}$. It requires
a~trivial, but tedious, computation to check that simple shear
flow of Euler/Navier-Stokes fluid~\eqref{F3031}
solves the balance equations in $\RR\times(-L,L)$ together
with boundary condition~\eqref{F32} on $\{y=\pm L\}$ provided that
\begin{align}
  \label{F33}
  C &=
  \frac{3Q}{4L^3}
  \Bigg[
    \left( 1 - \frac{\sqrt{2}\delta_*L^2 + 2v_*L}{|Q|} \right)^+
    \!\!
    - \frac{3\nu_*}{3\nu_*\!+\!\gamma_*L}
      \left( 1 - \frac{\sqrt{2}\delta_*L^2 + 2v_*L +\frac{2s_*L^2}{3\nu_*}}{|Q|} \right)^+
  \Bigg]
  \!\!
\end{align}
and $p_0\in\RR$ is arbitrary. If $C$ given by formula~\eqref{F33} is
zero, then all flows which fulfill \eqref{D18D19}, \eqref{F32},
and \eqref{F18} are solutions; more precisely if $C=0$ all the solutions
are given by
\begin{subequations}
\label{F34353637}
\begin{align}
  \label{F34}
  |u'|&\leq\sqrt{2}\delta_*
  \quad\text{ a.e. in }\RR, \\
  \gamma_*|u|&\leq\gamma_*v_*
  \quad\text{ a.e. on }\{|y|=L\}, \\
  \label{F36}
  \int_{-L}^L u \,\d y &= Q, \\
  \label{F37}
  p &= p_0,
\end{align}
\end{subequations}
with some Lipschitz continuous $u:[-L,L]\rightarrow\RR$ and
$p_0\in\RR$. Family \eqref{F3031}, \eqref{F33} and family
\eqref{F34353637} represent in fact all
possible simple shear flow solutions of motions of the Euler/Navier-Stokes fluid
subject to \noslip/\navierslip{} or \freeslip/\navierslip{} boundary
conditions~\eqref{F32}. We summarize combinations of bulk and boundary
activation criterions in Table~\ref{tab_simple_flows}.

\begin{table}
\centering
\caption{Solutions for simple shear flows of the Euler/Navier-Stokes fluid
         in combination with different activated and classical
         boundary conditions. The middle column contains $|\bv_{\btau}|$-$|\bs|$
         diagrams of boundary response (on the left) and $|\bD|$-$|\bS|$
         diagrams of bulk response (on the right). The solid segments and the
         circles (colored red in the electronic version) mark the
         part of the response being attained in the specific case.
         Note that \noslip/\navierslip{} (contrary to \freeslip/\navierslip)
         admits a~mode with the activation threshold exceeded in the bulk
         with the boundary under activation threshold.
         Also note that the \freeslip{} condition
         admits only Euler mode, frictionless solutions.
         For the Navier-Stokes limit just let $\delta_*\coloneqq0$.}
\renewcommand{\arraystretch}{1.8}
\setlength\arrayrulewidth{.1pt}
\footnotesize
\vspace{-6pt}
\begin{tabular}{p{0.25\linewidth}|cc|p{0.45\linewidth}}
\noalign{\vskip 0mm}
\multicolumn{4}{l}{ \textbf{\freeslip/\navierslip:} \eqref{F32} \& $s_*=0$ }
\\ \hline
$|Q|\leq\sqrt{2}\delta_*L^2+2v_*L$
&
\eulerlinone
&
\eulerlinone
&
\eqref{F34353637}
\\\hline
$|Q|\geq\sqrt{2}\delta_*L^2+2v_*L$
&
\eulerlintwo
&
\eulerlintwo
&
\eqref{F3031}, $C=$
\newline
\RaggedLeft
$
  \frac{\gamma_*L}{3\nu_*+\gamma_*L}
  \frac{3\left(|Q|-\sqrt{2}\delta_*L^2-2v_*L\right)}{4L^3}
  \frac{Q}{|Q|}
$
\\\hline
\multicolumn{1}{c}{} &
\multicolumn{1}{c}{\footnotesize bdry} &
\multicolumn{1}{c}{\footnotesize bulk} &
\multicolumn{1}{c}{} \\
\noalign{\vskip 0mm}
\multicolumn{4}{l}{ \textbf{\noslip/\navierslip:} \eqref{F32} \& $v_*=0$ }
\\\hline
$|Q|\leq\sqrt{2}\delta_*L^2$
&
\rigidlinrest
&
\eulerlinone
&
\eqref{F34353637}
\\\hline
$\sqrt{2}\delta_*L^2\leq|Q|\leq\sqrt{2}\delta_*L^2+\frac{2s_*L^2}{3\nu_*}$
&
\rigidlinone
&
\eulerlinthree
&
\eqref{F3031},
$
  C=
  \frac{3\left(|Q|-\sqrt{2}\delta_*L^2\right)}{4L^3}
  \frac{Q}{|Q|}
$
\\\hline
$|Q|\geq\sqrt{2}\delta_*L^2+\frac{2s_*L^2}{3\nu_*}$
&
\rigidlintwo
&
\eulerlinfour
&
\eqref{F3031}, $C=$
\newline
\RaggedLeft
$
  {\bigl[
  \frac{\gamma_*\!L}{3\nu_*\!+\gamma_*\!L}
  \frac{3\left(|Q|-\sqrt{2}\delta_*\!L^2\right)}{4L^3}
  \!+\!
  \frac{3\nu_*}{3\nu_*\!+\gamma_*\!L}
  \frac{s_*}{2\nu_*\!L}
  \bigr]}\!
  \frac{Q}{|Q|}
$\!\!\!
\\\hline
\multicolumn{1}{c}{} &
\multicolumn{1}{c}{\footnotesize bdry} &
\multicolumn{1}{c}{\footnotesize bulk} &
\multicolumn{1}{c}{} \\
\noalign{\vskip 0mm}
\multicolumn{4}{l}{ \textbf{\freeslip:} \eqref{F32} \& $s_*=0 \;\&\; v_*\rightarrow\infty$ }
\\
\noalign{\vskip -2.5mm}
\multicolumn{4}{l}{ \textbf{\freeslip:} \eqref{F32} \& $s_*=0 \;\&\; \gamma_*=0$ }
\\\hline
$Q\in\RR$
&
\eulermove
&
\eulerlinone
&
\eqref{F34353637}
\\\hline
\multicolumn{1}{c}{} &
\multicolumn{1}{c}{\footnotesize bdry} &
\multicolumn{1}{c}{\footnotesize bulk} &
\multicolumn{1}{c}{} \\
\noalign{\vskip 0mm}
\multicolumn{4}{l}{ \textbf{\noslip:} \eqref{F32} \& $v_*=0 \;\&\; s_*\rightarrow\infty$ }
\\
\noalign{\vskip -2.5mm}
\multicolumn{4}{l}{ \textbf{\noslip:} \eqref{F32} \& $v_*=0 \;\&\; \gamma_*\rightarrow\infty$ }
\\\hline
$|Q|\leq\sqrt{2}\delta_*L^2$
&
\rigidrest
&
\eulerlinone
&
\eqref{F34}, $u(\pm L)=0$, \eqref{F36}, \eqref{F37}
\\\hline
$|Q|\geq\sqrt{2}\delta_*L^2$
&
\rigidmove
&
\eulerlintwo
&
\eqref{F3031},
$
  C=
  \frac{3\left(|Q|-\sqrt{2}\delta_*L^2\right)}{4L^3}
  \frac{Q}{|Q|}
$
\\\hline
\multicolumn{1}{c}{} &
\multicolumn{1}{c}{\footnotesize bdry} &
\multicolumn{1}{c}{\footnotesize bulk} &
\multicolumn{1}{c}{} \\
\noalign{\vskip 0mm}
\multicolumn{4}{l}{ \textbf{\navierslip:} \eqref{F32} \& $v_*=0 \;\&\; s_*=0$ }
\\\hline
$|Q|\leq\sqrt{2}\delta_*L^2$
&
\nsrest
&
\eulerlinone
&
\eqref{F34353637}
\\\hline
$|Q|\geq\sqrt{2}\delta_*L^2$
&
\nsmove
&
\eulerlintwo
&
\eqref{F3031},
$
  C=
  \frac{\gamma_*L}{3\nu_*+\gamma_*L}
  \frac{3\left(|Q|-\sqrt{2}\delta_*L^2\right)}{4L^3}
  \frac{Q}{|Q|}
$
\\\hline
\multicolumn{1}{c}{} &
\multicolumn{1}{c}{\footnotesize bdry} &
\multicolumn{1}{c}{\footnotesize bulk} &
\multicolumn{1}{c}{} \\
\end{tabular}
\label{tab_simple_flows}
\end{table}

\section{Mathematical analysis of flows
         of activated Euler fluids} \label{Ch2}
Long-time and large-data existence theory (within the context of weak solutions)
for a~broad class of fluids described by implicit constitutive relation~\eqref{1.21} has been
developed in~\cite{BGMRS2012,BGMS2012}. These works deal with internal flows of
incompressible fluids with monotone responses, asymptotically behaving as
$|\bS|=\mathcal{O}\left(|\bD|^{r-1}\right)$ as $|\bD|\rightarrow\infty$ or
$|\bD|=\mathcal{O}\left(|\bS|^{\frac{r}{r-1}-1}\right)$ as $|\bS|\rightarrow\infty$
with $\frac65<r<\infty$ (or $\frac{2d}{d+2}<r<\infty$ in general dimension $d$)\footnote{%
  In fact, in~\cite{BGMRS2012,BGMS2012} the results are established even in
  a~more general setting replacing the Lebesgue spaces by the Orlicz spaces.
}. In order to get a~pressure\footnote{%
  Subtle difference between \emph{thermodynamic pressure}, the \emph{mean
  normal stress} (the latter usually referred to as pressure in mathematical
  literature on incompressible fluids), and the Lagrange multiplier is not to
  be discussed in this paper and we refer interested reader
  to~\cite{M_KRR_2010,KRR_pressure}. Henceforth we refer to the Lagrange
  multiplier that enforces the incompressibility condition as the ``pressure''.
} which is integrable over the space-time cylinder
in the unsteady case, the theory is developed with a~boundary condition
allowing some kind of \slip. The overview of the problem concerning the
connection between the integrability of the pressure and a~specific boundary
condition is given in~\cite{Frehse2003}; see also the original studies
\cite{koch-solonnikov-2001,koch-solonnikov-2002}
within the context of the Navier-Stokes equations with variable viscosity,
and in~\cite{bulicek-malek-rajagopal-2007} in the context of incompressible
fluids with pressure-dependent viscosity.
Existence theory when the \navierslip{} boundary condition is enforced has recently been
extended to the \noslipalt/\navierslipalt{}
boundary condition in~\cite{bulcek-malek-a:2016,bulcek-malek-b:2016}.
The theory for unsteady flows subject to the \noslip{} boundary condition
for standard power-law fluids
(generating \emph{strictly} monotone operators) was originally developed in \cite{DiRuWo06}
using the concept of local pressures introduced earlier in \cite{wolf07}. The theory
for the \noslip{} boundary condition can also be found in a~more recent study~\cite{BDS13},
where the divergence-free Lipschitz approximations of divergence-free Bochner-Sobolev functions are
constructed and analyzed (from the point of view of dependence of their
mathematical properties on approximation parameters). With such constructions,
the analysis of the problem can be performed without introducing the notion
of pressure.
On the other hand, in various fluid problems (heat-conducting fluids, fluids with pressure
dependent viscosity, etc.) one cannot avoid the pressure from the mathematical formulations. The reason the pressure
has, within the context of weak solutions, better properties for slip-type boundary conditions
(including \noslipalt/\navierslipalt{})
than for the \noslip{} is due to the different function space settings: for \noslip{} one works with Sobolev spaces vanishing on
the boundary while for slip-type boundary conditions one only requires that
the normal component of the velocity vanishes. The latter setting is compactible with the Helmholtz decomposition, which provides
the advantage.
It is also worth emphasizing that, in the case of the \noslip{} boundary condition,
there is also a remarkable difference in introducing the presssure for steady flow problems and
for evolutionary problems. While for the time-independent problems, one can identify the pressure using
a variant of de Rham's theorem, one cannot use the same tool directly for time-dependent problems
since the time derivative is not, in general, a distribution. We will pay attention to contruction
of the pressure for various boundary conditions in what follows.

In this section we will provide an~existence theory for steady and
unsteady flows of activated Euler fluids considering various types of
behavior after activation and various types of boundary conditions.
More specifically, we will study the system\footnote{%
  We assume that density $\rho$ is constant and we replace
  $\frac{p}{\rho}$ merely by $p$ throughout the whole section.
  Note that such $p$, although customarily called ``pressure''
  in mathematical fluid dynamics literature, is actually the
  mean normal stress scaled by the (constant) density.
}
\begin{subequations}
\label{eq:sys}
\begin{align}
  \label{eq:sysdiv}
  \diver \bv &= 0
  & &\text{ in }(0,T)\times\Omega ,\\
  \pp{\bv}{t} + \diver(\bv\otimes\bv) - \diver\bS &= -\nabla p + \bb
  & &\text{ in }(0,T)\times\Omega ,\\
  \label{eq:sysstress}
  \bS &= \smash{2\nu_* \left( |\bD|-\delta_* \right)^+ \mathcal{S}(|\bD|) \tfrac{\bD}{|\bD|}}
  & &\text{ in }(0,T)\times\Omega ,\\
  \label{eq:sysbcnor}
  \bv\cdot\bn &= 0
  & &\text{ on }(0,T)\times\partial\Omega ,\\
  \label{eq:sysbctan}
  \bh(\bs,\bv_{\btau}) &= \b0
  & &\text{ on }(0,T)\times\partial\Omega ,\\
  \label{eq:sysic}
  \bv(0,\cdot) &= \bv_0
  & &\text{ in }\Omega.
\end{align}
\end{subequations}
Here $\mathcal{S}:[0,\infty)\rightarrow[0,\infty)$ is supposed to be of
the following forms:
either
\begin{align*}
  \mathcal{S}\equiv1
\end{align*}
giving the Euler/Navier-Stokes fluid~\eqref{D7}, or
\begin{align*}
  \mathcal{S}(d)=\bigl(\tfrac{d}{d_*}\bigr)^{r-2}
  \quad\text{ or }\quad
  \mathcal{S}(d)=\Bigl(A+\bigl(\tfrac{d}{d_*}\bigr)^2\Bigr)^\frac{r-2}{2},
  \quad A>0,
\end{align*}
leading to the Euler/power-law fluid~\eqref{D8}, or
\begin{align*}
  \mathcal{S}(d)=1+A\bigl(\tfrac{d}{d_*}\bigr)^{r-2},
  \quad r>2, \quad A>0,
\end{align*}
leading to the Euler/Ladyzhenskaya fluid.

It is not difficult to verify (see Appendix~\ref{app:examples}) that the graph
$\cG\subset\RR^{3\times3}_\mathrm{sym}\times\RR^{3\times3}_\mathrm{sym}$
defined through
\begin{align*}
  (\bS,\bD)\in\cG \text{ if and only if } \bS\text{ and }\bD\text{ fulfill }\eqref{eq:sysstress}
\end{align*}
is a~maximal monotone $r$-graph, i.e., $\cG$ has the following properties:
\begin{enumerate}[($\cG$1)]
  \item \label{Gzero}
    $(\bO,\bO)\in\cG$;
  \item \label{Gmonotone}
    $(\bS_1-\bS_2)\mathbin{:}(\bD_1-\bD_2)\geq0$ for all
    $(\bS_1,\bD_1)\in\cG$ and $(\bS_2,\bD_2)\in\cG$;
  \item \label{Gmaximal}
    if $\bS,\bD\in\RR^{3\times3}_\mathrm{sym}$ satisfy
    $(\bS-\tilde\bS)\mathbin{:}(\bD-\tilde\bD)\geq0$ for all
    $(\tilde\bS,\tilde\bD)\in\cG$,
    then $(\bS,\bD)\in\cG$;
  \item \label{Gcoercive}
    there exist $r\in(1,\infty)$, $\alpha,\beta\in(0,\infty)$ such that
    $\bS\mathbin{:}\bD\geq\alpha\bigl(|\bS|^{r'}+|\bD|^r\bigr)-\beta$ whenever
    $(\bS,\bD)\in\cG$ and $\frac1r+\frac{1}{r'}=1$.
\end{enumerate}

Suitable choices of the function $\bh(\bs,\bv_{\btau})$ cover the boundary
conditions \eqref{F10}, \eqref{F11}, \eqref{F12}, \eqref{F14},
\eqref{F16}. Analogous to the above setting, we require that the graph
$\cB\subset\RR^3\times\RR^3$ defined through
\begin{align*}
  (\bs,\bv)\in\cB \Leftrightarrow \bs\text{ and }\bv\text{ fulfill }\eqref{eq:sysbctan}
\end{align*}
is a~maximal monotone $2$-graph, i.e., $\cB$ has the following properties:
\begin{enumerate}[($\cB$1)]
  \item \label{Bzero}
    $(\b0,\b0)\in\cB$;
  \item \label{Bmonotone}
    $(\bs_1-\bs_2)\cdot(\bv_1-\bv_2)\geq0$ for all
    $(\bs_1,\bv_1)\in\cB$ and $(\bs_2,\bv_2)\in\cB$;
  \item \label{Bmaximal}
    if $\bs,\bv\in\RR^3$ satisfy
    $(\bs-\tilde\bs)\cdot(\bv-\tilde\bv)\geq0$ for all
    $(\tilde\bs,\tilde\bv)\in\cB$,
    then $(\bs,\bv)\in\cB$;
  \item \label{Bcoercive}
    there are $\tilde\alpha,\tilde\beta\in(0,\infty)$ such that
    $\bs\cdot\bv\geq\tilde\alpha\left(|\bs|^2+|\bv|^2\right)-\tilde\beta$ for all
    $(\bs,\bv)\in\cB$.
\end{enumerate}
The requirement ($\cB$\ref{Bcoercive}) can be easily verified for
the boundary conditions~\eqref{F12}, \eqref{F14} and~\eqref{F16}.

Note that there is no boundary term in the weak formulation of the problem
in the case of the \freeslip{} condition~\eqref{F10}; this condition does not
invalidate the analysis. On the other hand the \noslip{} boundary
condition~\eqref{F11} needs to be treated separately.

We wish to remark that we restrict ourselves to $2$-graphs here
merely for simplicity. We could include a~more general
boundary condition involving $q$-graphs with $1<q<\infty$ or even
graphs which describe non-monotone relations $\bh(\bs,\bv_{\btau})=\b0$.
This extension would be possible due to compactness of~$\bv^n_{\btau}$
in $L^1((0,T)\times\partial\Omega)$ (in the sense of traces of Sobolev
functions). The reason we prefer to consider only $2$-graphs is twofold.
First, we wish to present different slipping mechanisms in their simplest
forms and, second, we did not want to increase the number of model parameters
in order to make the description of the results accessible to the reader.
See, however, \cite{bulicek-malek-2019-kolmogorov}
for a~general treatment of complex boundary conditions.

Apart from the general purpose of this paper we are further motivated
to study the problem~\eqref{eq:sys} for the following
reasons.
\begin{enumerate}
  \item
    The most studied systems of PDEs (partial differential equations) in fluid
    mechanics are the Euler
    equations (when $\bS=\bO$, or $\delta_*\rightarrow\infty$ in~\eqref{eq:sysstress})
    and the Navier-Stokes equations (when $\bS=2\nu_*\bD$, or $\delta_*=0$
    and $\mathcal{S}\equiv1$ in~\eqref{eq:sysstress}).
    The system of PDEs considered here is placed
    between them, as $\delta_*\in(0,\infty)$. While
    \eqref{eq:sysdiv}--\eqref{eq:sysstress} can, particularly for $\delta_*$ large,
    share several features associated with the physics of the Euler fluid
    (or the Euler equations), we will document that the mathematical properties of the
    flows described by \eqref{eq:sys} are similar to
    those described by the Navier-Stokes equations. This is important as
    recent achievements in the mathematical theory of the Euler equations
    considered in a~reasonable physical setting show that the equations
    exhibit pathological solutions within the framework of weak solutions with
    bounded (kinetic) energy (see
    \cite{delellis-szekelyhidi-2010,wiedemann-2011}).

    Fluids described by~\eqref{eq:sysstress} seem to have been completely overlooked both
    in physics and mathematical fluid dynamics literature; this may well be due to the
    fact that such behavior has not been observed. Below, we will
    focus on filling this lacuna and on developing the mathematical foundations
    associated with the problem \eqref{eq:sys}.
  \item
    It is worth noticing that the activated Euler fluids characterized
    by~\eqref{eq:sysstress} represent the models dual to the Bingham fluids
    that are obtained by interchanging the role of $\bD$ and $\bS$
    in~\eqref{eq:sysstress}. A~mathematical theory for Bingham fluids,
    in the spirit of the theory developed here, is given
    in~\cite{bulcek-malek-b:2016,bulcek-malek-a:2016,maringova-zabensky-2018},
    where the reader can also find more references concerning the earlier results
    on the analysis of flows of the Bingham fluids and their generalizations.
  \item
    The set-up of the problem considered here will be also used to show how different types
    of boundary conditions can be treated (while restricting ourselves to
    internal flows). We will also focus on the relation between the
    considered boundary conditions and the properties of the mean normal stress $p$.
  \item
    Since the operator $-\diver\bS$ is elliptic and
    degenerates for $|\bS|\leq \delta_*$, the theory presented below can
    be viewed as an~approach for studying degenerate problems.
  \item
    Finally, the constitutive relation~\eqref{eq:sysstress} is regularized by
    \begin{align*}
      \bS^{\epsilon}(\bD) = \left(
        \epsilon |\bD|^{q-2} +
        2\nu_* \frac{\left(|\bD|-\delta_*\right)^+}{|\bD|} \mathcal{S}(|\bD|)
      \right)\bD
    \end{align*}
    with $\epsilon>0$ and $q\geq2$ large enough.
    This explicit regularization allows us to proceed
    explicitly in the subsequent analysis.
\end{enumerate}

\bigskip
\subsection{Function spaces}

In what follows we assume that $\Omega\subset\RR^3$ is a~domain,
i.e., a~bounded open connected set.
For $1\leq p\leq\infty$, $k\in\mathbb{N}$,
$L^p(\Omega)$ and $W^{k,p}(\Omega)$ denote
the standard Lebesgue and Sobolev space respectively,
i.e., spaces of measurable functions of finite norm
  \begin{alignat*}{2}
    \|f\|_{L^p(\Omega)}
    &= \|f\|_{p,\Omega}
    &&= \Bigl( \int_\Omega |f|^p \Bigr)^{\frac1p}, \\
    \|f\|_{W^{k,p}(\Omega)}
    &= \|f\|_{k,p,\Omega}
    &&= \sum_{j=0}^{k} \||\nabla^j f|\|_{p,\Omega}
    \quad |\nabla^j f|
    = \Bigl( \sum_{|\alpha|=j} |\mathrm{D}^\alpha f|^2 \Bigr)^\frac12
  \end{alignat*}
respectively. When there is no risk of confusion the
subscript $\Omega$ can be omitted.
Often we will use symbol $L^p(\Omega)^{3\times 3}_\mathrm{sym}$
to denote functions with values in symmetric tensors with
$L^p$-integrable components. Bold-face symbols $\bW^{k,p}$
and $\bL^p$ will denote vector-valued Sobolev and Lebesgue
functions respectively. Parentheses $(\cdot,\cdot)_\Omega$ will
denote duality pairing in $L^p(\Omega)$ and $L^\frac{p}{p-1}(\Omega)$
including vector and tensor-valued case; the subscript $\Omega$
will be typically dropped where is no danger of confusion.
Analogously, angle brackets $\left<\cdot,\cdot\right>_{V^*,V}$ denote
a~duality paring between spaces $V^*$ and $V$, where $V^*$ denotes the
dual of $V$; the subscript can be omitted. For any vector-valued function
$\bv$, the symmetric part of the gradient is defined through
$\bD\bv\coloneqq\tfrac12\bigl(\nabla\bv+(\nabla\bv)^\top\bigr)$.
We use notation $L^q(0,T;X)$ and $\mathcal{C}^k(I;X)$
to denote Bochner spaces of functions with values in the Banach
space $X$ and $k$-times continuously differentiable functions
on interval $I\subset\RR$ with values in~$X$ respectively;
$\mathcal{C}^k_0(I;X)$ denotes functions from $\mathcal{C}^k(I;X)$
which are compactly supported in~$I$.

The function spaces relevant to the problems that are being
investigated vary depending on the type of boundary conditions.
Two cases are being considered separately.
First, the case of the \noslip{} boundary
condition (the \noslip{} case, in short) and then other boundary
conditions that involve various kinds of slipping mechanisms
(the \slip{} case, in short), cf. Table~\ref{tab_summary_bc}.

\subsubsection{No-slip case}

We consider the space of compactly supported smooth functions and its
subspace of solenoidal functions:
\begin{equation*}
  \Cinfz \coloneqq \left\{
    \bv: \Omega\rightarrow\RR^3;\, \bv \text{ smooth};\,
    \supp\bv\subset\Omega
  \right\},
  \quad
  \Cinfzdiv \coloneqq \left\{
    \bv\in\Cinfz;\, \diver\bv=0
  \right\}
\end{equation*}
and their closures in $L^p$-norm, $W^{1,p}$-norm (with $1<p<\infty$)
and $W^{3,2}$-norm:
\begin{align*}
  \begin{aligned}
    \Wz{p} &\coloneqq
    \overline{ \Cinfz }^{\|\cdot\|_{W^{1,p}}},
    \\
    \Lndiv{p} &\coloneqq
    \overline{ \Cinfzdiv }^{\|\cdot\|_{L^p}},
    \\
    \Wzdiv{p} &\coloneqq
    \overline{ \Cinfzdiv }^{\|\cdot\|_{W^{1,p}}},
  \end{aligned}
  &&
  \begin{aligned}
    \Vz &\coloneqq
    \overline{ \Cinfz }^{\|\cdot\|_{W^{3,2}}},
    \\
    \Vzdiv &\coloneqq
    \overline{ \Cinfzdiv }^{\|\cdot\|_{W^{3,2}}}.
  \end{aligned}
\end{align*}
As a~consequence of the Poincar\'{e} and Korn inequalities, see
\cite[Corollary~6.31]{adams-fournier-2003},
\cite[Theorem 5.15]{diening-ruzicka-schumacher-2010},
the following norms are equivalent on~$\Wz{p}$
(and $\Wzdiv{p}\subset\Wz{p}$) for $1<p<\infty$:
\begin{align}
  \label{eq:Wzequivnorms}
  \|\bD\bv\|_p \leq \|\nabla\bv\|_p \leq \|\bv\|_{1,p}
  \leq C_\mathrm{P} \|\nabla\bv\|_p \leq C_\mathrm{P} C_\mathrm{K} \|\bD\bv\|_p
  \quad \text{for all } \bv\in\Wz{p},
\end{align}
with $\|\nabla\bv\|_p\coloneqq\||\nabla\bv|\|_p$,
$\|\bD\bv\|_p\coloneqq\||\bD\bv|\|_p$;
the constant $C_\mathrm{P}>0$ that appears due to the Poincar\'{e} inequality
depends on $p$ and $\Omega$, while the constant $C_\mathrm{K}>0$
that appears due to the Korn inequality depends only on $p$.

Note that for a~domain $\Omega$ (without further regularity
assumption on the smoothness of $\partial \Omega$)
we have the embedding $ \Vzdiv \hookrightarrow \Vz
\hookrightarrow W^{1,\infty}(\Omega)^3 $.
If additionally $\Omega$ is a~$C^{0,1}$ domain, i.e., $\Omega$
is a~domain with Lipschitz boundary~$\partial\Omega$,
then the following characterization holds true:
\begin{equation}
  \label{eq:defW1inf0}
  \begin{aligned}
    \Wz{p} &= \left\{
      \bv\in W^{1,p}(\Omega)^3;\,
      \bv=\b0 \text{ on }\partial\Omega\text{ in the sense of traces}
    \right\},
    \\
    \Wzdiv{p} &= \left\{
      \bv\in \Wz{p};\, \diver\bv=0
    \right\};
  \end{aligned}
\end{equation}
moreover we use~\eqref{eq:defW1inf0} as a~definition of
$ \Wz{\infty} $ and $ \Wzdiv{\infty} $
in the case that $p=\infty$ and $\Omega$ is a~$C^{0,1}$ domain.

We can occasionally denote the norm on $\bigl(\Wz{p}\bigr)^*$,
the topological dual of $\Wz{p}$, by $\|\cdot\|_{-1,p'}$.

\subsubsection{Slip case}

Here we assume $\Omega$ is a~$C^{0,1}$ domain. We denote by
$\bn:\partial\Omega\rightarrow\RR^3$ the unit outer normal
vector to~$\partial\Omega$.
The space of smooth vector-valued functions with
vanishing normal component on the boundary and its solenoidal
subspace are then introduced through:
\begin{align*}
  \Cinfn &\coloneqq \left\{
    \bv: \Omega\rightarrow\RR^3;\, \bv \text{ smooth};\,
    \bv\cdot\bn=0 \text{ on } \partial\Omega
  \right\},
  \\
  \Cinfndiv &\coloneqq \left\{
    \bv\in \Cinfn;\, \diver\bv=0
  \right\}.
\end{align*}
Since $\partial\Omega$ is Lipschitz
we can define the following spaces with vanishing normal trace:
\begin{equation*}
  \begin{aligned}
    \Vn &\coloneqq \left\{
      \bv \in W^{3,2}(\Omega)^3;\,
      \bv\cdot\bn=0 \text{ on } \partial\Omega
    \right\},
    \\
    \Vndiv &\coloneqq \Vn \cap \Lndiv2,
  \end{aligned}
\end{equation*}
and subsequently, for $1<p\leq\infty$,
\begin{equation*}
  \Wn{p} \coloneqq
  \overline{ \Vn }^{\|\cdot\|_{W^{1,p}}},
  \quad
  \Wndiv{p} \coloneqq
  \overline{ \Vndiv }^{\|\cdot\|_{W^{1,p}}}.
\end{equation*}
The condition $\bv\cdot\bn=0$ on $\partial\Omega$ for $\Omega$ bounded
is sufficient for validity of the Poincar\'{e} inequality:\footnote{%
  To verify it, assume for the sake of contradicion,
  that there is $\{\bv^j\}_{j=1}^\infty\subset\Wn{p}$
  with $\|\nabla\bv^j\|_p\rightarrow0$ and $\|\bv^j\|_{1,p}=1$. Relying on
  the compact Sobolev embedding, it follows that there is a~(not relabeled)
  subsequence $\{\bv^j\}_{j=1}^\infty$ which converges strongly
  in $\bL^p$ to some $\bv\in \bL^p$. This implies
  that $\{\bv^j\}_{j=1}^\infty$ is a~Cauchy sequence in $\Wn{p}$ and
  converges in $\Wn{p}$ to $\bv\in\Wn{p}$ with $\nabla\bv=\bO$.
  Hence $\bv$ is constant and by virtue of the boundedness of $\Omega$ and
  the boundary condition it follows that $\bv=\b0$, which is a~contradiction.
} for $1<p<\infty$ there exists $C_\mathrm{P}>0$
depending on $p$ and $\Omega$ such that
\begin{align}
  \label{eq:Wnpoincare}
  \|\nabla\bv\|_p \leq \|\bv\|_{1,p}
  \leq C_\mathrm{P} \|\nabla\bv\|_p
  \qquad \text{for all } \bv\in\Wn{p}.
\end{align}
For the steady problem
we will consider two inequalities of the Korn type depending on whether
the type of considered boundary conditions leads to the control
of the trace of $\bv$ on the boundary or not\footnote{%
  A~priori estimates for a~steady problem subject to
  a~\slip{}-type condition given by a~maximal monotone $2$-graph
  with~($\cB$\ref{Bcoercive}) will ensure control
  of~$\|\bv_{\btau}\|_{2,\partial\Omega}$. On the other hand, under
  the \freeslip{} condition~\eqref{F10} there is~no a~priori control over
  the tangential velocity~$\|\bv_{\btau}\|_{2,\partial\Omega}$ and
  rigid motions, if admissible in~$\Wn{\infty}$, prevent one to obtain
  a~steady solution.
}.
In the first case, it follows from
\cite[Lemma~1.11]{bulicek-malek-rajagopal-2007} and~\eqref{eq:Wnpoincare}:
for $1<p<\infty$ there exists $C_\mathrm{K}>0$ depending on $p$ and $\Omega$
such that
\begin{align}
  \label{eq:Wnkorntrace}
  \|\nabla\bv\|_p
  \leq
  C_\mathrm{K} \bigl( \|\bD\bv\|_p + \|\bv\|_{2,\partial\Omega} \bigr)
  \qquad \text{for all } \bv\in\Wn{p}
  \text{ with } \bv_{\btau}\in\bL^2(\partial\Omega).
\end{align}
The second situation when $\bs=\b0$ on $\partial\Omega$
requires us to rule out domains that admit nontrivial rigid motions.
We say that $\Omega$ is axisymmetric if there exists a~rigid body
motion tangential to boundary, i.e., there is $\bv\in\Wn{\infty}$
with $\bD\bv=\bO$ and $\nabla\bv\neq\bO$ constant in $\Omega$.
In the other words, there is $\bv\in\Wn{\infty}$ of the form
$\bv(\bx) = \bQ(\bx-\bx_0)$ for some $\bQ\subset\RR^{3\times 3}$
non-zero skew-symmetric matrix and constant $\bx_0\in\RR^3$.
From~\cite[Theorem~11, Remark~12]{bauer-pauly-2016} it follows
that if $\Omega$ is not axisymmetric and $1<p<\infty$ there
exists $C_\mathrm{K}>0$ depending on $\Omega$ and $p$ such that
\begin{align}
  \label{eq:Wnkornaxi}
  \|\nabla\bv\|_p
  \leq
  C_\mathrm{K} \|\bD\bv\|_p
  \qquad \text{for all } \bv\in\Wn{p}.
\end{align}
For the unsteady case we will use the following Korn-type
inequality:\footnote{%
  This is a~consequence of another Korn-type inequality:
  \begin{align}
    \label{eq:WnkornLp}
    \|\nabla\bv\|_p
    \leq
    C'_\mathrm{K} \bigl( \|\bD\bv\|_p + \|\bv\|_p \bigr)
    \qquad \text{for all } \bv\in\bW^{1,p}(\Omega);
  \end{align}
  see \cite[Theorem~1.10]{MNRR96}.
  To verify~\eqref{eq:WnkornL1}, assume that there is
  $\{\bv^j\}_{j=1}^\infty\subset\Wn{p}$ such that $\|\nabla\bv^j\|_p=1$
  and $\|\bD\bv^j\|_p + \|\bv^j\|_1 \rightarrow0$. Poincar\'{e}
  inequality~\eqref{eq:Wnpoincare} implies that $\{\bv^j\}_{j=1}^\infty$
  is bounded in $\Wn{p}$. By virtue of its reflexivity there is
  a~(not relabeled) subsequence such that $\bv^j\rightharpoonup\bv$
  weakly in~$\Wn{p}$ and by the Sobolev embedding $\bv^j\rightarrow\bv$
  in~$\bL^p$. On the other hand $\bv^j\rightarrow\b0$ in~$\bL^1$ hence
  by uniqueness of the limit we conclude $\bv=\b0$. Summarizing, the right-hand side
  of~\eqref{eq:WnkornLp} (with $\bv^j$ in place of $\bv$) goes to zero
  but the left-hand side is equal to unity, which is a~contradiction.

  We can see that~\eqref{eq:WnkornL1} in fact holds independently
  of the considered boundary condition, i.e., for all $\bW^{1,p}$.
}
\begin{align}
  \label{eq:WnkornL1}
  \|\nabla\bv\|_p
  \leq
  C_\mathrm{K} \bigl( \|\bD\bv\|_p + \|\bv\|_1 \bigr)
  \qquad \text{for all } \bv\in\Wn{p}.
\end{align}

\bigskip
\subsection{Analysis of steady flows} \label{sec:steady}

In this section, we investigate internal flows that are independent
of time. Under such circumstances, the governing system of equations
takes the form
\begin{subequations}
\label{S.1234}
\begin{align}
  \diver \bv &= 0
  & &\text{ in }\Omega ,\\
  \diver\left( \bv\otimes\bv - \bS \right) &= -\nabla p + \bb
  & &\text{ in }\Omega ,\\
  \left(\bS,\bD\bv\right) &\in \cG
  & &\text{ in }\Omega ,\\
  \bv\cdot\bn &= 0
  & &\text{ on }\partial\Omega ,
\end{align}
\end{subequations}
where $\cG$ is a~maximal monotone $r$-graph, fulfilling
($\cG$\ref{Gzero})--($\cG$\ref{Gcoercive}), of the form
\begin{align}
  \label{S.3*}
  \cG \coloneqq
  \left\{
    (\bS,\bD) \in\RR^{3\times 3}_\mathrm{sym}\times\RR^{3\times 3}_\mathrm{sym}
    ;\,
    \bS = 2\nu_* \left(|\bD|-\delta_*\right)^+ \mathcal{S}(|\bD|) \tfrac{\bD}{|\bD|}
  \right\}.
\end{align}

We will distinguish two cases depending on the class of
boundary conditions considered.
The first case concerns the \noslip{} condition, i.e.,
\begin{align}
  \label{S.5}
  \bv_{\btau} &= \b0
  & &\text{ on }\partial\Omega .
\end{align}
The second case includes all other boundary conditions involving
tangential part of the normal traction; it refers to either
\begin{align*}
  \bs &= \b0
  & &\text{ on }\partial\Omega
\end{align*}
or
\begin{align*}
  \left(\bs,\bv_{\btau}\right) &\in \cB
  & &\text{ on }\partial\Omega ,
\end{align*}
where $\cB$ is a~maximal monotone $2$-graph fulfilling
($\cB$\ref{Bzero})--($\cB$\ref{Bcoercive}).

\subsubsection{No-slip case}

Let $\Omega\subset\RR^3$ be a~domain, $r\geq\frac65$,
$\bb\in\bigl(\Wz{r}\bigr)^*$ and $\cG$ be a~maximal monotone $r$-graph
specified in~\eqref{S.3*}. We say that
\begin{align*}
  (\bv,\bS) \in \Wzdiv{r} \times L^{r'}(\Omega)^{3\times 3}_\mathrm{sym}
\end{align*}
is a~weak solution to \eqref{S.1234}, \eqref{S.3*}, \eqref{S.5} if
\begin{align}
  \label{S.9}
  (\bS,\bD\bfi) - (\bv\otimes\bv,\nabla\bfi) = \left<\bb,\bfi\right>
  \quad\text{ for all } \bfi\in\Cinfzdiv
\end{align}
and
\begin{align}
  \label{S.10}
  (\bS,\bD\bv) \in \cG
  \quad\text{ a.e. in } \Omega;
\end{align}
equivalently, we can require that \eqref{S.9} holds for all
$\bfi\in\Wzdiv{\frac{3r}{5r-6}}\cap\Wzdiv{r}$.

\begin{Theorem}
  Let $\Omega\subset\RR^3$ be a~domain.
  Let $r>\frac65$,
  $\bb\in\bigl(\Wzdiv{r}\bigr)^*$ and $\cG$ be
  a~maximal monotone $r$-graph of the form~\eqref{S.3*}.
  Then there is a~weak solution
  $
  (\bv,\bS) \in \Wzdiv{r} \times L^{r'}(\Omega)^{3\times 3}_\mathrm{sym}
  $
  to \eqref{S.1234}, \eqref{S.5} which fulfills~\eqref{S.10} and
  \begin{equation}
  \begin{aligned}
    \label{eq:SAmomentum}
    (\bS,\bD\bfi) - (\bv\otimes\bv,\bD\bfi)
    =&\;
    \left<\bb,\bfi\right>
    \\
    &\text{for all } \bfi\in \biggl\{\begin{array}{ll}
      \Wzdiv{r}               & \text{if } r\geq\frac95, \\
      \Wzdiv{\frac{3r}{5r-6}} & \text{if } r\in\bigl(\frac65,\frac95\bigr).
    \end{array}
  \end{aligned}
  \end{equation}
  Assume additionally that
  $\Omega$ is a~$C^{0,1}$ domain and
  $\bb\in\bigl(\Wz{r}\bigr)^*$. Then there is
  \begin{align}
    \label{eq:SApressure}
    p\in \biggl\{\begin{array}{ll}
      L^{r'}(\Omega)                & \text{if } r\geq\frac95 \\
      L^{\frac{3r}{2(3-r)}}(\Omega) & \text{if } r\in\bigl(\frac65,\frac95\bigr)
    \end{array},
    \quad
    \int_\Omega p\,\d x=0,
  \end{align}
  such that
  \begin{equation}
  \begin{aligned}
    \label{S.11}
    (\bS,\bD\bfi) - (\bv\otimes\bv,\nabla\bfi)
    =&\;
    (p,\diver\bfi) + \left<\bb,\bfi\right>
    \\
    &\text{for all } \bfi\in \biggl\{\begin{array}{ll}
      \Wz{r}               & \text{if } r\geq\frac95, \\
      \Wz{\frac{3r}{5r-6}} & \text{if } r\in\bigl(\frac65,\frac95\bigr).
    \end{array}
  \end{aligned}
  \end{equation}
\end{Theorem}
\begin{proof}
  \textbf{The case $r\geq\frac{9}{5}$.}
  Since $\Wzdiv{r}$ is separable, there is a~countable basis
  denoted by $\left\{\bom^r\right\}_{r=1}^\infty$. For $N\in\mathbb{N}$
  arbitrary, but fixed, we first look for the vector
  $\bc^N=\left(c^N_1,\ldots,c^N_N\right)\in\RR^N$ such that
  \begin{align*}
    \bv^N(x) \coloneqq \sum_{r=1}^N c^N_r \bom^r(x)
  \end{align*}
  satisfies the system of $N$ nonlinear equations
  \begin{align}
    \label{S.12}
    \left(\bS^N,\bD\bom^r\right) - \left(\bv^N\otimes\bv^N,\nabla\bom^r\right)
    = \left<\bb,\bom^r\right>,
    \quad r=1,\ldots,N,
  \end{align}
  where
  \begin{align}
    \label{S.13}
    \bS^N \coloneqq \bS\left(\bD\bv^N\right) \coloneqq
    2\nu_* \left(\left|\bD\bv^N\right|-\delta_*\right)^+
    \mathcal{S}\left(\left|\bD\bv^N\right|\right)
    \tfrac{\bD\bv^N}{\left|\bD\bv^N\right|}.
  \end{align}
  Introducing the (continuous) mapping $\vec{P}^N:\RR^N\rightarrow\RR^N$ as
  \begin{align*}
    \Bigl(\vec{P}^N\bigl(\bc^N\bigr)\Bigr)_r \coloneqq
    \left(\bS^N,\bD\bom^r\right) - \left(\bv^N\otimes\bv^N,\nabla\bom^r\right)
    -\left<\bb,\bom^r\right>,
    \quad r=1,\ldots,N,
  \end{align*}
  then
  \begin{align}
    \label{S.135}
    \vec{P}^N\bigl(\bc^N\bigr)\cdot\bc^N = \left(\bS^N,\bD\bv^N\right) -\left<\bb,\bv^N\right>.
  \end{align}
  It follows from~($\cG$\ref{Gcoercive}) and~\eqref{S.13} that
  \begin{align}
    \label{S.14}
    \vec{P}^N\bigl(\bc^N\bigr)\cdot\bc^N > 0 \quad\text{ for } |\bc^N| \text{ sufficiently large}.
  \end{align}
  As a~consequence of Brouwer's fixed-point theorem (see \cite[p.~53]{lions1969}),
  \eqref{S.14} implies the existence of $\bc^N$ fulfilling
  $\vec{P}^N\bigl(\bc^N\bigr)=\b0$, i.e., \eqref{S.12}~holds,
  and, by~\eqref{S.135},
  \begin{align}
    \label{S.136}
    \left(\bS^N,\bD\bv^N\right) = \left<\bb,\bv^N\right>.
  \end{align}
  This together with ($\cG$\ref{Gcoercive}), \eqref{eq:Wzequivnorms}
  and Young's inequality leads to
  \begin{align*}
    \|\bS^N\|_{r'} + \|\nabla\bv^N\|_r \leq
    c_1 \|\bb\|_{\bigl(\Wz{r}\bigr)^*} + c_2.
  \end{align*}
  This implies the existence of
  $\bv\in \Wzdiv{r}$ and $\bS\in L^{r'}(\Omega)^{3\times 3}_\mathrm{sym}$
  such that for suitable (not relabeled) subsequences
  \begin{subequations}
  \label{S.131}
  \begin{align}
    \label{S.131a}
    \bv^N    &\rightharpoonup \bv    & \text{weakly in } & \Wzdiv{r},      \\
    \bD\bv^N &\rightharpoonup \bD\bv & \text{weakly in } & L^r(\Omega)^{3\times 3}_\mathrm{sym}, \\
    \bS^N    &\rightharpoonup \bS    & \text{weakly in } & L^{r'}(\Omega)^{3\times 3}_\mathrm{sym},
    \intertext{
      as $N\rightarrow\infty$.
      Consequently, as $W^{1,q}_0(\Omega)$ is compactly embedded into
      $L^2(\Omega)$ for any $q>\tfrac65$,
      we also have
    }
    \bv^N &\rightarrow \bv & \text{strongly in } & \bL^2
    \text{ as } N\rightarrow\infty.
  \end{align}
  \end{subequations}
  Then, \eqref{S.12} leads to, for $r\geq\tfrac95$,
  \begin{align}
    \label{S.16}
    \left(\bS,\bD\bom^s\right) - \left(\bv\otimes\bv,\nabla\bom^s\right)
    = \left<\bb,\bom^s\right>,
    \quad s=1,2,\ldots.
  \end{align}
  Note that the restriction $r\geq\frac95$ is due to the requirement
  that for $s\in\mathbb{N}$ arbitrary
  \begin{align*}
    \int_\Omega (\bv\otimes\bv)\mathbin{:}\nabla\bom^s \,\d x < \infty
    \quad\text{for } \bv,\bom^s\in \Wz{r}.
  \end{align*}
  Hence~\eqref{S.16} implies that
  \begin{align}
    \label{S.165}
    \left(\bS,\bD\bom\right) - \left(\bv\otimes\bv,\nabla\bom\right)
    = \left<\bb,\bom\right>
    \quad \text{ for all }\bom\in \Wzdiv{r},
  \end{align}
  which completes the proof of~\eqref{eq:SAmomentum} for $r\geq\tfrac95$.
  Taking $\bom=\bv$ in~\eqref{S.165} one obtains
  \begin{align}
    \label{S.17}
    \left(\bS,\bD\bv\right) = \left<\bb,\bv\right>.
  \end{align}
  Taking the limit with $N\rightarrow\infty$ in~\eqref{S.136},
  we conclude from~\eqref{S.17} and~\eqref{S.131a} that
  \begin{align*}
    \lim_{N\rightarrow\infty} \left(\bS^N,\bD\bv^N\right) =
    \left(\bS,\bD\bv\right).
  \end{align*}
  By virtue of the graph convergence lemma (see Lemma~\ref{lemma:conv}
  in Appendix) this implies together with~\eqref{S.131} that
  $\left(\bS,\bD\bv\right)\in\cG$, i.e., \eqref{S.10} holds.

  \noindent
  \textbf{The case $r\in\bigl(\frac65,\frac95\bigr)$.}
  In this case, we consider, for $\epsilon>0$, the following
  approximating problem:
  \begin{equation}
  \label{L1234}
  \begin{aligned}
    -\diver\left( \bS + \epsilon\,\bD\bv - \bv\otimes\bv \right) &= -\nabla p + \bb
    && \text{in }\Omega, \\
    \diver \bv &= 0
    && \text{in }\Omega, \\
    (\bS,\bD\bv) &\in \cG
    && \text{in }\Omega, \\
    \bv &= \b0
    && \text{on }\partial\Omega.
  \end{aligned}
  \end{equation}
  For fixed~$\epsilon$, the existence of a~weak solution to
  \eqref{L1234} follows from the above proof for the case $r\geq\frac95$.
  More precisely, following step-by-step the proof of existence via
  the Galerkin method used above we can show that, for $\epsilon>0$ fixed,
  there is
  $(\bv^\epsilon,\bS^\epsilon)\in \Wzdiv2 \times L^{r'}(\Omega)^{3\times 3}_\mathrm{sym}$
  such that\footnote{%
    In order to verify that
    \begin{align*}
      \limsup_{N\rightarrow\infty} \int_\Omega \bS^N\mathbin{:}\bD\bv^N
      \leq \int_\Omega \bS\mathbin{:}\bD\bv,
    \end{align*}
    one uses the identities
    \begin{align*}
      \int_\Omega \bS^N\mathbin{:}\bD\bv^N + \epsilon\int_\Omega \bigl|\bD\bv^N\bigr|^2
    &= \bigl< \bb, \bv^N \bigr>, \\
      \int_\Omega \bS\mathbin{:}\bD\bv + \epsilon\int_\Omega \bigl|\bD\bv\bigr|^2
    &= \bigl< \bb, \bv \bigr>,
    \end{align*}
    the weak lower semicontinuity of the $L^2$-norm of $\bD\bv^N$,
    and the inequality
    \begin{align*}
      \liminf_{n\rightarrow\infty} a_n + \limsup_{n\rightarrow\infty} b_n
      \leq
      \limsup_{n\rightarrow\infty} (a_n + b_n)
    \end{align*}
    applied to any sequences $\{a_n\}_{n=1}^\infty$, $\{b_n\}_{n=1}^\infty$
    with $a_n\geq0$, $b_n\geq0$ for all $n\in\mathbb{N}$.

  }
  \begin{subequations}
  \begin{align}
    \label{L5}
    \left(
      \bS^\epsilon + \epsilon\,\bD\bv^\epsilon - \bv^\epsilon\otimes\bv^\epsilon,
      \bD\bfi
    \right)
    &= \left< \bb, \bfi \right>
    \text{ for all } \bfi\in \Wzdiv2, \\
    (\bS^\epsilon,\bD\bv^\epsilon)
    &\in \cG
    \text{ a.e. in } \Omega.
  \end{align}
  \end{subequations}
  Moreover, taking $\bfi=\bv^\epsilon$ in~\eqref{L5},
  \begin{align*}
    (\bS^\epsilon,\bD\bv^\epsilon) + \epsilon\,\|\bD\bv^\epsilon\|_2^2
    = \left< \bb, \bv^\epsilon \right>.
  \end{align*}
  The last identity together with the assumption~($\cG$\ref{Gcoercive})
  and Korn's inequality~\eqref{eq:Wzequivnorms}
  implies the following estimate
  \begin{align}
    \label{L8}
    \|\bS^\epsilon\|_{r'}^{r'}
    + \|\bD\bv^\epsilon\|_r^r + \|\nabla \bv^{\epsilon}\|_r^r
    + \epsilon\,\|\bD\bv^\epsilon\|_2^2
    \leq c_1 \|\bb\|_{\bigl(\Wz{r}\bigr)^*}^{r'} + c_2.
  \end{align}
  This implies the existence of
  $
    (\bv,\bS)\in \Wzdiv{r}\times L^{r'}(\Omega)^{3\times 3}_\mathrm{sym}
  $
  such that, for a~suitable vanishing subsequence
  $\{\epsilon_n\}_{n=1}^\infty$ and
  $(\bv^n,\bS^n)\coloneqq(\bv^{\epsilon_n},\bS^{\epsilon_n})$,
  \begin{subequations}
  \begin{align}
    \label{L9}
    \bS^n    &\rightharpoonup \bS    && \text{weakly in   } L^{r'}(\Omega)^{3\times 3}_\mathrm{sym}, \\
    \label{L10}
    \bD\bv^n &\rightharpoonup \bD\bv && \text{weakly in   } L^{r }(\Omega)^{3\times 3}_\mathrm{sym}, \\
    \label{L10a}
    \bv^n &\rightharpoonup \bv && \text{weakly in   } \Wzdiv{r}.
    \intertext{
      As $r<\tfrac{9}{5}$ $(<3)$ and $W^{1,r}_0(\Omega)$ is compactly
      embedded into $L^q(\Omega)$ for any $q\in[1,\tfrac{3r}{3-r})$,
      we also have
    }
    \label{L11}
    \bv^n &\rightarrow        \bv && \text{strongly in } \bL^{q }
    \text{ for all } q\in\bigl[1,\tfrac{3r}{3-r}\bigr).
    \intertext{
      Limit~\eqref{L11} provides the strong convergence
    }
    \label{L111}
     \bv^n &\rightarrow \bv  && \text{strongly in } \bL^{2}
  \end{align}
  \end{subequations}
  provided that
  $\frac{3r}{3-r}>2$, which gives the bound stated in the formulation
  of the theorem, namely $r>\frac65$. Consequently, $\bv$ and $\bS$
  fulfill~\eqref{eq:SAmomentum}. It remains to prove that
  $(\bS,\bD\bv)\in\cG$ a.e. in~$\Omega$. By the graph convergence
  lemma (Lemma~\ref{lemma:conv}), it is enough to show that
  \begin{align*}
    \limsup_{n\rightarrow\infty}\, \bigl(\bS^n,\bD\bv^n\bigr) \leq \bigl(\bS,\bD\bv\bigr).
  \end{align*}
  To prove it, we first subtract~\eqref{eq:SAmomentum} from~\eqref{L5}
  to obtain
  \begin{align}
    \label{L12}
    \begin{multlined}
      \left(\bS^n-\bS,\bD\bfi\right)
      +\epsilon_n\,\left(\bD\bv^n,\bD\bfi\right)
      +\left(\bv\otimes\bv-\bv^n\otimes\bv^n,\bD\bfi\right)
      = 0
      \\
      \text{ for all } \bfi\in\Wzdiv{2}\cap\Wzdiv{\frac{3r}{5r-6}}.
    \end{multlined}
  \end{align}
  Let $B\subset\RR^3$ be an~arbitrary ball of radius $R$ such that
  $\tfrac{B}{2}\subset B\subset2B\subset\Omega$ and $\chi\in\D(B)$ be such
  that $\chi=1$ in~$\tfrac{B}{2}$, $0\leq\chi\leq1$ in~$B$, and $|\nabla\chi|\leq CR^{-1}$
  in~$B$. Then we set
  \begin{align*}
    \bu^n \coloneqq \chi \bigl( \bv^n - \bv \bigr) - \bh^n,
  \end{align*}
  where $\bh^n\in \Wz{r}(B)$ solves
  \begin{align*}
    \diver\bh^n = \nabla\chi\cdot \bigl( \bv^n - \bv \bigr)
    \qquad \text{in } B.
  \end{align*}
  (Note that the compatibility condition
  $\int_B \nabla\chi\cdot \bigl( \bv^n - \bv \bigr) = \int_B \diver\bigl(\chi ( \bv^n - \bv) \bigr) = 0$ is
  met.)
  In fact, there is a~continuous linear operator
  $\mathcal{B}:\{q\in L^p(B),\,\int_B q=0\}\rightarrow\Wz{p}(B):
  g\mapsto\bu$ such that $\diver\bu=g$,
  cf. Remark~\ref{rem:bog}~(\ref{prop:bog}).
  Consequently, $\bu^n$ extended by zero in $\Omega\setminus B$ fulfill
  $\diver\bu^n=0$ in~$\Omega$.
  Next we consider divergence-free Lipschitz approximations $\bu^{n,k}\in\Wzdiv{\infty}(2B)$
  to $\bu^n\in\Wzdiv{r}(B)$ from Lemma~\ref{lemma:truncsteady}.
  Taking $\bfi\coloneqq\bu^{n,k}$ in~\eqref{L12} and letting $n\rightarrow\infty$,
  we conclude, using in particular properties~(\ref{lemma:truncsteady:lambda}),
  (\ref{lemma:truncsteady:bound}), and~(\ref{lemma:truncsteady:weakstarconv})
  of Lemma~\ref{lemma:truncsteady},
  $\epsilon_n\rightarrow0$, \eqref{L10}, and~\eqref{L111},
  that
  \begin{align}
    \label{L13}
    \lim_{n\rightarrow\infty} \left( \bS^n-\bS, \bD\bu^{n,k} \right)_{2B}
    = 0
    \quad \text{for all } k\in\mathbb{N}.
  \end{align}
  \begingroup\setlength\emergencystretch{\hsize}\hbadness=10000
  On the other hand, using H\"older's inequality, \eqref{L8},
  and properties
  (\ref{lemma:truncsteady:bound}), (\ref{lemma:truncsteady:smallness})
  of Lemma~\ref{lemma:truncsteady}, there is a $C>0$ independent of $n$ and $k$ such that
  for all $n,k\in\mathbb{N}$
  \begin{equation}
    \label{L131}
    \left| \left( \bS^n-\bS, \bD\bu^{n,k} \right)_{\mathcal{O}^{n,k}} \right|
    \leq \|\bS^n-\bS\|_{r'}
    \|\bD\bu^{n,k}\|_{\infty,\mathcal{O}^{n,k}}
    |\mathcal{O}^{n,k}|^{\frac{1}{r}}
    \leq C\, 2^{-\frac{k}{r}}.
  \end{equation}
  \endgroup
  As $\bu^n$ vanishes outside of~$B$ and $\bu^{n,k} = \bu^{n}$ in
  $2B\setminus\mathcal{O}^{n,k}$ (see the property (\ref{lemma:truncsteady:lipapp})
  of Lemma~\ref{lemma:truncsteady}) we observe that
  \begin{equation}
    \label{L132}
    \left( \bS^n-\bS, \bD\bu^{n} \right)_{B\setminus\mathcal{O}^{n,k}}
    =
    \left( \bS^n-\bS, \bD\bu^{n,k} \right)_{2B\setminus\mathcal{O}^{n,k}}
    \quad \text{for all } n,k\in\mathbb{N}.
  \end{equation}
  Using \eqref{L13}--\eqref{L132}
  we conclude that
  \begin{align}
    \label{L133}
    \limsup_{n\rightarrow\infty}
    \bigl| \bigl(\bS^n-\bS,\bD\bu^n\bigr)_{B\setminus\mathcal{O}^{n,k}} \bigr|
    \leq C\, 2^{-\frac{k}{r}}
    \quad \text{for all } k\in\mathbb{N}.
  \end{align}
  It directly follows from
  the definition of~$\bu^n$
  that
  \begin{align*}
    \bD\bu^n
    = \chi\,\bD(\bv^n-\bv) + \tfrac12 \bigl[ \nabla\chi\otimes(\bv^n-\bv)
    + (\bv^n-\bv)\otimes\nabla\chi \bigr]
    - \bD\bh^n,
  \end{align*}
  where $\bh^n:=\mathcal{B}\bigl(\nabla\chi\cdot(\bv^n-\bv)\bigr)$.
  It then follows from~\eqref{L133},
  the properties of the operator~$\mathcal{B}$, and~\eqref{L11} that
  \begin{align}
    \label{L135}
    \limsup_{n\rightarrow\infty}
    \left|
    \int_{B\setminus\mathcal{O}^{n,k}}
    (\bS^n-\bS)\mathbin{:}(\bD\bv^n-\bD\bv)\, \chi
    \right|
    \leq C\, 2^{-\frac{k}{r}}
    \quad \text{for arbitrary } k\in\mathbb{N}.
  \end{align}
  Next, we set $\bS^*\coloneqq2\nu_*\,
  (|\bD\bv|-\delta_*)^{+}\, \mathcal{S}(|\bD\bv|)\, \bD\bv/|\bD\bv|$. Then, by virtue
  of definition~\eqref{S.3*}, it holds that
  $(\bS^*,\bD\bv)\in\cG$ and $\bS^*\in L^{r'}(\Omega)^{3\times 3}_\mathrm{sym}$.
  Hence \eqref{L135} and \eqref{L10} imply
  \begin{align}
    \label{L136}
    \limsup_{n\rightarrow\infty}
    \left|
    \int_{B\setminus\mathcal{O}^{n,k}}
    (\bS^n-\bS^*)\mathbin{:}(\bD\bv^n-\bD\bv)\, \chi
    \right|
    \leq C\, 2^{-\frac{k}{r}}
    \quad \text{for arbitrary } k\in\mathbb{N}.
  \end{align}
  Since $(\bS^*,\bD\bv)\in\cG$ and $(\bS^n,\bD\bv^n)\in\cG$
  the monotonicity property~($\cG$\ref{Gmonotone}) of~$\cG$ implies
  that the integrand in~\eqref{L136} is non-negative.
  Using H\"older's inequality and the properties (\ref{lemma:truncsteady:lambda}),
  (\ref{lemma:truncsteady:smallness}) of Lemma~\ref{lemma:truncsteady}
  we obtain
  \begin{align*}
    \limsup_{n\rightarrow\infty}
    \int_{B}
    \bigl( (\bS^n-\bS^*)\mathbin{:}(\bD\bv^n-\bD\bv) \bigr)^\frac12
    \chi^\frac12
    \leq C\, 2^{-\frac{k}{r}}
    \quad \text{for arbitrary } k\in\mathbb{N}.
  \end{align*}
  This leads to
  \begin{align*}
    \limsup_{n\rightarrow\infty}
    \int_{\frac{B}{2}} \bigl( (\bS^n-\bS^*)\mathbin{:}(\bD\bv^n-\bD\bv) \bigr)^\frac12
    \leq C\, 2^{-\frac{k}{r}}
    \quad \text{for arbitrary } k\in\mathbb{N}.
  \end{align*}
  Let us set $g^n\coloneqq (\bS^n-\bS^*)\mathbin{:}(\bD\bv^n-\bD\bv)$.
  Clearly $g^n\geq0$ and
  $g^n\rightarrow0$ almost everywhere in~$\tfrac{B}{2}$. But as $B$ is
  arbitrary, we conclude that
  \begin{align}
    \label{L14}
    g^n\rightarrow0 \text{ almost everywhere in } \Omega.
  \end{align}
  Since $\{g^n\}_{n=1}^\infty$ is bounded in $L^1(\Omega)$ and has the pointwise
  limit~\eqref{L14}, Corollary~\ref{col:strongbiting}, a~consequence
  of the biting lemma (Lemma~\ref{lemma:biting}), then implies
  existence of a~subsequence $\{g^{n_j}\}_{j=1}^\infty$,
  and a~sequence of sets $\{E_k\}_{k=1}^\infty$ with
  $\Omega\supset E_1\supset E_2\supset\ldots$, $|E_k|\rightarrow0$ such that for
  all $k\in\mathbb{N}$
  \begin{align*}
    g^{n_j} \rightarrow 0 \text{ strongly in } L^1(\Omega\setminus E_k).
  \end{align*}
  From the definition of $g^{n_j}$ we conclude that
  \begin{align*}
    \limsup_{j\rightarrow\infty}
     \int_{\Omega\setminus E_k} (\bS^{n_j}-\bS^*)\mathbin{:}(\bD\bv^{n_j}-\bD\bv) = 0,
  \end{align*}
  which implies as a~consequence of~\eqref{L9} and \eqref{L10}
  that for all $k\in\mathbb{N}$
  \begin{align*}
    \limsup_{j\rightarrow\infty} \int_{\Omega\setminus E_k} \bS^{n_j}\mathbin{:}\bD\bv^{n_j}
     = \int_{\Omega\setminus E_k} \bS\mathbin{:}\bD\bv.
  \end{align*}
  Since $|E_k|\rightarrow0$, we can conclude from the graph convergence
  lemma (Lemma~\ref{lemma:conv}) that
  \begin{align*}
    (\bS,\bD\bv)\in\cG\text{ almost everywhere in } \Omega
  \end{align*}
  so that \eqref{S.10} holds and the first part of the theorem is proved.

  \noindent
  \textbf{On the pressure.}
  Setting
  \begin{align*}
    \left<\bF,\bfi\right> \coloneqq
    \left(\bS,\bD\bfi\right) - \left(\bv\otimes\bv,\nabla\bfi\right)
    - \left<\bb,\bfi\right>
  \end{align*}
  we observe that
  \begin{align*}
    \left<\bF,\bfi\right>=0
    \quad\text{for all }\bfi\in\Cinfzdiv
  \end{align*}
  and
  \begin{align}
    \label{S.19}
    \bF\in \Biggl\{\begin{array}{ll}
      \Bigl(\Wz{r}\Bigr)^*               & \text{if } r\geq\frac95, \\
      \Bigl(\Wz{\frac{3r}{5r-6}}\Bigr)^* & \text{if } r\in\bigl(\frac65,\frac95\bigr).
    \end{array}
  \end{align}
  By the de~Rham theorem, see~\cite[Theorem~2.1]{ag94}, there is
  $p\in\bigl(\D(\Omega)\bigr)^*$ such that
  \begin{align}
    \label{S.20}
    \left<\bF,\bfi\right> = \left<-\nabla p,\bfi\right>
    \quad\text{for all }\bfi\in\Cinfz.
  \end{align}
  Since $\Omega$ is~$C^{0,1}$ domain, the Ne\v{c}as theorem
  (see Lemma~\ref{lemma:necas}, Remark~\ref{rem:necas})
  together with~\eqref{S.19} and~\eqref{S.20}
  implies~\eqref{eq:SApressure} and~\eqref{S.11}.
\end{proof}

\subsubsection{Slip case}

In this part we replace the \noslip{} boundary condition either by
\begin{align}
  \label{3.101}
  \bv\cdot\bn = 0
  \quad\text{ and }\quad
  \bs = \b0
  \quad\text{ on }\partial\Omega
\end{align}
or by
\begin{equation}
  \label{3.102}
  \begin{gathered}
    \bv\cdot\bn = 0
    \quad\text{ and }\quad
    (\bs,\bv_{\btau}) \in \cB
    \quad\text{ on }\partial\Omega,
    \\
    \text{ where } \cB \text{ fulfills the conditions }
    (\cB\ref{Bzero})\text{--}(\cB\ref{Bcoercive}).
  \end{gathered}
\end{equation}
We prove the following result.

\begin{Theorem} \label{SS}
  Let $\Omega\subset\RR^3$ be a~$C^{1,1}$~domain\footnote{%
    In the case of the boundary condition~\eqref{3.101},
    the required regularity of the boundary can be weakened at the cost
    of losing the information concerning the pressure.
  }. Let further $r>\frac65$, $\bb\in \bigl(\Wn{r}\bigr)^*$, $\cG$ be a~maximal
  monotone $r$-graph of the form~\eqref{S.3*}.
  \bigskip
  \begin{enumerate}[(i)]
    \item \emph{(Boundary condition~\eqref{3.102})}
      Let $\cB$ be a~maximal monotone $2$-graph.
      Then there is a~weak solution
      \begin{align*}
        (\bv,\bS,\bs) \in
        \Wndiv{r} \times L^{r'}(\Omega)^{3\times 3}_\mathrm{sym} \times \bL^2(\partial\Omega)
      \end{align*}
      to \eqref{S.1234} and \eqref{3.102} such that
      \begin{align}
        \label{SS_G}
        \bigl(\bS,\bD\bv\bigr) &\in \cG \text{ a.e. in }\Omega, \\
        \label{SS_B}
        \bigl(\bs,\bv_{\btau}   ) &\in \cB \text{ a.e. in }\partial\Omega,
      \end{align}
      and
        \begin{multline}
        \label{SS_div}
          \makebox[\leftmargin]{}
          \bigl(\bS,\bD\bfi\bigr) - \bigl(\bv\otimes\bv,\bD\bfi\bigr)
          + \bigr(\bs,\bfi\bigr)_{\partial\Omega}
          =\; \bigr<\bb,\bfi\bigr>
          \\
          \text{for all } \bfi\in\biggl\{\begin{array}{ll}
            \Wndiv{r}               & \text{if } r\geq\frac95, \\
            \Wndiv{\frac{3r}{5r-6}} & \text{if } r\in\bigl(\frac65,\frac95\bigr).
          \end{array}
        \end{multline}
    \item \emph{(Boundary condition~\eqref{3.101})}
      Assume that $\Omega$ is not axisymmetric.
      Then there is a~weak solution
      \begin{align*}
        (\bv,\bS) \in
        \Wndiv{r} \times L^{r'}(\Omega)^{3\times 3}_\mathrm{sym}
      \end{align*}
      to \eqref{S.1234} and \eqref{3.101} such that
      \eqref{SS_G} and \eqref{SS_div} with $\bs=\b0$ hold true.
  \end{enumerate}
  \bigskip
  In addition, there is
  \begin{align}
    \label{SS_p}
    p\in \biggl\{\begin{array}{ll}
      L^{r'}(\Omega)                & \text{if } r\geq\frac95 \\
      L^{\frac{3r}{2(3-r)}}(\Omega) & \text{if } r\in\bigl(\frac65,\frac95\bigr)
    \end{array},
    \quad
    \int_\Omega p\,\d x=0
  \end{align}
  such that
  \begin{multline}
    \label{SS_momentum}
    \bigl(\bS,\bD\bfi\bigr) - \bigl(\bv\otimes\bv,\bD\bfi\bigr)
    + \bigr(\bs,\bfi\bigr)_{\partial\Omega}
    =\;
    \bigl(p,\diver\bfi\bigr) + \bigr<\bb,\bfi\bigr>
    \\
    \text{for all } \bfi\in \biggl\{\begin{array}{ll}
      \Wn{r}               & \text{if } r\geq\frac95, \\
      \Wn{\frac{3r}{5r-6}} & \text{if } r\in\bigl(\frac65,\frac95\bigr).
    \end{array}
  \end{multline}
\end{Theorem}
\begin{proof}
  The proof of existence of $\bv$, $\bS$, $\bs$ with \eqref{SS_G}, \eqref{SS_B},
  and \eqref{SS_div} follows the same scheme as in the case of the \noslip{} boundary
  condition. The only differences are
  \begin{enumerate}[(i)]
    \item
      due to a~different choice of the
      function spaces for the velocity
      as $\Wzdiv{r}$ is replaced by
      $\Wndiv{r}$,
    \item
  due to the presence of the term
      $\int_{\partial\Omega} \bs\cdot\bfi$
      in the weak formulation of the
      balance of linear momentum,
    \item
      and due to the necessity to verify validity
      of the boundary condition
      $\bh(\bs,\bv_{\btau})\allowbreak=\b0$.
  \end{enumerate}
  Note that in the case $\bs=\b0$ on $\partial\Omega$, the last two
  differences disappear. In the case of the \noslipalt/\navierslipalt{} boundary condition
  $\bv_{\btau}=\frac{1}{\gamma_*}\frac{(|\bs|-\sigma_*)^+}{|\bs|}\bs$
  we start with its approximation
  \begin{align*}
    \bv_{\btau} = \frac{1}{\gamma_*}\frac{(|\bs|-\sigma_*)^+}{|\bs|}\bs
             + \frac{\epsilon}{\gamma_*}\bs,
    \quad \epsilon>0.
  \end{align*}
  Other boundary conditions of the type~\eqref{3.102} employ similar approximations.

  For $r\geq\frac95$, the proof then proceeds as in the case of \noslip{}
  boundary conditions up to the use of the appropriate form of Korn's inequality
  (\eqref{eq:Wnkornaxi}~for boundary condition~\eqref{3.101} assuming the domain is
  not axisymmetric,
  or \eqref{eq:Wnkorntrace}~for boundary condition~\eqref{3.102})
  and the following points. The additional term
  $\int_{\partial\Omega} \bs^N\cdot\bfi$ in the balance of linear momentum is
  treated using the weak convergence $\bs^N\rightharpoonup\bs$
  in $\bL^2(\partial\Omega)$. Since $W^{1,r}(\Omega)
  \smash[t]{{}\overset{\operatorname{tr}}{\hookrightarrow}{}}
  W^{1-\frac1r,r}(\partial\Omega) \hookrightarrow\hookrightarrow
  L^2(\partial\Omega)$ provided $r>\frac32$ and $\Omega$ is a~$C^{0,1}$ domain,
  the trace operator maps $\Wn{r}$ into  $\bL^2(\partial\Omega)$ compactly
  (i.e.,
  $\Wn{r} \smash[t]{{}\overset{\operatorname{tr}}{\hookrightarrow\hookrightarrow}{}} \bL^2(\partial\Omega)$)
  even for $r>\tfrac32$, and consequently
  $\bv^N_{\btau}\rightarrow\bv_{\btau}$
  strongly in~$\bL^2(\partial\Omega)$ and thus
  \begin{align}
    \label{N201}
    \int_{\partial\Omega} \bs^N\cdot\bv^N_{\btau}
    \rightarrow
    \int_{\partial\Omega} \bs\cdot\bv_{\btau}.
  \end{align}
  This implies (see Lemma \ref{lemma:conv}) that $(\bs,\bv_{\btau})\in\cB$~a.e.
  on~$\partial\Omega$.

  For $r\in\bigl(\frac65,\frac95\bigr)$, the validity of
  $(\bS,\bD\bv)\in\cG$ a.e. in $\Omega$ can be established in the same way
  as in the \noslip{} case since the proof is based on local (interior)
  arguments. It remains to show that $(\bs,\bv_{\btau})\in\cB$ a.e. on $\partial\Omega$.
  For $r>\tfrac32$, it follows from \eqref{N201} and Lemma~\ref{lemma:conv}. To use
  Lemma~\ref{lemma:conv} also for $r\in\bigl(\frac65,\frac32\bigr]$, it suffices to show that
  \begin{align*}
    \limsup_{N\rightarrow\infty} \int_{\partial\Omega} \bs^N\cdot\bv^N
    \leq \int_{\partial\Omega} \bs\cdot\bv.
  \end{align*}
  However, for $r>1$: $W^{1,r}(\Omega)
  \overset{\operatorname{tr}}{\hookrightarrow\hookrightarrow}
  L^1(\partial\Omega)$, the strong convergence $\bv^N_{\btau}\rightarrow\bv_{\btau}$
  in $\bL^1(\partial\Omega)$ together with Egorov's theorem implies
  that for any $\delta>0$ there is $\mathcal{U}_\delta\subset\partial\Omega$
  such that $|\partial\Omega\setminus\mathcal{U}_\delta|<\delta$ and
  $\bv^N_{\btau}\rightarrow\bv_{\btau}$ strongly in $\bL^\infty(\mathcal{U}_\delta)$. Hence
  \begin{align*}
    \limsup_{N\rightarrow\infty} \int_{\mathcal{U}_\delta} \bs^N\cdot\bv^N_{\btau}
    =
    \int_{\mathcal{U}_\delta} \bs\cdot\bv_{\btau}.
  \end{align*}
  Consequently, by the graph convergence lemma (Lemma~\ref{lemma:conv}
  in Appendix), $(\bs,\bv_{\btau})\in\cB$
  a.e. in $\mathcal{U}_\delta$. As $\delta>0$ was arbitrary, we conclude
  that $(\bs,\bv_{\btau})\in\cB$ a.e. on $\partial\Omega$.
  The proof of the first part of the theorem is complete.

  It remains to prove the existence of pressure~\eqref{SS_p}
  fulfilling~\eqref{SS_momentum}. Let us define a~linear
  functional $\bF$ through the relation
  \begin{align*}
    \bigl<\bF,\bfi\bigr>
      \coloneqq \bigl<\bb,\bfi\bigr>
      - \bigl(\bS-\bv\otimes\bv,\bD\bfi\bigr)
      - \bigr(\bs,\bfi\bigr)_{\partial\Omega}
      \quad\text{ for any }\bfi\in \Cinfn.
  \end{align*}
  From~\eqref{SS_div} we can see that $\bF\in \bigl(\Wn{q}\bigr)^*$ where
  $q = r$ if $r\geq\tfrac95$ and $q = \frac{3r}{5r-6}$ if $r\in(\tfrac65,\tfrac95)$
  and
  \begin{align}
    \label{SS_Fdiv}
    \bigl<\bF,\bfi\bigl> = 0
    \quad\text{ for all } \bfi\in \Wndiv{q}.
  \end{align}
  Now consider a~variational problem to find $p\in L^{q'}(\Omega)$ with
  $\int_\Omega p=0$ such that
  \begin{align}
    \label{SS_defp}
    \bigl(p,-\Delta\phi\bigr) = \bigl<\bF,\nabla\phi\bigr>
    \quad\text{ for all } \phi\in W^{2,q}(\Omega)
    \text{ with } \nabla\phi\in\Wn{q}.
  \end{align}
  As a~consequence of the $C^{1,1}$ smoothness of the domain, one can employ
  Lemma~\ref{lemma:neumannregularity} to conclude that~\eqref{SS_defp} is
  equivalent with the problem:
  find $p\in L^{q'}(\Omega)$ such that $\int_\Omega p=0$ and
  \begin{align}
    \label{SS_defp2}
    \bigl(p,q) = \bigl<\bF,\nabla A^{-1}q\bigr>
    \quad\text{ for all } q\in L^q(\Omega)
    \text{ with } \int_\Omega q = 0
  \end{align}
  where $A^{-1}$ is the solution operator for the
  Neumann-Poisson problem~\eqref{SS_laplace}.
  The problem~\eqref{SS_defp2} has a~unique solution by virtue of
  Lemma~\ref{lemma:neumannregularity}.
  Thus we have constructed $p$ with properties~\eqref{SS_p}.
  To verify~\eqref{SS_momentum}, consider a~test function $\bfi\in\Wn{q}$.
  With the Helmholtz decomposition (see Corollary~\ref{cor:helmholtz} in
  Appendix) $\bfi=\nabla\phi+\bfi_0$ with $\nabla\phi\in\Wn{q}$
  and $\bfi_0\in\Wndiv{q}$ we can immediately obtain,
  using~\eqref{SS_Fdiv} and~\eqref{SS_defp}, that
  \begin{align*}
    \bigl<\bF,\bfi\bigr> = \bigl<\bF,\nabla\phi\bigr> + \bigl<\bF,\bfi_0\bigr>
    = \bigl(p,-\Delta\phi\bigr) = -\bigl(p,\diver\bfi\bigr).
  \end{align*}
  This proves~\eqref{SS_momentum}. The proof of Theorem~\ref{SS}
  is thus complete.
\end{proof}

\bigskip
\subsection{Analysis of unsteady flows} \label{sec:unsteady}

In this section, we investigate unsteady internal flows, i.e.,
flows governed by \eqref{eq:sys}. Again,
we treat separately two cases: the \noslip{} boundary condition
and the boundary conditions allowing \slip.

\subsubsection{No-slip case}

We first provide an~existence result for the \noslip{} case, i.e., we investigate
the system~\eqref{eq:sysdiv}--\eqref{eq:sysbcnor}, \eqref{eq:sysic},
and $\bv_{\btau}=\b0$ on $(0,T)\times\partial\Omega$ as
a~special case of~\eqref{eq:sysbctan}.

\begin{Theorem} \label{th:noslip}
  Let $T\in(0,\infty)$, $\Omega\subset\RR^3$ be a~domain and
  $Q\coloneqq (0,T)\times\Omega$.
  Let $r>\frac65$,
  $\bb\in L^{r'}\bigl(0,T;(\Wz{r})^*\bigr)$ and $\bv_0\in\Lndiv2$.
  Let further
  $\cG\subset\RR^{3\times 3}_\mathrm{sym}\times\RR^{3\times 3}_\mathrm{sym}$
  be a~maximal monotone $r$-graph of the form~\eqref{S.3*} fulfilling
  ($\cG$\ref{Gzero})--($\cG$\ref{Gcoercive}).
  Then there exists a~pair $(\bv,\bS)$:
  \begin{subequations}
  \begin{align}
    \label{eq:P1}
    \bv &\in L^\infty(0,T;\Lndiv{2})
             \cap L^r(0,T;\Wzdiv{r}),
    \\
    \bS &\in L^{r'}(Q)^{3\times 3}_\mathrm{sym}
  \end{align}
  satisfying
  \begin{gather}
    \label{eq:P3}
    \lim_{t\rightarrow0+} \int_\Omega |\bv(t,\cdot)-\bv_0|^2 = 0,
    \\
    \begin{aligned}
      \label{eq:P4}
      \int_Q \bS\mathbin{:}\bD\bw
      =&\;
      \int_0^T \left< \bb, \bw \right>
      + \int_Q \bv\otimes\bv\mathbin{:}\bD\bw
      + \int_Q \bv\cdot\pp{\bw}{t}
      + \int_\Omega \bv_0\cdot\bw(0,\cdot)
      \\
      &\text{for all } \bw\in \D\bigl([0,T);\Wzdiv{q}\bigr),\,
      q=\max\bigl\{r,\tfrac{5r}{5r-6}\bigr\},
    \end{aligned}
    \\
    \label{eq:P5}
    \bS = 2\nu_* \left( |\bD\bv|-\delta_* \right)^+
          \mathcal{S}(|\bD\bv|) \tfrac{\bD\bv}{|\bD\bv|}
    \quad\text{almost everywhere in } Q.
  \end{gather}
  \end{subequations}
  Also
  \begin{subequations}
  \label{eq:cont}
  \begin{alignat}{2}
    \label{eq:contweak}
    \bv &\in \mathcal{C}\bigl([0,T]; \bL^2_\mathrm{weak}\bigr)
    && \text{ if } r>\tfrac{6}{5},
    \\
    \label{eq:contstrong}
    \bv &\in \mathcal{C}\bigl([0,T]; \bL^2\bigr)
    && \text{ if } r\geq\tfrac{11}{5},
  \end{alignat}
  \end{subequations}
  i.e., $t\mapsto(\bv(t),\bfi)$ is continuous on $[0,T]$ for all $\bfi\in\bL^2$ if $r>\tfrac65$
  and $t\mapsto\|\bv(t)\|_2$ is continuous on $[0,T]$ if $r\geq\tfrac{11}{5}$, respectively.
  Moreover, the energy inequality holds:
  \begin{align}
    \begin{gathered}
      \label{eq:P6}
      \int_\Omega \frac{|\bv(t,\cdot)|^2}{2}
      + \int_0^t \int_\Omega \bS\mathbin{:}\bD\bv
      \leq
      \int_\Omega \frac{|\bv_0|^2}{2}
      + \int_0^t \left<\bb, \bv\right>
      \qquad
      \text{for all } t\in[0,T];
    \end{gathered}
  \end{align}
  if $r\geq\frac{11}{5}$, \eqref{eq:P6} becomes equality.

  \begingroup\setlength\emergencystretch{\hsize}\hbadness=10000
  In addition, if $\Omega$ is a~$C^{0,1}$ domain
  with sufficiently small Lipschitz constant (smallness depending only on $r$)
  or $\Omega$ is any $C^1$ domain,
  then there are $P^1\in L^\infty(0,T;L^6(\Omega))$,
  $P^1(t,\cdot)$ harmonic in $\Omega$ for almost every $t\in(0,T)$
  and $p^2\in L^{q'}(Q)$ with
  $q=\max\bigl\{r,\frac{5r}{5r-6}\bigr\}$ such that
  \begin{align} \label{eq:P5p}
    \begin{aligned}
      \int_Q \bS\mathbin{:}\bD\bom
      &= \int_0^T \bigl<\bb,\bom\bigr>
      + \int_Q \bv\otimes\bv\mathbin{:}\bD\bom
      + \int_Q \bv\cdot\pp{\bom}{t}
      + \int_\Omega \bv_0\cdot\bom(0,\cdot)
      \\& - \int_Q P^1 \diver\pp{\bom}{t}
      + \int_Q p^2\diver\bom
      \begin{gathered}
        \qquad\text{for all } \bom\in \D\bigl([0,T);\Wz{q}\bigr),\\
         q=\max\bigl\{r,\tfrac{5r}{5r-6}\bigr\}.
      \end{gathered}
    \end{aligned}
  \end{align}
  Functions $P^1$ and $p^2$ can be chosen such that
  $\int_\Omega P^1(t,\cdot)=\int_\Omega p^2(t,\cdot)=0$
  for almost every $t\in(0,T)$.
  If, in addition, $\Omega$ is a~$C^{1,1}$ domain then it holds
  $\nabla P^1\in L^\infty(0,T;L^2(\Omega)) \cap L^\frac{5r}{3}(Q)$.
  \endgroup
\end{Theorem}
\begin{Rem}
      We could define the weak solution to the problem considered differently.
      We could say that $\bv$ is a~weak solution to the problem if
      $\bv$ fulfills~\eqref{eq:P1},
      \eqref{eq:P3} and
      \begin{multline*}
        \begin{aligned}
          \int_Q 2\nu_* \left( |\bD\bv|-\delta_* \right)^+
                 \mathcal{S}(|\bD\bv|) \tfrac{\bD\bv}{|\bD\bv|}\mathbin{:}\bD\bw
          & =
          \int_0^T \left< \bb, \bw \right>
          + \int_Q \bv\otimes\bv\mathbin{:}\bD\bw
          \\ &
          + \int_Q \bv\cdot\pp{\bw}{t}
          + \int_\Omega \bv_0\cdot\bw(0,\cdot)
        \end{aligned}
        \\
        \text{for all } \bw\in \D\bigl([0,T);\Wzdiv{q}\bigr),\,
        q=\max\bigl\{r,\tfrac{5r}{5r-6}\bigr\}.
      \end{multline*}
\end{Rem}

\begin{proof}[Proof of Theorem~\ref{th:noslip}]
  We shall distinguish two cases (that can be also identified via behavior
  of the total dissipation of energy with respect to scaling invariance of
  the governing equations, see~\cite{handbook2005}):
  the subcritical/critical case $r\geq\frac{11}{5}$
  and the supercritical case $r\in\bigl(\frac65,\frac{11}{5}\bigr)$.
  The problem can be analyzed in an~arbitrary spatial dimension $d$;
  then the supercritical case corresponds to
  $r\in\bigl(\frac{2d}{d+2},1+\frac{2d}{d+2}\bigr)$
  and the subcritical/critical case to $r\geq1+\frac{2d}{d+2}$.
  Note that the case $r=2$ (including the Euler/Navier-Stokes fluid)
  belongs to the supercritical case in any spatial dimension $d>2$.

  \noindent
  \textbf{The case $r\geq\frac{11}{5}$.}
  \textbf{Step 1. Galerkin approximations.}
  We first construct a~finite-dimensional approximation to the problem by
  the Galerkin method. To proceed, we consider an~auxiliary eigenvalue
  problem to find $\lambda\in\RR$ and $\bom\in\Vzdiv
  \hookrightarrow W^{1,\infty}(\Omega)^3$ satisfying
  \begin{align}
    \label{eq:p1}
    \llparen\bom,\bfi\rrparen = \lambda (\bom,\bfi)
    \text{ for all }\bfi\in\Vzdiv,
  \end{align}
  where $(\cdot,\cdot)$ is a~scalar product in $\bL^2$ and
  $\llparen\cdot,\cdot\rrparen$ is a~scalar product in
  $\Vzdiv$, i.e.,
  $\llparen\bom,\bfi\rrparen\coloneqq
    (\nabla^3\bom,\nabla^3\bfi) + (\bom,\bfi)
  $.
  It is known, see for example~\cite[Appendix A.4]{MNRR96}, that there
  exist eigenvalues $\left\{\lambda_m\right\}_{m=1}^\infty$
  and corresponding eigenfunctions $\left\{\bom^m\right\}_{m=1}^\infty$
  for the eigenvalue problem~\eqref{eq:p1}
  such that $0<\lambda_1\leq\lambda_2\leq\ldots$, $\lambda_m\rightarrow\infty$
  as $m\rightarrow\infty$, $(\bom^m,\bom^n)=\delta_{mn}$,
  $\llparen\frac{\bom^m}{\sqrt{\lambda_m}},
           \frac{\bom^n}{\sqrt{\lambda_n}}\rrparen=\delta_{mn}$.
  Furthermore, the mappings
  $\vec{P}^N:\Vzdiv\rightarrow H^N\coloneqq\lin\{\bom^1,
    \bom^2,\ldots,\bom^N\}$
  defined by $\vec{P}^N\bv\coloneqq\sum_{i=1}^N(\bv,\bom^i)\,\bom^i$
  are continuous orthonormal projectors in $\bL^2$,
  $\Vzdiv$ and $\bigl(\Vzdiv\bigr)^*$,
  in particular
  \begin{gather}
    \label{eq:p2}
    \|\vec{P}^N\|_{\mathcal{L}(\bL^2)} \leq 1, \qquad
    \|\vec{P}^N\|_{\mathcal{L}\bigl(\Vzdiv\bigr)} \leq 1, \qquad
    \|\vec{P}^N\|_{\mathcal{L}\bigl(\bigl(\Vzdiv\bigr)^*\bigr)} \leq 1.
  \end{gather}

  Galerkin approximations $\bv^N(t)\in H^N$ of the form
  $\bv^N(t,x)=\sum_{j=1}^N c_j^N(t)\,\bom^j(x)$ are introduced
  in such a~way that the coefficients $\bc^N=(c_1^N,c_2^N,\ldots,c_N^N)$
  fulfill
  \begin{subequations}
  \begin{align}
    \label{eq:p3a}
    \Bigl(\frac{\partial\bv^N}{\partial t},\bom^j\Bigr)
    -\left(\bv^N\otimes\bv^N,\nabla\bom^j\right)
    +\left(\bS(\bD\bv^N),\bD\bom^j\right)
    &= \bigl<\bb,\bom^j\bigr>
    \quad j=1,2,\ldots,N,
    \\
    \bv^N(0,\cdot) &= \vec{P}^N\bv_0,
    \label{eq:p3}
  \end{align}
  \end{subequations}
  where
  \begin{gather*}
    \bS(\bD) \coloneqq 2\nu_* \left( |\bD|-\delta_* \right)^+
                       \mathcal{S}(|\bD|) \tfrac{\bD}{|\bD|}.
  \end{gather*}
  Since the mappings
  $\bz\mapsto\bz\otimes\bz$ and $\bz\mapsto\bS(\bD\bz)$ are continuous,
  the Carath\'eodory theory for systems of ordinary differential equations
  implies local existence of a~solution $\bc^N$ solving~\eqref{eq:p3}.
  Global existence then follows from the fact that
  \begin{align*}
    \sup_{t\in(0,T)}|\bc^N(t)|_{\RR^N}<\infty.
  \end{align*}
  This piece of information is
  a~simple consequence of the orthogonality of the basis
  $\{\bom^j\}_{j=1}^\infty$ and a~priori estimates that will follow,
  see~\eqref{eq:p55} below.

  \noindent
  \textbf{Step 2. Uniform estimates and their consequences.}
  Multiplying~\eqref{eq:p3} by $c_j^N(t)$, taking the sum over
  $j=1,2,\ldots,N$, using the fact $(\bz\otimes\bz,\nabla\bz)=0$
  for $\bz$ with $\diver\bz=0$, $\bz\cdot\bn=0$ on $\partial\Omega$,
  we obtain
  \begin{align}
    \label{eq:p4}
    \frac12\frac{\d}{\d t}\|\bv^N\|_2^2 + \bigl(\bS(\bD\bv^N),\bD\bv^N\bigr)
    = \bigl<\bb,\bv^N\bigr>.
  \end{align}
  Since $\cG$ is an $r$-graph fulfilling~($\cG$\ref{Gcoercive}),
  we conclude that for all $t\in(0,T]$
  \begin{align}
    \label{eq:p5}
    \|\bv^N(t)\|_2^2
    + \alpha \int_0^t \left( \|\bS^N\|_{r'}^{r'} + \|\bD\bv^N\|_r^r \right)
    \leq C \left(\beta,
                 \|\bv_0\|_2^2,
                 \|\bb\|_{\left(L^{r}(0,T;\Wzdiv{r})\right)^*}
           \right).
  \end{align}
  Using the orthogonality of $\{\bom^j\}_{j=1}^N$ in $\bL^2$,
  this, in particular, implies that
  \begin{align}
    \label{eq:p55}
    \sup_{t\in(0,T)}|\bc^N(t)|_{\RR^N}<\infty
  \end{align}
  so that the proof of global-in-time existence of
  $\bc^N:[0,T]\rightarrow\RR^N$ is complete.

  Furthermore, Korn's inequality, see~\eqref{eq:Wzequivnorms},
  the Gagliardo-Nirenberg interpolation\footnote{%
    For $r\in\Bigl[\tfrac65,3\Bigr)$ the inequality~\eqref{eq:interp} follows
    immediately from the interpolation inequality
    $\|\bu\|_\frac{5r}{3}\leq\|\bu\|_2^\frac{2}{5}\|\bu\|_\frac{3r}{3-r}^\frac{3}{5}$
    and the Gagliardo-Nirenberg embedding.
    For arbitrary $r\geq2$, consider the following interpolation inequality
    \begin{align*}
      \|\bu\|_\frac{5r}{3}
      \leq C_r \|\bu\|_r^\frac{5r-6}{5r} \|\nabla\bu\|_r^\frac{6}{5r}
      \quad \text{ for all } \bu\in\Wz{r},\; r\geq\tfrac65;
    \end{align*}
    see \cite[Theorem~5.8]{adams-fournier-2003}.
    Using the H\"older interpolation
    $\|\bu\|_r\leq\|\bu\|_2^\frac{4}{5r-6}\|\bu\|_\frac{5r}{3}^\frac{5r-10}{5r-6}$
    we obtain
    \begin{align*}
      \|\bu\|_\frac{5r}{3}
      \leq C_r \|\bu\|_2^\frac{4}{5r}
               \|\bu\|_\frac{5r}{3}^{1-\frac{2}{r}}
               \|\nabla\bu\|_r^\frac{6}{5r},
    \end{align*}
    which yields the desired estimate~\eqref{eq:interp} after rearrangement.
  }
  \begin{align}
    \label{eq:interp}
    \|\bu\|_\frac{5r}{3} \leq C_r \|\bu\|_2^\frac{2}{5} \|\nabla\bu\|_r^\frac{3}{5}
    \quad \text{ for all } \bu\in\Wz{r},\; r\geq\tfrac65,
  \end{align}
  together with~\eqref{eq:p5} imply
  \begin{multline*}
    \int_0^T \|\bv^N\|_\frac{5r}{3}^\frac{5r}{3}
      {}\leq
      C \esssup_{t\in(0,T)} \|\bv^N(t)\|_2^{\frac{2r}{3}}
      \int_0^T \|\bD\bv^N\|_r^r
    \\
    \leq C \left(\beta,
                 \|\bv_0\|_2^2,
                 \|\bb\|_{\left(L^{r}(0,T;\Wzdiv{r})\right)^*}
           \right).
  \end{multline*}

  Finally, since for all $\bfi\in L^s(0,T;\Vzdiv)$
  \begin{align*}
    \int_0^T \left( \frac{\partial\bv^N}{\partial t},\bfi \right)
    = \int_0^T \left( \frac{\partial\bv^N}{\partial t},\bP^N\bfi \right),
  \end{align*}
  it follows from~\eqref{eq:p3}, \eqref{eq:p5}, the fact that
  $2r'=\frac{2r}{r-1}\leq\frac{5r}{3}\Leftrightarrow r\geq\frac{11}{5}$,
  and~\eqref{eq:p2}, that
  \begin{align*}
    \left\lVert \frac{\partial\bv^N}{\partial t} \right\rVert
      _{\Bigl(L^{r\fixsupscriptmore}\bigl(0,T;\Vzdiv\bigr)\Bigr)^*}
    &\coloneqq \sup_{\substack{
      \bfi\in L^{r\fixsupscriptmore}\bigl(0,T;\Vzdiv\bigr) \\
      \|\bfi\|_{L^{r\fixsupscript}\bigl(0,T;\Vzdiv\bigr)} \leq1
    }}
    \left| \int_0^T \left( \frac{\partial\bv^N}{\partial t},\bfi \right) \right|
    \\ &\leq C \left(\beta,
                     \|\bv_0\|_2^2,
                     \|\bb\|_{\Bigl(L^{r\fixsupscriptmore}\bigl(\Wzdiv{r}\bigr)\Bigr)^*}
               \right).
  \end{align*}
  Consequently, there are (not relabeled) subsequences so that
  \begin{align}
    \nonumber
    \bv^N &\stackrel{*}{\rightharpoonup} \bv
    && *\text{-weakly in } L^\infty(0,T;\bL^2) ,\\
    \nonumber
    \bD\bv^N &\rightharpoonup \bD\bv
    && \text{weakly in } L^r(0,T;L^r(\Omega)^{3\times3}_\mathrm{sym}) ,\\
    \nonumber
    \nabla\bv^N &\rightharpoonup \nabla\bv
    && \text{weakly in } L^r(0,T;L^r(\Omega)^{3\times3}) ,\\
    \nonumber
    \bS^N &\rightharpoonup \bS
    && \text{weakly in } L^{r'}(0,T;L^{r'}(\Omega)^{3\times3}_\mathrm{sym}) ,\\
    \label{eq:p8}
    \partial_t\bv^N &\rightharpoonup \partial_t\bv
    && \text{weakly in } \bigl( L^r(0,T;\Vzdiv) \bigr)^* ,\\
    \label{eq:p9}
    \bv^N &\rightarrow \bv
    && \text{strongly in } L^q(0,T;\bL^q)
    \text{ for all } q\in\bigl[1,\tfrac{5r}{3}\bigr),
  \end{align}
  \begingroup\setlength\emergencystretch{\hsize}\hbadness=10000
  where the last limit~\eqref{eq:p9} follows from the Aubin-Lions
  compactness lemma applied to
  $\Vzdiv \hookrightarrow \Wzdiv{r}
  \hookrightarrow\hookrightarrow \Lndiv{r} \hookrightarrow
  \bigl( \Vzdiv \bigr)^*$.
  \par\endgroup

  Finally, letting $N\rightarrow\infty$ in~\eqref{eq:p3} for $j\in\mathbb{N}$
  arbitrary but fixed, one concludes that $(\bv,\bS)$ satisfy
  \begin{align*}
    \int_0^T \left(
      \left<\frac{\partial\bv}{\partial t},\bom\right>
      -\left(\bv\otimes\bv,\nabla\bom\right)
      +\left(\bS,\bD\bom\right)
    \right) \phi(t) \d t
    =  \int_0^T \left<\bb,\bom\right>\, \phi(t) \d t
  \end{align*}
  valid for all $\phi\in \D(-\infty,\infty)$ and
  $\bom\in\Vzdiv$. Since the space $\Vzdiv$
  is dense in $\Wzdiv{r}$ for all $r\geq1$,
  and functions of the form
  $\phi(t)\bom(x)$ are dense in $L^r(0,T; \Wzdiv{r})$,
  and finally, for $r\geq\frac{11}{5}$,
  \begin{align*}
    \left| \int_0^T \left(\bv\otimes\bv,\nabla\psi\right) \right|
    &\leq \left( \int_0^T \|\bv\|_{2r'}^{2r'} \right)^\frac{1}{r'}
          \left( \int_0^T \|\nabla\psi\|_{r}^{r} \right) \\
    &\leq C \left( \int_0^T \|\bv\|_\frac{5r}{3}^\frac{5r}{3} \right)^\frac{1}{r'}
            \left( \int_0^T \|\nabla\psi\|_{r}^{r} \right)
    < +\infty,
  \end{align*}
  we deduce that $(\bv,\bS)$ satisfies
  \begin{multline}
    \label{eq:p10}
    \int_0^T \left< \partial_t\bv,\bom\right>
    +\int_Q \bS\mathbin{:}\bD\bom
    = \int_0^T \left<\bb,\bom\right>
    + \int_Q\left(\bv\otimes\bv\right)\mathbin{:}\bD\bom
    \\
    \text{ for all }\bom\in L^r(0,T;\Wzdiv{r}).
  \end{multline}
  This implies that
  $\partial_t\bv\in L^{r'}\bigl(0,T;(\Wzdiv{r})^*\bigr)$.
  The embedding
  \begin{align}
    \label{eq:p11}
    \Bigl\{\bu\in L^r(0,T;\Wzdiv{r});\,\allowbreak
    \partial_t\bu\in L^{r'}\bigl(0,T;(\Wzdiv{r})^*\bigr)\Bigr\}
    \hookrightarrow \mathcal{C}([0,T];\bL^2)
  \end{align}
  proves~\eqref{eq:contstrong}.
  Inserting $\bom\coloneqq\bv\chi_{(0,t)}$ into~\eqref{eq:p10},
  we obtain the energy equality~\eqref{eq:P6}
  with $\bv(0)$ instead of~$\bv_0$. Once we show that
  $\bv(0)=\bv_0$ we obtain \eqref{eq:P3}, \eqref{eq:P6}
  and infer~\eqref{eq:P4} from~\eqref{eq:p10}. In order to show that
  $\bv(0)=\bv_0$ we first observe that for any $\phi\in\D([0,T))$ it holds
  \begin{align}
    \label{eq:p110}
    \int_0^T \bigl(\partial_t\bv^N,\bom^j\bigr)\,\phi =
    -\int_0^T \bigl(\bv^N,\bom^j\bigr)\,\phi' -(\bv^N(0),\bom^j)\,\phi(0)
    \quad \text{for } j=1,2,\ldots,N.
  \end{align}
  By virtue of~\eqref{eq:p3} and the orthonormality of~$\vec{P}^N$,
  the equality
  $(\bv^N(0),\bom^j)=(\vec{P}^N\bv_0,\bom_j)\allowbreak{}=(\bv_0,\bom_j)$
  holds
  for $j=1,2,\ldots,N$. Hence the right-hand side of~\eqref{eq:p110}
  with $N\rightarrow\infty$ tends to
  \begin{align}
    \label{eq:p111}
    -\int_0^T \bigl(\bv,\bom^j\bigr)\,\phi' -(\bv_0,\bom^j)\,\phi(0).
  \end{align}
  On the other hand, by virtue of~\eqref{eq:p8}, the left-hand side of~\eqref{eq:p110}
  converges to
  \begin{align}
    \label{eq:p112}
    \int_0^T \bigl<\partial_t\bv,\bom^j\bigr>\,\phi =
    -\int_0^T \bigl(\bv,\bom^j\bigr)\,\phi' -(\bv(0),\bom^j)\,\phi(0),
  \end{align}
  where the last equality is justified by the inclusions
  $\bom^j\in \Vzdiv \hookrightarrow \Wzdiv{r}$
  and $\partial_t\bv\in L^{r'}\bigl(0,T;(\Wzdiv{r})^*\bigr)$
  together with~\eqref{eq:p11}.
  As $\phi(0)$ is arbitrary, comparing~\eqref{eq:p111}
  and~\eqref{eq:p112} yields
  \begin{align*}
    (\bv_0,\bom^j) = (\bv(0),\bom^j)
    \qquad j=1,2,\ldots.
  \end{align*}
  As $\{\bom_j\}_{j=1}^\infty$ is basis in~$\bL^2$, we conclude
  that $\bv(0)=\bv_0$.

  It remains to show~\eqref{eq:P5}.

  \noindent
  \textbf{Step 3. Attainment of the constitutive equation.}
  To prove~\eqref{eq:P5}, we wish to use the graph convergence lemma
  (see Lemma~\ref{lemma:conv} in Appendix).
  To apply this lemma, we need to show that
  \begin{align*}
    \limsup_{N\rightarrow\infty} \int_Q \bS^N\mathbin{:}\bD\bv^N \leq \int_Q \bS\mathbin{:}\bD\bv.
  \end{align*}
  However~\eqref{eq:p9} implies, in particular, that
  \begin{align*}
    \bv^N(t) \rightarrow \bv(t)
    \quad\text{ in } \bL^2 \text{ for almost all } t\in[0,T].
  \end{align*}
  Integrating~\eqref{eq:p4} from $0$ to such $t$'s, and letting
  $N\rightarrow\infty$, one concludes
  \begin{align*}
    \tfrac12\|\bv(t)\|_2^2
    +\lim_{N\rightarrow\infty} \int_0^t \int_\Omega \bS^N\mathbin{:}\bD\bv^N
    = \int_0^t \left<\bb,\bv\right> + \tfrac12 \|\bv_0\|_2^2.
  \end{align*}
  By comparing this identity with~\eqref{eq:P6}
  (which is an equality as $r\geq\frac{11}{5}$), we conclude
  \begin{align*}
    \lim_{N\rightarrow\infty}
    \int\limits_{\mathclap{(0,t)\times\Omega}} \bS^N\mathbin{:}\bD\bv^N
    =
    \int\limits_{\mathclap{(0,t)\times\Omega}} \bS\mathbin{:}\bD\bv
    \qquad \text{for a.a. } t\in[0,T].
  \end{align*}
  The graph convergence lemma (Lemma~\ref{lemma:conv}) then implies that
  $\bS$ and $\bD\bv$ fulfill~\eqref{eq:P5}. The proof for
  $r\geq\frac{11}{5}$ is thus complete.

  \noindent
  \textbf{The case $r\in\left(\frac65,\frac{11}{5}\right)$.}
  \textbf{Step 1. Approximations and their validity.}
  For $\epsilon>0$, we look for $(\bv^\epsilon,\bS^\epsilon)$ such that
  \begin{subequations}
  \label{eq:p121314}
  \begin{align}
    \bv^\epsilon &\in L^\infty\bigl(0,T;\Lndiv2\bigr) \cap L^\frac{11}{5}\bigl(0,T;\Wzdiv{\frac{11}{5}}\bigr),
    \\
    \partial_t\bv^\epsilon &\in L^\frac{11}{5}\bigl(0,T;(\Wzdiv{\frac{11}{5}})^*\bigr),
    \\
    \bS^\epsilon &\in L^{r'}(Q)^{3\times 3}_\mathrm{sym}
  \end{align}
  \end{subequations}
  satisfy
  \begin{gather}
    \label{eq:p15}
    \begin{multlined}
      \int_Q \bS^\epsilon\mathbin{:}\bD\bfi
      + \epsilon \int_Q |\bD\bv^\epsilon|^\frac15\, \bD\bv^\epsilon\mathbin{:}\bD\bfi
      \\
      = \int_0^T \left<\bb,\bfi\right>
      + \int_Q \left(\bv^\epsilon\otimes\bv^\epsilon\right)\mathbin{:}\bD\bfi
      + \int_Q \bv^\epsilon\cdot\pp{\bfi}{t}
      + \int_\Omega \bv_0\cdot\bfi(0,\cdot)
      \\
      \text{for all }
      \bfi\in \D\bigl([0,T);\Wzdiv{\frac{11}{5}}\bigr),
    \end{multlined}
    \\
    \label{eq:p16}
    \bS^\epsilon = 2\nu_* \left( |\bD\bv^\epsilon|-\delta_* \right)^+
            \mathcal{S}(|\bD\bv^\epsilon|) \tfrac{\bD\bv^\epsilon}{|\bD\bv^\epsilon|}
    \quad\text{ almost everywhere in }Q,
    \\
    \shortintertext{and}
    \label{eq:p17}
    \bv^\epsilon\in \mathcal{C}\bigl([0,T]; \bL^2\bigr),
    \qquad
    \bv^\epsilon(0) = \bv_0.
  \end{gather}

  The existence of $(\bv^\epsilon,\bS^\epsilon)$ fulfilling \eqref{eq:p121314}--\eqref{eq:p17}
  for arbitrary but fixed $\epsilon>0$ can be proved in the same way as
  the existence of a~weak solution to the problem
  for the case $r\geq\frac{11}{5}$. In addition, by taking
  $\bfi\coloneqq\bv^\epsilon\chi^n$ in~\eqref{eq:p15}
  with a~sequence $\chi^n\in\D([0,T))$, $0\leq\chi^n\leq1$,
  $\chi^n\rightarrow\chi_{(0,t)}$, $0\leq t\leq T$
  and letting $n\rightarrow\infty$,
  using Lebesgue's dominated convergence theorem, we obtain
  \begin{align}
    \label{eq:p21}
    \tfrac12 \|\bv^\epsilon(t)\|_2^2 -\tfrac12 \|\bv_0\|_2^2
    +\epsilon \int_0^t \|\bD\bv^\epsilon\|_\frac{11}{5}^\frac{11}{5}
    +\int_0^t \int_\Omega \bS^\epsilon\mathbin{:}\bD\bv^\epsilon
    = \int_0^t \left<\bb,\bv^\epsilon\right>
  \end{align}
  for all $t\in[0,T]$
  where $\bD\bv^\epsilon$ and $\bS^\epsilon$ satisfy~\eqref{eq:p16}.

  \noindent
  \textbf{Step 2. Estimates uniform with respect to $\epsilon$ and their consequences.}
  Since
  ${\bS\mathbin{:}\bD}\geq{\alpha\bigl(|\bD|^r+|\bS|^{r'}\bigr)-\beta}$
  for all
  $(\bS,\bD)$ fulfilling~\eqref{eq:p16}, i.e., $(\bS,\bD)\in\cG$, the
  energy identity~\eqref{eq:p21} implies that
  \begin{subequations}
  \label{eq:p2223242526}
  \begin{align}
    \label{eq:p22}
    \{\bv^\epsilon;\epsilon>0\}
    &\text{ is bounded in } L^\infty\bigl(0,T;\bL^2\bigr) ,\\
    \{\bv^\epsilon;\epsilon>0\}
    &\text{ is bounded in } L^r\bigl(0,T;\Wzdiv{r}\bigr) ,\\
    \{\bD\bv^\epsilon;\epsilon>0\}
    &\text{ is bounded in } L^r\bigl(0,T;L^r(\Omega)^{3\times 3}_\mathrm{sym}\bigr) ,\\
    \label{eq:p25}
    \{\epsilon^\frac{5}{11}\bD\bv^\epsilon;\epsilon>0\}
    &\text{ is bounded in } L^\frac{11}{5}\bigl(0,T;\bL^\frac{11}{5}\bigr) ,\\
    \{\bS^\epsilon;\epsilon>0\}
    &\text{ is bounded in } L^{r'}\bigl(0,T;L^{r'}(\Omega)^{3\times 3}_\mathrm{sym}\bigr) .
  \end{align}
  \end{subequations}
  Using the fact that
  \begin{align*}
    \int_0^T \left<\partial_t\bv^\epsilon,\bfi\right>
    = -\int_Q \bv^\epsilon\cdot\partial_t\bfi
    -\int_\Omega \bv_0
    \cdot\bfi(0,\cdot)
  \end{align*}
  for all $\bfi\in\D\bigl([0,T);\Wzdiv{\frac{11}{5}}\bigr)$
  we conclude from~\eqref{eq:p15} and
  \eqref{eq:p2223242526} that
  \begin{align*}
    \left\{\partial_t\bv^\epsilon;\epsilon>0\right\}
    \text{ is bounded in }
    \Bigl(L^\frac{5r}{5r-6}\bigl(0,T;\Wz{\frac{5r}{5r-6}}(\Omega)\bigr)\Bigr)^*.
  \end{align*}
  These estimates together with Korn's inequality~\eqref{eq:Wzequivnorms}
  and the Aubin-Lions lemma
  lead to the existence of $\bv$ and $\bS$ such that for suitable
  sequences $(\bv^m,\bS^m)\coloneqq(\bv^{\epsilon_m},\bS^{\epsilon_m})$
  and $m\rightarrow\infty$ the following convergences hold:
  \begin{align}
    \nonumber
    \bv^m &\stackrel{*}{\rightharpoonup} \bv
    & *\text{-weakly in }& L^\infty\bigl(0,T;\bL^2\bigr) ,\\
    \nonumber
    \bv^m &\rightharpoonup \bv
    & \text{weakly in }& L^r\bigl(0,T;\Wzdiv{r}\bigr) ,\\
    \label{eq:p35}
    \bD\bv^m &\rightharpoonup \bD\bv
    & \text{weakly in }& L^r\bigl(0,T;L^r(\Omega)^{3\times3}_\mathrm{sym}\bigr) ,\\
    \label{eq:p36}
    \bS^m &\rightharpoonup \bS
    & \text{weakly in }& L^{r'}\bigl(0,T;L^{r'}(\Omega)^{3\times3}_\mathrm{sym}\bigr) ,\\
    \nonumber
    \bv^m &\rightarrow \bv
    & \text{strongly in }& L^q\bigl(0,T;\bL^q\bigr)
    \text{ for all } q\in\bigl[1,\tfrac{5r}{3}\bigr).
    \intertext{
      As $r> \tfrac65$, it is $\tfrac{5r}{3}>2$. Consequently
    }
    \nonumber
    \bv^m &\rightarrow \bv
    & \text{strongly in }& L^2\bigl(0,T;{\bL}^2\bigr), \\
    \label{newP1}
    \bv^m(t) &\rightarrow \bv(t)
    & \text{strongly in }& \bL^2 \text{ for almost all } t\in[0,T].
  \end{align}
  Also \eqref{eq:p25} implies that
  \begin{align}
    \label{eq:p28}
    \epsilon_m \int_Q |\bD\bv^m|^\frac15\; |\bD\bv^m\mathbin{:}\bD\bfi|
    \leq \epsilon_m^\frac{5}{11}
    \left( \epsilon_m \int_Q |\bD\bv^m|^\frac{11}{5} \right)^\frac{6}{11}
    \|\bD\bfi\|_{\frac{11}{5},Q}
    \xrightarrow{m\rightarrow\infty}0.
  \end{align}
  Let us consider~\eqref{eq:p15} with
  $\bfi\in \D\bigl([0,T);\Wzdiv{q}\bigr)$, $q=\frac{5r}{5r-6}$,
  and let $m\rightarrow\infty$. Then we easily arrive to~\eqref{eq:P4}.

  \noindent
  \textbf{Step 3. Attainment of the constitutive equation~\eqref{eq:P5}.}
  In order to apply the graph convergence lemma (Lemma~\ref{lemma:conv}) we need to proceed
  in a~more subtle way as $\bw\coloneqq\bv$ is not an~admissible test
  function in~\eqref{eq:P4}. For $\bu^m\coloneqq\bv^m-\bv$, the following
  identity holds
  \begin{align}
    \label{eq:p29}
    \begin{multlined}
      -\int_Q (\bv^m-\bv)\cdot\partial_t\bw
      +\int_Q (\bS^m-\bS)\mathbin{:}\bD\bw
      +\epsilon_m \int_Q |\bD\bv^m|^\frac15\;\bD\bv^m\mathbin{:}\bD\bom
      \\
      =\int_Q (\bv^m\otimes\bv^m-\bv\otimes\bv)\mathbin{:}\bD\bw
      \\
      \text{for all } \bw\in \Cinf([0,T]; \Cinfzdiv)
    \end{multlined}
  \end{align}
  and
  \begin{align*}
    \bu^m &\stackrel{*}{\rightharpoonup} \b0
    && *\text{-weakly in } L^\infty(0,T;\Lndiv2) ,\\
    \bu^m &\rightharpoonup \b0
    && \text{weakly in } L^r(0,T;\Wzdiv{r}) ,\\
    \bu^m &\rightarrow \b0
    && \text{strongly in } L^q(0,T;\bL^q)
    \text{ for all } q\in\bigl[1,\tfrac{5r}{3}\bigr) .
  \end{align*}

  We also observe that besides~\eqref{eq:p28}
  \begin{align*}
    (\bS^m-\bS) \eqqcolon \bH_1^m
    &\rightharpoonup \bO
    &&\text{ weakly in } L^{r'}(Q)^{3\times 3}_{\mathrm{sym}} ,\\
    \Biggl(
    \begin{aligned}
      &\epsilon_m |\bD\bv^m|^\frac15\;\bD\bv^m+ \\
      &+(\bv\otimes\bv-\bv^m\otimes\bv^m)
    \end{aligned}
    \Biggr) \eqqcolon \bH_2^m
    &\rightarrow \bO
    &&\text{ strongly in } L^\sigma(Q) ^{3\times 3}_{\mathrm{sym}}
    \text{ for some } \sigma\in\bigl(1,\tfrac{5r}{6}\bigr)
  \end{align*}
  and we rewrite~\eqref{eq:p29} as
  \begin{align*}
    \int_Q \bu^m\cdot\partial_t\bw
    =\int_Q \bH_1^m\mathbin{:}\bD\bw
    +\int_Q \bH_2^m\mathbin{:}\bD\bw.
  \end{align*}
  Let $Q_0=I_0\times B_0\subset Q$ be any space-time cylinder.
  Then, applying the Lipschitz truncation Lemma~\ref{lemma:trunc},
  property~(\ref{lemma:trunc:conseq}), we conclude,
  using the above convergences, that
  \begin{align}
    \label{eq:p100}
    \limsup_{m\rightarrow\infty} \left|
      \int_{\frac{1}{6}Q_0\setminus\mathcal{O}^{m,k}}
      (\bS^m-\bS)\mathbin{:}(\bD\bv^m-\bD\bv)\,\xi
    \right|
    \leq C\, 2^{-\frac{k}{r}}
  \end{align}
  with some $\xi\in\D\bigl(\frac{Q_0}{6}\bigr)$ such that
  $\chi_{\frac18Q_0} \leq \xi \leq \chi_{\frac16Q_0}$.
  \begingroup\setlength\emergencystretch{\hsize}\hbadness=10000
  Pick $\bS^*$ so that $(\bS^*,\bD\bv)\in\cG$. By virtue
  of definition~\eqref{S.3*} such $\bS^*$ is defined uniquely and
  fulfills $\bS^*\in L^{r'}(0,T;L^{r'}(\Omega)^{3\times 3}_\mathrm{sym})$.
  Hence, \eqref{eq:p100} and \eqref{eq:p35} imply
  \begin{align}
    \label{eq:p101}
    \limsup_{m\rightarrow\infty} \left|
      \int_{\frac{1}{6}Q_0\setminus\mathcal{O}^{m,k}}
      (\bS^m-\bS^*)\mathbin{:}(\bD\bv^m-\bD\bv)\,\xi
    \right|
    \leq C\, 2^{-\frac{k}{r}}.
  \end{align}
  \endgroup
  As both $(\bS^*,\bD\bv)\in\cG$ and $(\bS^m,\bD\bv^m)\in\cG$,
  the monotonicity property~($\cG$\ref{Gmonotone}) of~$\cG$
  implies that the integrand
  in~\eqref{eq:p101} is non-negative.
  This, together with properties (\ref{lemma:trunc:lambda}),
  (\ref{lemma:trunc:smallness})
  of the truncation lemma (Lemma~\ref{lemma:trunc}) and H\"older's inequality,
  implies that
  \begin{align*}
    \limsup_{m\rightarrow\infty}
        \int_{\frac18Q_0} \bigl| (\bS^m-\bS^*)\mathbin{:}(\bD\bv^m-\bD\bv) \bigr|^\frac12
    \leq C\, 2^{-\frac{k}{r}}.
  \end{align*}
  Set $g^m\coloneqq(\bS^m-\bS^*)\mathbin{:}(\bD\bv^m-\bD\bv)$.
  Clearly $g^m\geq0$ and
  $g^m\rightarrow0$ almost everywhere in $\frac18Q_0$. But as $Q_0$ is
  arbitrary, we conclude that
  \begin{align}
    \label{eq:p34}
    g^m\rightarrow0 \text{ almost everywhere in } Q.
  \end{align}
  Since $\{g^m\}_{m=1}^\infty$ is bounded in $L^1(Q)$ and has pointwise
  limit~\eqref{eq:p34}, Corollary~\ref{col:strongbiting}, a~consequence
  of the biting lemma (Lemma~\ref{lemma:biting}), ensures the
  existence of a~subsequence $\{g^{m_j}\}_{j=1}^\infty$
  and a~sequence of sets $\{E_k\}_{k=1}^\infty$ with
  $Q\supset E_1\supset E_2\supset\ldots$, $|E_k|\rightarrow0$ such that for
  all $k\in\mathbb{N}$
  \begin{align*}
    g^{m_j} \rightarrow 0 \text{ strongly in } L^1(Q\setminus E_k).
  \end{align*}
  From the definition of $g^{m_j}$ we conclude that
  \begin{align*}
    \limsup_{j\rightarrow0}
     \int_{Q\setminus E_k} (\bS^{m_j}-\bS^*)\mathbin{:}(\bD\bv^{m_j}-\bD\bv) = 0,
  \end{align*}
  which implies with the help of~\eqref{eq:p35} and \eqref{eq:p36}
  that for all $k\in\mathbb{N}$
  \begin{align*}
    \limsup_{j\rightarrow0} \int_{Q\setminus E_k} \bS^{m_j}\mathbin{:}\bD\bv^{m_j}
     = \int_{Q\setminus E_k} \bS\mathbin{:}\bD\bv.
  \end{align*}
  Since $|E_k|\rightarrow0$, we can conclude from the graph convergence
  lemma (Lemma~\ref{lemma:conv}) that
  \begin{align*}
    (\bS,\bD\bv)\in\cG\text{ almost everywhere in }Q
  \end{align*}
  so that \eqref{eq:P5} holds.

  \noindent
  \textbf{Step 4. The energy inequality and the initial condition.}
  Since $(\bS^*,\bD\bv)\in\cG$ and $\cG$ is monotone, we first observe that
  $g^m=(\bS^m-\bS^*)\mathbin{:}(\bD\bv^m-\bD\bv)\geq0$.
  It thus follows from~\eqref{eq:p34},
  Fatou's lemma applied to $\{g^m\}_{m=1}^\infty$,
  and from the weak convergences
  \eqref{eq:p35} and \eqref{eq:p36} that
  \begin{align*}
    \mathop{\int_0^t} \int_\Omega \bS\mathbin{:}\bD\bv
    \leq
    \liminf_{m\rightarrow\infty}
    \mathop{\int_0^t} \int_\Omega \bS^m\mathbin{:}\bD\bv^m.
  \end{align*}
  It is then easy to conclude the energy inequality~\eqref{eq:P6}
  from~\eqref{eq:p21} and~\eqref{newP1} for almost all $t\in (0,T]$.

  Now we prove ~\eqref{eq:contweak}. Taking $\bw = \bfi\chi^n$ in~\eqref{eq:P4} with
  $\chi^n\in\D((0,T))$, $0\leq\chi^n\leq1$,
  $\chi^n\rightarrow\chi_{(t_0,t)}$,
  $0<t_0\leq t\leq T$ and $\bfi\in\Wzdiv{q}$, and letting then $n\rightarrow\infty$,
  using Lebesgue's dominated convergence theorem, we obtain
  \begin{align*}
    \bigl(\bv(t)-\bv(t_0),\bfi\bigr)
    = \int_{t_0}^t (\bv\otimes\bv,\bD\bfi) + \langle\bb,\bfi\rangle - (\bS,\bD\bfi).
  \end{align*}
  As the integrand (at the right-hand side) is integrable over $(0,T)$, then necessarily
  \begin{align*}
    \lim_{t\rightarrow t_0} \bigl(\bv(t)-\bv(t_0),\bfi\bigr) = 0
    \qquad \text{for all } \bfi\in\Wzdiv{q}.
  \end{align*}
  This together with the fact that $\bv\in L^\infty(0,T;\bL^2)$ implies  that
  $\bv\in \mathcal{C}\bigl([0,T]; \bL^2_\mathrm{weak}\bigr)$. Once we know this, we can extend
  validity of the energy inequality~\eqref{eq:P6} to all $t\in (0,T]$ by approximating the
  exceptional $t$'s by $t_n\to t$ as $n\to \infty$, whereas $t_n$'s are time instants for which~\eqref{eq:P6} holds.
  Consequently, \eqref{eq:P6} holds for all $t\in (0,T]$.

  By the same argument as above, this time with $\chi^n\in\D([0,T))$, $\chi^n\rightarrow\chi_{[0,t)}$,
  $0<t<T$, we also obtain $\bv(t) \rightharpoonup \bv_0$ weakly in  $\bL^2$ as $t\rightarrow0+$.
  This also implies $\liminf_{t\to0+}\|\bv(t)\|_2 \geq \|\bv_0\|_2$. But
  from~\eqref{eq:P6} we have
  $\limsup_{t\to0+}\|\bv(t)\|_2 \leq \|\bv_0\|_2$,
  hence $\|\bv(t)\|_2\to\|\bv_0\|_2$ as $t\to0+$.
  This together with the weak convergence implies
  that
  $\bv(t) \rightarrow \bv_0$ strongly in  $\bL^2$ as $t\rightarrow0+$
  and~\eqref{eq:P3} is thus proved.

  \noindent
  \textbf{On the pressure.}
  Let us consider for fixed $t\in(0,T)$ the functionals
  \begin{equation} \label{P98}
  \begin{aligned}
    \bigl<\bF^1(t),\bfi\bigr> &\coloneqq
      \int_\Omega \bigl( \bv(t,\cdot) - \bv_0 \bigr)\cdot\bfi, \\
    \bigl<\bF^2(t),\bfi\bigr> &\coloneqq
      \int_\Omega \begingroup\textstyle\int\endgroup_0^t \bigl( \bS-\bv\otimes\bv \bigr) \mathbin{:} \bD\bfi
      - \Bigl< \begingroup\textstyle\int\endgroup_0^t\bb, \bfi \Bigr>
  \end{aligned}
  \end{equation}
  for $\bfi\in \Wz{q}$ with $q\coloneqq\max\bigl\{r,\frac{5r}{5r-6}\bigr\}$.
  Clearly $\bF^1(t),\bF^2(t)\in \bigl(\Wz{q}\bigr)^*$ for almost every $t\in(0,T)$.
  Testing~\eqref{eq:P4} by $\bw^j\in \D\bigl([0,T);\Cinfzdiv\bigr)$
  such that $\bw^j\rightarrow\bw$ and
  \begin{align*}
    \bw(s, x) = \left\{\begin{array}{ll}
      \bfi(x) & s\in[0,t), \\
      \b0     & s\in[t,T)
    \end{array}\right.
  \end{align*}
  with arbitrary $\bfi\in \Cinfzdiv$ and comparing with~\eqref{P98}
  we obtain
  \begin{align} \label{P99}
    \Bigl<\bigl(\bF^1+\bF^2\bigr)(t),\bfi\Bigr> = 0
    \qquad\text{for all } \bfi\in\Cinfzdiv
    \text{ and a.e. } t\in(0,T).
  \end{align}
  Now consider the Stokes problems
  \begin{subequations} \label{P101}
  \begin{align}
    \label{P101a}
    - \Delta\bU^1 + \nabla P^1 &= \bF^1,
    & \diver\bU^1 &= 0
    \quad\text{in }Q,
    & \bU^1 &= \b0
    \quad\text{on }(0,T)\times\partial\Omega,
    \\
    - \Delta\bU^2 + \nabla P^2 &= \bF^2,
    & \diver\bU^2 &= 0
    \quad\text{in }Q,
    & \bU^2 &= \b0
    \quad\text{on }(0,T)\times\partial\Omega.
  \end{align}
  \end{subequations}
  By virtue of the assumptions on the domain, we conclude
  from~\eqref{eq:stokes_w1q} of
  Lemma~\ref{lemma:stokesregularity} that
  \begin{alignat*}{4}
    &\|\nabla\bU^1(t)\|_6
    &&{}+ \|P^1(t)\|_6
    &&\leq C \|\bv(t)-\bv_0\|_{-1,6}
    \leq C \|\bv(t)-\bv_0\|_2,
    \\
    &\|\nabla\bU^2(t)\|_{q'}
    &&{}+ \|P^2(t)\|_{q'}
    &&\leq C \|\bF^2(t)\|_{-1,q'}
    \\&&&
    &&\leq C \int_0^t \| \bS-\bv\otimes\bv \|_{q'} + C \int_0^t \| \bb \|_{-1,q'},
  \end{alignat*}
  which leads to $P^1\in L^\infty(0,T;L^6)$
  and $p^2\coloneqq\partial_t P^2\in L^{q'}(Q)$.
  Testing~\eqref{P101a} with $\nabla\phi$ for $\phi\in\D(\Omega)$ we get
  \begin{align*}
    \underbrace{(\nabla\bU^1(t),\nabla^2\phi)}_{(\diver\bU^1(t),\Delta\phi)=0}
    - (P^1(t),\Delta\phi)
    = (\bv(t)-\bv_0,\nabla\phi) = -(\diver(\bv(t)-\bv_0),\phi) = 0
  \end{align*}
  so that $(P^1(t),\Delta\phi)=0$ for all $\phi\in\D(\Omega)$ and Weyl's lemma
  (cf. \cite[Lemma 2]{Weyl1940}, \cite[Chapter 10]{Garding1997}) yields
  that $P^1(t)$ is harmonic.

  Testing~\eqref{P101} by $\bfi\in \Wzdiv{q}$
  we obtain, by using~\eqref{P99},
  \begin{align*}
    (\nabla(\bU^1(t)+\bU^2(t)),\nabla\bfi) = 0
    \qquad \text{ for all }\bfi\in \Wzdiv{q}
    \text{ and a.e. } t\in(0,T),
  \end{align*}
  which shows together with $\bU^1(t)+\bU^2(t) \in \Wzdiv{q'}$
  that $\bU^1+\bU^2=\b0$. Now we are in a~position to sum up~\eqref{P101},
  test by $\partial_t\bw$ with $\bw\in \D([0,T);\Cinfzdiv)$,
  integrate over~$Q$, and use the facts shown above and $P^2(0)=0$ to obtain~\eqref{eq:P5p}.

  Furthermore, when $\Omega$ is a~$C^{1,1}$ domain, \eqref{eq:stokes_w2q} from
  Lemma~\ref{lemma:stokesregularity} yields
  \begin{alignat*}{3}
    &\|\nabla^2\bU^1(t)\|_2
    &&{}+ \|\nabla P^1(t)\|_2
    &&\leq C \|\bv(t)-\bv_0\|_2,
    \\
    &\|\nabla^2\bU^1(t)\|_\frac{5r}{3}
    &&{}+ \|\nabla P^1(t)\|_\frac{5r}{3}
    &&\leq C \|\bv(t)-\bv_0\|_\frac{5r}{3}
  \end{alignat*}
  so that $\esssup_{t\in(0,T)}\|\nabla P^1(t)\|_2 \leq
  C\esssup_{t\in(0,T)}\|\bv(t)-\bv_0\|_2$ and
  $\int_0^T \|\nabla P^1\|_\frac{5r}{3}^\frac{5r}{3}
  \leq \int_0^T \|\bv-\bv_0\|_\frac{5r}{3}^\frac{5r}{3} \leq C$
  and the proof is complete.
\end{proof}

\subsubsection{Slip case}

Here we consider the boundary condition
\begin{align} \label{US0}
  \bv\cdot\bn=0 \quad \text{and} \quad \bs=\b0
  \qquad \text{on } (0,T)\times\partial\Omega
\end{align}
or the boundary condition
\begin{align} \label{US1}
  \bv\cdot\bn=0 \quad \text{and} \quad (\bs,\bv_{\btau})\in\cB
  \qquad \text{on } (0,T)\times\partial\Omega
\end{align}
where $\cB$ fulfills ($\cB$\ref{Bzero})--($\cB$\ref{Bcoercive}).
The following result holds.
\begin{Theorem} \label{th:slip}
  Let $T\in(0,\infty)$, $\Omega\subset\RR^3$ be a~$C^{0,1}$~domain,
  and $Q\coloneqq (0,T)\times\Omega$. Let $r>\frac65$,
  $\bb\in L^{r'}\bigl(0,T;(\Wn{r})^*\bigr)$, and $\bv_0\in\Lndiv2$.
  Let $\cG\subset\RR^{3\times 3}_\mathrm{sym}\times\RR^{3\times 3}_\mathrm{sym}$
  be a~maximal monotone $r$-graph of the form~\eqref{S.3*} fulfilling
  ($\cG$\ref{Gzero})--($\cG$\ref{Gcoercive}).
  \bigskip
  \begin{enumerate}[(i)]
    \item \emph{(Boundary condition~\eqref{US1})}
      Let $\cB\subset\RR^{3}\times\RR^{3}$
      be a~maximal monotone $2$-graph fulfilling
      ($\cB$\ref{Bzero})--($\cB$\ref{Bcoercive}).
      Then there exists a~triplet $(\bv,\bS,\bs)$ satisfying
      \begin{subequations}
      \label{eq:slipsol}
      \begin{align}
        \label{eq:P1as}
        \bv &\in L^\infty(0,T;\Lndiv{2})
                 \cap L^r(0,T;\Wndiv{r}),
        \\
        \label{eq:P1bs}
        \bS &\in L^{r'}(Q)^{3\times 3}_\mathrm{sym},
        \\
        \bs &\in \bL^2\bigl((0,T)\times\partial\Omega\bigr),
      \end{align}
      and
      \begin{gather}
        \label{eq:P3s}
        \lim_{t\rightarrow0+} \int_\Omega |\bv(t,\cdot)-\bv_0|^2 = 0,
        \\
        \begin{multlined}
          \label{eq:P4s}
          \int_Q \bS\mathbin{:}\bD\bw
          + \smashoperator{\int_{(0,T)\times\partial\Omega}} \bs\cdot\bw
          =
          \int_0^T \left< \bb, \bw \right>
          + \int_Q \bv\otimes\bv\mathbin{:}\bD\bw
          + \int_Q \bv\cdot\pp{\bw}{t}
          \\
          + \int_\Omega \bv_0\cdot\bw(0,\cdot)
          \begin{gathered}
            \qquad\text{for all } \bw\in \D\bigl([0,T);\Wndiv{q}\bigr),\\
            q=\max\bigl\{r,\tfrac{5r}{5r-6}\bigr\} \qquad
          \end{gathered}
        \end{multlined}
        \\
        \label{eq:P5s}
        \bS = 2\nu_* \left( |\bD\bv|-\delta_* \right)^+
              \mathcal{S}(|\bD\bv|) \tfrac{\bD\bv}{|\bD\bv|}
        \quad\text{almost everywhere in } Q,
      \end{gather}
      \end{subequations}
      \begin{gather}
        \label{eq:P7s}
        (\bs,\bw_{\btau}) \in \cB \text { almost everywhere in } (0,T)\times\partial\Omega.
      \end{gather}
    \item \emph{(Boundary condition~\eqref{US0})}
      \begingroup\setlength\emergencystretch{\hsize}\hbadness=10000
      There exists a~couple $(\bv,\bS)$ satisfying
      \eqref{eq:P1as}, \eqref{eq:P1bs}, \eqref{eq:P3s},
      \eqref{eq:P5s}, and \eqref{eq:P4s} with $\bs=\b0$.
      \par\endgroup
  \end{enumerate}
  \bigskip
  Moreover, it holds
  \begin{align*}
    \bv \in \mathcal{C}\bigl([0,T]; \bL^2_\mathrm{weak}\bigr)
    \text{ if } r>\tfrac{6}{5},
    \qquad
    \bv \in \mathcal{C}\bigl([0,T]; \bL^2\bigr)
    \text{ if } r\geq\tfrac{11}{5}.
  \end{align*}
  Also, the following energy inequality holds:
  \begin{align}
    \begin{gathered}
      \label{eq:P6s}
      \int_\Omega \frac{|\bv(t,\cdot)|^2}{2}
      + \int_0^t \int_\Omega \bS\mathbin{:}\bD\bv
      + \int_0^t \int_{\partial\Omega} \bs\cdot\bv
      \leq
      \int_\Omega \frac{|\bv_0|^2}{2}
      + \int_0^t \left<\bb, \bv\right>
      \\
      \shoveright{\text{for all } t\in(0,T];}
    \end{gathered}
  \end{align}
  if $r\geq\frac{11}{5}$, then \eqref{eq:P6s} becomes an~equality.

  In addition, if $\Omega$ is a~$C^{1,1}$ domain,
  then there is $p\in L^{q'}(Q)$ with
  $q=\max\bigl\{r,\frac{5r}{5r-6}\bigr\}$
  such that $\int_\Omega p(t,\cdot)=0$ for almost every $t\in(0,T)$ and
  \begin{align*}
    \begin{aligned}
      \int_Q \bS\mathbin{:}\bD\bom
      + \smashoperator{\int_{(0,T)\times\partial\Omega}} \bs\cdot\bw
      &= \int_0^T \bigl<\bb,\bom\bigr>
      + \int_Q \bv\otimes\bv\mathbin{:}\bD\bom
      + \int_Q \bv\cdot\pp{\bom}{t}
      + \int_\Omega \bv_0\cdot\bom(0,\cdot)
      \\
      &+ \int_Q p\diver\bom
      \shoveright{\quad\text{for all } \bom\in \D\bigl([0,T);\Wn{q}\bigr).}
    \end{aligned}
  \end{align*}
\end{Theorem}
\begin{Rem}
      In the case of the boundary condition~\eqref{US0}
      we could define the weak solution to the problem considered differently.
      We could say that $\bv$ is a~weak solution to the problem if
      $\bv$ fulfills~\eqref{eq:P1as}, \eqref{eq:P3s}, and
      \begin{multline*}
        \begin{aligned}
          \int_Q 2\nu_* \left( |\bD\bv|-\delta_* \right)^+
                 \mathcal{S}(|\bD\bv|) \tfrac{\bD\bv}{|\bD\bv|}\mathbin{:}\bD\bw
          & =
          \int_0^T \left<\bb, \bw\right>
          + \int_Q \bv\otimes\bv\mathbin{:}\bD\bw
          \\ &
          + \int_Q \bv\cdot\pp{\bw}{t}
          + \int_\Omega \bv_0\cdot\bw(0,\cdot)
        \end{aligned}
        \\
        \text{for all } \bw\in \D\bigl([0,T);\Wndiv{q}\bigr),\,
        q=\max\bigl\{r,\tfrac{5r}{5r-6}\bigr\}.
      \end{multline*}
\end{Rem}

\begin{proof}[Proof of Theorem~\ref{th:slip}]
  We focus only on the details in which the proof differs from the proof of
  Theorem~\ref{th:noslip}. Note however that a~remarkable difference
  concerns the pressure: for the \noslip{} boundary condition the pressure is
  not integrable up to the boundary; here, for $C^{1,1}$ domains, we
  establish the existence of the pressure belonging to  $L^s(Q)$ for some
  $s>1$. This concerns in particular the \noslip/\navierslip{} boundary
  condition which ``approximates'' well the \noslip{} boundary condition
  and in addition its mathematical theory admits integrable pressure.

  Regarding the case $r\geq\tfrac{11}{5}$, the main departures from the problem with
  the \noslip{} boundary condition is due to the choice of function spaces and due to the
  formulation of the eigenvalue problem that generates the basis for Galerkin
  approximations. Here, we look for $\lambda\in\RR$ and $\bom\in\Vndiv
  \hookrightarrow W^{1,\infty}(\Omega)^3$ satisfying
  \begin{align*}
    \llparen\bom,\bfi\rrparen = \lambda (\bom,\bfi)
    \text{ for all }\bfi\in\Vndiv,
  \end{align*}
  where $(\cdot,\cdot)$ is again the scalar product in $\bL^2$ and
  $\llparen\cdot,\cdot\rrparen$ is a~scalar product in
  $\Vndiv$ defined through
  $\llparen\bom,\bfi\rrparen\coloneqq
    (\nabla^3\bom,\nabla^3\bfi) + (\bom,\bfi)
    + (\bom_{\btau},\bfi_{\btau})_{\partial\Omega}
  $.
  The properties of the eigenfunctions are the same as in~\eqref{eq:p1} and
  consequently, for the \freeslip{} boundary condition~\eqref{US0} there is no
  other change in the proof.

  If the other slip-type conditions are considered, then we regularize the
  boundary conditions as in the time independent case. Independent of the
  approximation parameter, we, in addition to standard uniform estimates, know
  that $\{\bs^n\}_{n=1}^\infty$ and $\{\bv^n_{\btau}\}_{n=1}^\infty$ are
  bounded in $L^2(0,T; \bL^2(\partial\Omega))$. Furthermore, as $W^{1,r}(\Omega)$
  compactly embeds into $W^{\frac{1}{r},q}(\partial \Omega)$ for all $q<r$, we
  conclude that
  \begin{align*}
    \bv^N_{\btau} \rightarrow \bv_{\btau}
    \quad \text{strongly in } L^r\bigl( 0,T; \bL^1(\partial\Omega) \bigr).
  \end{align*}
  Then (up to a~subsequence which we do not relabel)
  \begin{align*}
    \bv^N_{\btau} \rightarrow \bv_{\btau}
    \quad \text{a.e. on } (0,T)\times\partial\Omega
  \end{align*}
  and by Egorov's theorem, for any $\delta>0$,
  \begin{align*}
    \bv^N_{\btau} \rightarrow \bv_{\btau}
    \quad \text{strongly in } L^\infty(\mathcal{U}_\delta)
  \end{align*}
  where $\mathcal{U}_\delta\subset(0,T)\times\partial\Omega$ is such that
  $|(0,T)\times\partial\Omega\setminus\mathcal{U}_\delta|<\delta$.
  The last convergence implies that
  \begin{align*}
    \limsup_{N\rightarrow\infty}
    \int_{\mathcal{U}_\delta} \bs^N\mathbin{\cdot}\bv^N_{\btau}
    =
    \int_{\mathcal{U}_\delta} \bs\mathbin{\cdot}\bv_{\btau}.
  \end{align*}
  Consequently, by Lemma~\ref{lemma:conv}, $(\bs,\bv_{\btau})\in\cB$
  a.e.~on $\mathcal{U}_\delta$. This is true for all $r>1$ and
  gives~\eqref{eq:P7s}.

  If $r\in(\tfrac65,\tfrac95)$, the proof of $(\bS,\bD\bv)\in\cG$ is
  carried out as in the \noslip{} case, as the proof is based on local analysis
  in the interior of~$\Omega$.

  Finally, we reconstruct the pressure. We set
  $p = p_1 + p_2$ where $p_1\in L^\frac{5r}{6}(Q)$ solves
  \begin{align*}
    (p_1,-\Delta z) &= (\bv\otimes\bv,\nabla^2 z)
    &&\text{for all } z\in W^{2,\frac{5r}{5r-6}}
    \text{ with } \nabla z\in\Wn{\frac{5r}{5r-6}},
    \\
    \int_\Omega p_1(t) &= 0
    &&\text{for a.e. } t\in(0,T)
  \end{align*}
  and $p_2\in L^{r'}(Q)$ solves
  \begin{align*}
    (p_2,-\Delta z) &= \langle\bb,\nabla z\rangle
                     - (\bS,\bD\nabla z)
                     - (\bs,\nabla z)_{\partial\Omega}
    &&\text{for all } z\in W^{2,\frac{5r}{5r-6}}
    \text{, } \nabla z\in\Wn{\frac{5r}{5r-6}},
    \\
    \int_\Omega p_2(t) &= 0
    &&\text{for a.e. } t\in(0,T).
  \end{align*}
  Note that this is a~well-posed definition because of the $C^{1,1}$
  regularity of the domain~$\Omega$ and Lemma~\ref{lemma:neumannregularity}.
  Now consider a~test function $\bfi\in L^q(0,T;\Wn{q})$ and its
  Helmholtz decomposition using Corollary~\ref{cor:helmholtz}:
  \begin{align*}
    \bfi = \nabla\phi+\bfi_0
    \quad \text{with }
    \nabla\phi\in L^{q}(0,T;\Wn{q}),\,\bfi_0\in L^{q}(0,T;\Wndiv{q}).
  \end{align*}
  Then we have
  \begin{align*}
    &
    \Bigl\langle\pp{\bv}{t},\bfi\Bigr\rangle
    -(p,\diver\bfi)
    \\
    &\quad
    =
    \Bigl\langle\pp{\bv}{t},\nabla\phi+\bfi_0\Bigr\rangle
    -(p,\diver(\nabla\phi+\bfi_0))
    =
    \Bigl\langle\pp{\bv}{t},\bfi_0\Bigr\rangle
    +(p_1+p_2,-\Delta\phi)
    \\
    &\quad
    =
    (\bv\otimes\bv,\bD\bfi_0) - (\bS,\bD\bfi_0)
    - \langle\bs,\bfi_0\rangle_{\partial\Omega}
    + \langle\bb,\bfi_0\rangle + (\bv\otimes\bv,\nabla^2\phi)
    \\
    &\quad
    \phantom{=}
    + \langle\bb,\nabla\phi\rangle - (\bS,\bD\nabla\phi)
    - \langle\bs,\nabla\phi\rangle_{\partial\Omega}
    \\
    &\quad
    =
    (\bv\otimes\bv,\bD\bfi) - (\bS,\bD\bfi)
    - \langle\bs,\bfi_{\btau}\rangle_{\partial\Omega}
    + \langle\bb,\bfi\rangle
  \end{align*}
  for a.a. $t\in(0,T)$. Thus the theorem is proven.
\end{proof}

\section{Concluding remarks} \label{Ch3}
We have classified incompressible fluids that span the gamut from Euler fluids
-- Navier-Stokes fluid -- power-law fluids -- generalized power-law fluids --
stress power-law fluids -- to fluids that only undergo rigid motions, that can change
their constitutive character due to an activation criterion based on the value
of the norm of the symmetric part of the velocity gradient or the shear stress.
In the process we came across constitutive relations that have hitherto been
unrecognized but could possibly be useful. In the course of our investigation we have delineated
how an Euler fluid is different from a fluid that behaves like an Euler fluid
prior activation and behaves like a viscous fluid when the activation criterion
takes place.  The latter fluid would lead to governing equations that imbed the
boundary layer equations as a~special case, the philosophy behind the
development of the boundary layer equations and the equations governing the
activated fluid being totally different. We have touched upon one important aspect in this
study, namely the tremendously different properties that are exhibited by the
Euler fluids and the activated Euler fluids. It is known that while the Euler
fluid exhibits pathological features (such as existence of a~nontrivial solution
to internal flows with zero initial data and vanishing external body forces),
we have shown that the new class of activated Euler fluids admits a~weak solution that might be even
unique in its dependence of what kind of response occurs after activation.

A~classification similar to that presented here for incompressible fluids can be carried
out within the context of compressible fluids, where however the framework is
more complicated as there are two type of viscosities (bulk and shear) and
corresponding fluidities. This issues will be addressed in a subsequent study.

\section*{Acknowledgement}
The authors would like to thank Endre S\"uli and the anonymous referees
for stimulating comments and suggestions.

\appendix
\section{Auxiliary convergence tools}
In this section, we state, without proofs, several characterizations of
weak compactness in~$L^1$. Then, following~\cite{BDS13},
we summarize several properties of refined (divergence-free) Lipschitz approximations of
(divergence-free) Sobolev and Bochner-Sobolev functions. Next, we present a~convergence
lemma (proved recently in~\cite{bulcek-malek-b:2016}) regarding stability of maximal monotone constitutive equations
(maximal monotone $r$-graphs) with respect to weakly converging sequences.
Finally, we close this section by the Ne\v{c}as theorem and
Sobolev regularity results for the Neumann-Poisson problem and the
Stokes system.

In the following lemma, several assertions characterizing
weak compactness in $L^1$, namely the Dunford-Pettis
criterion (ii), uniform integrability (iii), and the de~la Vall\'e-Poussin
criterion (iv), are provided. The exact statement is taken from
\cite[p.~21, Theorem~10]{feireisl-karper-pokorny}.
\begin{Lem}[Characterization of weak compactness in $L^1$] \label{lemma:char}
  Let $Q\subset\RR^M$ be a~bounded measurable set and
  $\mathcal{V}\subset L^1(Q)$. Then the following conditions
  are equivalent:
  \begin{enumerate}[(i)]
    \item
      any sequence $\{v_n\}_{n=1}^\infty \subset \mathcal{V}$
      contains a~subsequence weakly converging in $L^1(Q)$;
    \item
      for any $\epsilon>0$ there exists $K>0$ such that for all $v\in\mathcal{V}$
      \begin{equation*}
        \int_{\{|v|\geq K\}} |v(y)| \d y \leq \epsilon;
      \end{equation*}
    \item \label{lemma:char:equi}
      for any $\epsilon>0$ there exists $\delta>0$ such that for
      all $v\in\mathcal{V}$
      and for any measurable set $M\subset Q$ such that $|M|<\delta$
      \begin{equation*}
        \int_M |v(y)| \d y < \epsilon;
      \end{equation*}
    \item
      there exists a~nonnegative function $\Phi\in \mathcal{C}([0,\infty))$
      fulfilling
      \begin{equation*}
        \lim_{z\rightarrow\infty} \frac{\Phi(z)}{z} = \infty,
      \end{equation*}
      such that
      \begin{equation*}
        \sup_{v\in\mathcal{V}} \int_Q \Phi(|v(y)|) \d y < \infty.
      \end{equation*}
  \end{enumerate}
\end{Lem}

Since $L^1$ is not reflexive, weak precompactness does not follow from
boundedness. Instead bounded sequences in $L^1$ can exhibit
local concentrations weakly converging only in the space of measures.
The next lemma ensures that these concentrations are located on arbitrarily
small sets and when removed (by ``biting''), bounded sets
are $L^1$-weak precompact on the complement (''unbitten''
part). See original reference~\cite{brooks-chacon-1980} and
also~\cite{BM89} for a~simple proof and other references.
\begin{Lem}[Biting lemma] \label{lemma:biting}
  Let $Q\subset\RR^M$ be bounded and measurable. Let
  $\{v_n\}_{n=1}^\infty$ be a~sequence bounded in $L^1(Q)$.
  Then there exist a~subsequence
  $\{v_{n_j}\}_{j=1}^\infty \subset \{v_n\}_{n=1}^\infty$,
  a~function $v\in L^1(Q)$, and a~sequence of measurable sets
  $\{E_k\}_{k=1}^\infty$, $Q\supset E_1 \supset E_2 \supset \ldots$,
  $|E_k|\rightarrow0$ such that for all $k\in\mathbb{N}$
  \begin{equation*}
    v_{n_j} \rightharpoonup v
    \quad\text{ weakly in } L^1(Q\setminus E_k) \text{ as } j\rightarrow\infty.
  \end{equation*}
\end{Lem}

In the following corollary of the preceding lemmas we establish strong
convergence in $L^1$ up to arbitrarily small sets for a~pointwise null
sequence bounded in $L^1$.
\begin{Col} \label{col:strongbiting}
  Let the assumptions of Lemma~\ref{lemma:biting} be fulfilled. Furthermore,
  assume that
  \begin{equation*}
	v_n\rightarrow 0
    \quad \text{ a.e. in } Q \text{ as } n\rightarrow\infty.
  \end{equation*}
  Then for the sequences $\{v_{n_j}\}_{j=1}^\infty$
  and $\{E_k\}_{k=1}^\infty$ from Lemma~\ref{lemma:biting}
  and for every $k\in\mathbb{N}$
  \begin{equation*}
    v_{n_j} \rightarrow 0
    \quad\text{ strongly in } L^1(Q\setminus E_k) \text{ as } j\rightarrow\infty.
  \end{equation*}
\end{Col}
\begin{proof}
  Let $k\in\mathbb{N}$ be fixed. The sequence
  $\{v_{n_j}\}_{j=1}^\infty$ provided by the biting
  lemma~\ref{lemma:biting} is weakly compact in $L^1(Q\setminus E_k)$ and by
  the lemma~\ref{lemma:char}, (ii), $\{v_{n_j}\}_{j=1}^\infty$ is
  uniformly continuous with respect to the Lebesgue measure on
  $Q\setminus E_k$. By the Vitali convergence theorem, the assertions follows.
\end{proof}

Lipschitz approximations of solenoidal Bochner-Sobolev functions is another
useful tool needed in the analysis of isochoric flows. There are several
variants: Acerbi and Fusco survey the basic properties of Lipschitz
approximations of Sobolev functions in~\cite{acerbi.e.fusco.n:approximation};
further refinements have been put into place, see~\cite{frehse.j.malek.j.ea:on,
diening.l.malek.j.ea:on}. The extension to evolutionary problems goes back
to~\cite{kinnunen-lewis-2000,kinnunen-lewis-2002,DiRuWo06}. Further extensions
have been established in~\cite{BGMS2012,BDS13}.

We first state the version \cite[Theorem~4.2]{BDS13},
which is suitable for analysis of steady problems.
\begin{Lem}[Divergence-free Lipschitz truncation of Sobolev functions]
  \label{lemma:truncsteady}
  Let $B\subset\RR^3$ be an~arbitrary ball.
  Let $r\in(1,\infty)$.
  Let $\{\bu^m\}_{m=1}^\infty \subset \Wzdiv{r}(B)$ be weakly converging to zero
  in $\Wzdiv{r}(B)$.

  Then there is a~double sequence
  $\{\lambda_{m,k}\}_{m,k=1}^\infty \subset (0,\infty)$ with
  \begin{enumerate}[(a)]
    \item \label{lemma:truncsteady:lambda}
      $ 2^{2^k} \leq \lambda_{m,k} \leq 2^{2^{k+1}-1} $,
  \end{enumerate}
  a~double sequence of functions
  $\{\bu^{m,k}\}_{m,k=1}^\infty$,
  a~double sequence
  $\{\mathcal{O}^{m,k}\}_{m,k=1}^\infty$
  of measurable subsets of~$2B$,
  a~constant $C>0$,
  and $k_0\in\mathbb{N}$
  such that for all $k \geq k_0$ it holds:
  \begin{enumerate}[(a)] \setcounter{enumi}{1}
    \item \label{lemma:truncsteady:lipapp}
      $ \bu^{m,k}\in \Wzdiv{\infty}(2B) $
      and
      $\bu^{m,k}=\bu^m$ in $2B\setminus\mathcal{O}^{m,k}$
      for all $m\in\mathbb{N}$,
    \item \label{lemma:truncsteady:bound}
      $ \|\nabla\bu^{m,k}\|_{L^\infty(2B)} \leq C\, \lambda_{m,k} $
      for all $m\in\mathbb{N}$,
    \item \label{lemma:truncsteady:strongconv}
      $ \bu^{m,k} \rightarrow 0 $
      strongly in $\bL^\infty(2B)$
      as $m\rightarrow\infty$,
    \item \label{lemma:truncsteady:weakstarconv}
      $ \nabla\bu^{m,k} \stackrel{*}{\rightharpoonup} 0 $
      weakly-$*$ in $L^\infty(2B)^{3\times 3}$
      as $m\rightarrow\infty$,
    \item \label{lemma:truncsteady:smallness}
      $ (\lambda_{m,k})^r |\mathcal{O}^{m,k}| \leq C\, 2^{-k} \|\nabla\bu^m\|_r^r $
      for all $m\in\mathbb{N}$.
  \end{enumerate}
\end{Lem}

Next we will formulate the assertion suitable for analysis of time-dependent
problems; the presented version is taken from~\cite{BDS13}.
\begin{Lem}[Divergence-free Lipschitz truncation of Bochner-Sobolev functions]
  \label{lemma:trunc}
  Let $Q_0=I_0\times B_0 \subset \RR\times\RR^3$ be
  a~space-time cylinder.
  Let $1<r<\infty$ with $r,r'>\sigma>1$, $\frac1r+\frac{1}{r'}=1$.
  Assume that there are sequences of functions
  $\{\bu^m\}_{m=1}^\infty$ and
  $\{\bH^m\}_{m=1}^\infty$ such that
  \begin{gather*}
    \diver\bu^m = 0 \quad\text{ a.e. in } Q_0, \\
    -\int_{Q_0} \bu^m \cdot \partial_t \bom = \int_{Q_0} \bH^m \mathbin{:} \nabla\bom
    \quad\text{ for all } \bom\in\D(I_0;\Cinfzdiv(B_0)), \\
    \begin{aligned}
      \bu^m &\text{ is uniformly bounded in } L^\infty(I_0; \bL^{\sigma}(B_0)), \\
      \bu^m &\rightharpoonup \b0 \text{ weakly   in } L^r(I_0; \bW^{1,r}(B_0)), \\
      \bu^m &\rightarrow     \b0 \text{ strongly in } \bL^\sigma(Q_0),
    \end{aligned}
  \end{gather*}
  and $\bH^m = \bH^m_1 + \bH^m_2$ satisfies
  \begin{align*}
    \bH^m_1 &\rightharpoonup \bO \text{ weakly   in } L^{r'}(Q_0)^{3\times 3}, \\
    \bH^m_2 &\rightarrow     \bO \text{ strongly in } L^\sigma(Q_0)^{3\times 3}.
  \end{align*}
  Then there is a~double sequence
  $\{\lambda^{m,k}\}_{m,k=1}^\infty \subset (0,\infty)$
  with
  \begin{enumerate}[(a)]
    \item \label{lemma:trunc:lambda}
      $ 2^{2^k} \leq \lambda_{m,k} \leq 2^{2^{k+1}} $,
  \end{enumerate}
  \begingroup\setlength\emergencystretch{\hsize}\hbadness=10000
  a~double sequence of functions
  $\{\bu^{m,k}\}_{m,k=1}^\infty \subset \bL^1(Q_0)$,
  a~double sequence
  $\{\mathcal{O}^{m,k}\}_{m,k=1}^\infty$
  of measurable subsets of~$Q_0$,
  $C>0$,
  and $k_0\in\mathbb{N}$
  such that for all $k \geq k_0$:
  \begin{enumerate}[(a)] \setcounter{enumi}{1}
    \item
      $ \bu^{m,k} \in L^s(\frac14I_0; \Wzdiv{s}(\frac16B_0)) $
      for all $s\in(1,\infty)$
      and all $m\in\mathbb{N}$,
    \item
      $\supp\bu^{m,k}\subset\frac16Q_0$
      for all $m\in\mathbb{N}$,
    \item
      $ \bu^{m,k} = \bu^m $ a.e. on $\frac18Q_0\setminus\mathcal{O}^{m,k}$
      for all $m\in\mathbb{N}$,
    \item
      $ \|\nabla\bu^{m,k}\|_{L^\infty(\frac14Q_0)} \leq C\,\lambda^{m,k}$
      for all $m\in\mathbb{N}$,
    \item
      $ \bu^{m,k} \rightarrow 0 $ strongly in $\bL^\infty(\frac14Q_0)$
      as $m\rightarrow\infty$,
    \item
      $ \nabla\bu^{m,k} \rightharpoonup 0 $ weakly in $L^s(\frac14Q_0)^{3\times 3}$
      for all $s\in(1,\infty)$
      as $m\rightarrow\infty$,
    \item \label{lemma:trunc:smallness}
      $ \limsup_{m\rightarrow\infty} (\lambda^{m,k})^r |\mathcal{O}^{m,k}|
        \leq C\, 2^{-k}$,
    \item
      $ \limsup_{m\rightarrow\infty} \bigl| \int_{\frac{1}{6}Q_0} \bH^m\mathbin{:}\nabla\bu^{m,k} \bigr|
        \leq C\,
        \limsup_{m\rightarrow\infty} (\lambda^{m,k})^r |\mathcal{O}^{m,k}| $,
    \item \label{lemma:trunc:conseq}
      $
        \limsup_{m\rightarrow\infty} \bigl|
          \int_{\frac{1}{6}Q_0\setminus\mathcal{O}^{m,k}} \bH^m_1 \mathbin{:} \nabla\bu^m \,\xi
        \bigr|
        \leq C\,2^{-\frac{k}{r}}
      $
      with some $\xi\in\D\bigl(\frac{Q_0}{6}\bigr)$
      that fulfills $\chi_{\frac18Q_0} \leq \xi \leq \chi_{\frac16Q_0}$.
  \end{enumerate}
  \endgroup
\end{Lem}

Another tool that we use is a~simple lemma concerning the stability of the
constitutive equations (represented as maximal monotone $r$-graphs) with
respect to weakly converging sequences; see~\cite{bulcek-malek-a:2016} for
a~short proof.
\begin{Lem}[Graph convergence lemma] \label{lemma:conv}
  Let $D\subset\RR^M$ be an arbitrary measurable set and let a~graph $\cG$
  fulfill the assumptions~($\cG$\ref{Gmonotone}) and~($\cG$\ref{Gmaximal})
  on page~\pageref{Gmonotone}.  Assume that, for some $r\in (1, \infty)$,
  \begin{gather*}
    (\bS^n,\bD^n)\in \cG \quad\text{almost everywhere in } D, \\
    \begin{aligned}
      \bD^n &\rightharpoonup \bD &&\text{weakly in } L^{r}(D)^{d\times d}, \\
      \bS^n &\rightharpoonup \bS &&\text{weakly in } L^{\frac{r}{r-1}}(D)^{d\times d},
    \end{aligned}
    \\
    \limsup_{n\to \infty} \int_{D} \bS^n \mathbin{:} \bD^n \le \int_{D} \bS \mathbin{:} \bD.
  \end{gather*}
  Then
  \begin{align*}
    (\bS,\bD)\in \cG \text{ almost everywhere in }D.
  \end{align*}
\end{Lem}

Next, we state a~theorem due to Ne\v{c}as~\cite{necas1967};
the following version is from~\cite[Corollary~2.5~ii)]{ag94}.
\begin{Lem}[Ne\v{c}as theorem] \label{lemma:necas}
  Let $\Omega\subset\RR^M$ be a~domain of class $C^{0,1}$.
  Let $r\in(1,\infty)$.
  Then there exists $\beta>0$ such that
  \begin{align}
    \label{eq:necas}
    \|\nabla q\|_{\bigl(\Wz{r}\bigr)^*} \coloneqq
    \sup_{\bfi\in\Wz{r}} \frac{(q,\diver\bfi)}{\|\nabla\bfi\|_r}
    \geq \beta \|q\|_{r'}
    \quad\text{for all } q\in L^{r'}(\Omega)
    \text{ with } \int_\Omega q=0.
  \end{align}
\end{Lem}
\begin{Rem} \label{rem:bog}
  Lemma~\ref{lemma:necas} is closely related to the results known as
  the Lions lemma (coined in~\cite{MR0123818}),
  the Babu\v{s}ka-Aziz inequality,
  or the Ladyzhenskaya-Babu\v{s}ka-Brezzi condition
  (see~\cite{Babuska1973,Brezzi1974}).
  If we set
  $L^p_0\coloneqq\bigl\{q\in L^p(\Omega),\,\int_\Omega q=0\bigr\}$,
  the statement of Lemma~\ref{lemma:necas} can also be rephrased as
  the following:
  \begin{enumerate}[(i)]
    \item
      the gradient operator $\nabla:L^{r'}_0\rightarrow (\Wz{r})^*$
      is injective with closed range,
    \item \label{prop:bog}
      the divergence operator $\diver:\Wz{r}\rightarrow L^r_0$ is surjective
      and has a~continuous right inverse, i.e., there is a~bounded linear
      operator $\mathcal{B}:L^r_0\rightarrow\Wz{r}$ such that
      \begin{align*}
        \diver\mathcal{B} \text{ is identity on } L^r_0.
      \end{align*}
  \end{enumerate}
  The operator $\mathcal{B}$ is usually called the Bogovski\u{i}
  operator due to the explicit construction by
  Bogovski\u{i}~\cite{Bogovskii1979,Bogovskii1980}.
  We refer the reader to~\cite{ag94}, where these relations are
  discussed in detail.
\end{Rem}
\begin{Rem} \label{rem:necas}
  It is shown in~\cite[Proposition~2.10~ii)]{ag94} that
  for the validity of the estimate of Lemma~\ref{lemma:necas} it is
  sufficient to assume a~priori that $q\in\bigl(\D(\Omega)\bigr)^*$,
  $\int_\Omega q=0$, and $\nabla q\in \bigl(\Wz{r}\bigr)^*$.
  Then, provided that $\Omega$ is Lipschitz,
  $q\in L^{r'}(\Omega)$ and the estimate~\eqref{eq:necas} holds.
\end{Rem}
\begin{Rem}
  The Lipschitz condition on $\Omega$ in
  Lemma~\ref{lemma:necas} can be weakened,
  see for example \cite{diening-ruzicka-schumacher-2010}.
\end{Rem}

Now we mention few regularity results for the Neumann problem
and the Stokes system.
\begin{Lem}[$W^{2,q}$-regularity of Neumann-Poisson problem]
  \label{lemma:neumannregularity}
  Let $\Omega\subset\RR^M$ be a~domain of class $C^{1,1}$.
  Let $1<q<\infty$ be given. Then there exists $C>0$
  such that for every $f\in L^q(\Omega)$ with $\int_\Omega f = 0$
  there is a~weak solution $u\in W^{1,q}(\Omega)$ of the problem
  \begin{subequations} \label{SS_laplace}
  \begin{align}
    -\Delta u                       &= f && \text{ in }\Omega, \\
    \frac{\partial u}{\partial\bn } &= 0 && \text{ on }\partial\Omega, \\
    \int_\Omega u                   &= 0
  \end{align}
  \end{subequations}
  fulfilling $u\in W^{2,q}(\Omega)$, $\nabla u\in\Wn{q}$,
  and
  \begin{align*}
    \|\nabla^2 u\|_q \leq C \|f\|_q.
  \end{align*}
\end{Lem}
\noindent
Proof of Lemma~\ref{lemma:neumannregularity} is outlined in
\cite[Remark~3.2]{AGG1997} and invokes \cite{ADN1959,Gr85,LionsMagenes1962}.
As a~consequence of the lemma, we get the following result concerning the
Helmholtz decomposition for functions from $\Wn{q}$.
\begin{Col}[Helmholtz decomposition] \label{cor:helmholtz}
  Let $\Omega\subset\RR^M$ be a~domain of class $C^{1,1}$.
  Let $1<q<\infty$ be given. Then there exists $C>0$ such
  that the following holds.
  For every $\bfi\in\Wn{q}$ there exists a~couple $(\phi,\bfi_0)$
  fulfilling
  \begin{gather*}
    \phi \in W^{2,q}(\Omega),
    \qquad
    \nabla\phi \in \Wn{q},
    \qquad
    \bfi_0 \in \Wndiv{q}, \\
    \bfi = \nabla\phi + \bfi_0,
    \qquad
    \|\nabla^2\phi\|_q + \|\nabla\bfi_0\|_q \leq C \|\nabla\bfi\|_q.
  \end{gather*}
\end{Col}
The following lemma contains certain regularity results for the Stokes system,
see~\cite[Theorem 2.1]{GSS1994}, \cite[Proposition 4.3]{ag94}.
\begin{Lem}[Regularity of the Stokes system]
  \label{lemma:stokesregularity}
  Let $\Omega\subset\RR^M$ be a~domain and $1<q<\infty$ be given.

  If $\Omega$ is of class $C^{0,1}$ with sufficiently small Lipschitz
  constant $L>0$ (i.e., $L\leq L_0$ with $L_0>0$ depending only
  on $M$ and $q$) or $\Omega$ is of class $C^1$
  then there exists $C_0>0$ (depending on $\Omega$, $M$, $q$)
  such that for every $\bb\in\bigl(\Wz{q'}\bigr)^*$
  there is a~unique weak solution $(\bv,p)\in\Wz{q}\times L^{q}(\Omega)$
  of the problem
  \begin{subequations} \label{eq:stokes}
  \begin{align}
    -\Delta \bv + \nabla p &= \bb && \text{ in }\Omega, \\
    \diver \bv             &= 0   && \text{ in }\Omega, \\
    \bv                    &= \b0 && \text{ on }\partial\Omega
  \end{align}
  \end{subequations}
  and the following estimate holds true
  \begin{align}
    \label{eq:stokes_w1q}
    \|\nabla\bv\|_q + \|p\|_q \leq C_0 \|\bb\|_{\bigl(\Wz{q'}\bigr)^*}.
  \end{align}
  Furthermore, if $\Omega$ is of class $C^{1,1}$
  then there exists $C_1>0$ (depending on $\Omega$, $M$, $q$)
  such that for every $\bb\in L^q(\Omega)^M$
  the unique weak solution $(\bv,p)\in\Wz{q}\times L^{q}(\Omega)$
  of the problem~\eqref{eq:stokes} fulfills additionally $\bv\in\bW^{2,q}$,
  $p\in W^{1,q}(\Omega)$, and admits the estimate
  \begin{align}
    \label{eq:stokes_w2q}
    \|\nabla^2\bv\|_q + \|\nabla p\|_q \leq C_1 \|\bb\|_q.
  \end{align}
\end{Lem}
\begin{proof}
  The first part of the lemma is exactly the statement \cite[Theorem
  2.1]{GSS1994}. This statement guarantees existence of unique
  $(\bv,p)\in\Wz{q}\times L^{q}(\Omega)$ and gives~\eqref{eq:stokes_w1q} under
  the aforementioned conditions.

  The second part, i.e., the inclusions $\bv\in\bW^{2,q}$, $p\in W^{1,q}(\Omega)$
  and the estimate~\eqref{eq:stokes_w2q} follow from \cite[Proposition
  4.3]{ag94}.  Remark~4.4 therein warns that Proposition~4.3 ibid. is not
  an~existence result and that~\eqref{eq:stokes_w2q} holds only if a~solution
  in the appropriate spaces exists. But we know this is the case due to the
  first part, by virtue of~\cite{GSS1994}.
\end{proof}
 
\section{Examples of maximal monotone graphs} \label{app:examples}
Let us consider a~graph
$\cG\subset\RR^{3\times 3}_\mathrm{sym}\times\RR^{3\times 3}_\mathrm{sym}$
characterized by the relationship
\begin{align}
  \label{A.1}
  (\bS,\bD)\in\cG
  \Leftrightarrow
  \bS = \frac{\left(|\bD|-\delta_*\right)^+}{|\bD|}\mathcal{S}(|\bD|)\bD,i
  \qquad \mathcal(t)^+ = \left\{\begin{array}{ll}
    t & t>0 \\
    0 & \text{otherwise}
  \end{array}\right.
\end{align}
with
\begin{subequations}
  \label{A.2}
  \begin{alignat}{2}
    \label{A.2a}
    &\text{either}\qquad&
    \mathcal{S}(d) &= (1 + d^2)^\frac{r-2}{2},
    \\
    \label{A.2b}
    &\text{or}\qquad&
    \mathcal{S}(d) &= 1 + d^{r-2}.
  \end{alignat}
\end{subequations}
We will prove the following statement.
\begin{Lem}
  The graph $\cG$ characterized by~\eqref{A.1} and~\eqref{A.2a}
  with some $\delta_*\geq0$ and $r\in(1,\infty)$ is
  a~maximal monotone $r$-graph fulfilling
  ($\cG$\ref{Gzero})--($\cG$\ref{Gcoercive}).

  The graph $\cG$ characterized by~\eqref{A.1} and~\eqref{A.2b}
  with some $\delta_*\geq0$ and $r\in(1,\infty)$ is
  a~maximal monotone $q$-graph fulfilling
  ($\cG$\ref{Gzero})--($\cG$\ref{Gcoercive})
  with $q=\max\{r,2\}$.
\end{Lem}
\begin{proof}
  \noindent(i). Clearly $(\bO,\bO)\in\cG$.

  \noindent(ii). Let $\bS=\frac{\left(|\bD|-\delta_*\right)^+}{|\bD|}\left(1+|\bD|^2\right)^\frac{r-2}{2}\bD$
  and $\bD_s\coloneqq\bD_2+s(\bD_1-\bD_2)$ for any
  $\bD_1,\bD_2\in\RR^{3\times 3}_\mathrm{sym}$. Then
  \begin{flalign*}
    \left(\bS(\bD_1)-\bS(\bD_2)\right)\mathbin{:}\left(\bD_1-\bD_2\right) &&
  \end{flalign*}
  \begin{align*}
    &= \left(\bD_1-\bD_2\right)
       \mathbin{:}\int_0^1\frac{\d}{\d s}
         \left[\frac{\left(|\bD_s|-\delta_*\right)^+}{|\bD_s|}
               \left(1+|\bD_s|^2\right)^\frac{r-2}{2}\bD_s
         \right]\,\d s \\
    &= \left|\bD_1-\bD_2\right|^2
       \,\int_0^1 \frac{\left(|\bD_s|-\delta_*\right)^+}{|\bD_s|}
                  \left(1+|\bD_s|^2\right)^\frac{r-2}{2}
                  \,\d s \\
    &+ \int_0^1\left(\bD_s\mathbin{:}(\bD_1-\bD_2)\right)^2
       \Biggl\{
         \mathcal{H}\left(|\bD_s|-\delta_*\right)\frac{\left(1+|\bD_s|^2\right)^\frac{r-2}{2}}{|\bD_s|^2}
         \\&\qquad
         +(r-2)\frac{\left(|\bD_s|-\delta_*\right)^+}{|\bD_s|}\left(1+|\bD_s|^2\right)^\frac{r-4}{2}
         -\frac{\left(|\bD_s|-\delta_*\right)^+}{|\bD_s|^3}\left(1+|\bD_s|^2\right)^\frac{r-2}{2}
       \Biggr\},
  \end{align*}
  where $\mathcal{H}(t)=1$ for $t>0$, $\mathcal{H}(t)=0$ otherwise.
  Since
  $\mathcal{H}\left(|\bD_s|-\delta_*\right)
    -\frac{\left(|\bD_s|-\delta_*\right)^+}{|\bD_s|}
    =\frac{\delta_*}{|\bD_s|}$
  if $|\bD_s|>\delta_*$ (otherwise it is zero),
  we observe that for $r\geq 2$
  \begin{align*}
    \left(\bS(\bD_1)-\bS(\bD_2)\right)\mathbin{:}\left(\bD_1-\bD_2\right)\geq0.
  \end{align*}
  If $r\in(1,2)$, then the property ($\cG$\ref{Gmonotone}) follows as well
  from the fact
  \begin{flalign*}
    &
    |\bD_1-\bD_2|^2 \,\int_0^1 \frac{\left(|\bD_s|-\delta_*\right)^+}{|\bD_s|}
                               \left(1+|\bD_s|^2\right)^\frac{r-2}{2} \,\d s
    &\\&
    -(2-r) \,\int_0^1 \left(\bD_s\mathbin{:}(\bD_1-\bD_2)\right)^2
                      \frac{\left(|\bD_s|-\delta_*\right)^+}{|\bD_s|}
                      \left(1+|\bD_s|^2\right)^\frac{r-4}{2} \,\d s
    &
  \end{flalign*}
  \begin{flalign*}
    &&
    \geq |\bD_1-\bD_2|^2 \,\int_0^1 \frac{\left(|\bD_s|-\delta_*\right)^+}{|\bD_s|}
                                    \left(1+|\bD_s|^2\right)^\frac{r-4}{2}
                                    \left(1+|\bD_s|^2-(2-r)|\bD_s|^2\right) \,\d s
    &\\&&
    \geq (r-1) |\bD_1-\bD_2|^2 \,\int_0^1 \frac{\left(|\bD_s|-\delta_*\right)^+}{|\bD_s|}
                                          \left(1+|\bD_s|^2\right)^\frac{r-2}{2} \,\d s
    \geq0.
    &
  \end{flalign*}
  This also implies the monotone property for the graph with
  $\mathcal{S}(d)=d^{r-2}$ and $r>1$. Consequently, the same is true for
  $\mathcal{S}(d)=1+d^{r-2}$.

  \noindent(iii). In order to show that the graph is a~maximal monotone graph,
  we note that the assumption:
  $(\bS,\bD)\in\RR^{3\times 3}_\mathrm{sym}\times\RR^{3\times 3}_\mathrm{sym}$
  \begin{align*}
    \left(\bS-\tilde\bS, \bD-\tilde\bD\right)\geq0
    \quad\text{ for all } (\tilde\bS,\tilde\bD)\in\cG
  \end{align*}
  implies, using~\eqref{A.1} and~\eqref{A.2}, that
  \begin{align}
    \label{A.3}
    \left(
      \bS-\frac{\left(|\tilde\bD|-\delta_*\right)^+}{|\tilde\bD|}\mathcal{S}(|\tilde\bD|)\tilde\bD
    \right)\mathbin{:}\left(
      \bD-\tilde\bD
    \right)\geq0
    \quad\text{ for all }\tilde\bD\in\RR^{3\times 3}_\mathrm{sym}.
  \end{align}
  Taking $\tilde\bD=\bD\pm\lambda\bA$, $\bA$ arbitrary, $\lambda>0$,
  we conclude from~\eqref{A.3} that
  \begin{align*}
    \mp\bA\mathbin{:}
    \left(
      \bS-\frac{\left(|\bD\pm\lambda\bA|-\delta_*\right)^+}{|\bD\pm\lambda\bA|}
          \mathcal{S}(|\bD\pm\lambda\bA|)(\bD\pm\lambda\bA)
    \right)\geq0.
  \end{align*}
  Letting $\lambda\rightarrow0+$, we finally obtain (using continuity
  of the involved functions)
  \begin{align*}
    \left(
      \bS-\frac{\left(|\bD|-\delta_*\right)^+}{|\bD|}
          \mathcal{S}(|\bD|)\bD
    \right)\mathbin{:}\bA=0
    \quad\text{ for all } \bA\in\RR^{3\times 3}_\mathrm{sym}.
  \end{align*}
  Hence $\bS$ and $\bD$ fulfill the right-hand side of~\eqref{A.1}
  and thus $(\bS,\bD)\in\cG$.

  \noindent(iv).
  Assume that $\delta_*>0$.
  For $d\geq0$ we have
  \begin{align} \label{A.4}
    \begin{gathered}
      \begin{multlined}
        \min\Bigl\{1, \tfrac{\bigl(1+(2\delta_*)^2\bigr)^\frac{r-2}{2}}{(2\delta_*)^{r-2}}\Bigr\}
        \mathcal{H}(d-2\delta_*) d^{r-1}
        \leq (1+d^2)^\frac{r-2}{2} d
        \leq (1+d^2)^\frac{r-1}{2}
        \\
        \leq (1+d)^{r-1}
        \leq \Bigl( (2\delta_*)^{-1} \max\{d,2\delta_*\}
                  + \max\{d,2\delta_*\} \Bigr)^{r-1},
      \end{multlined}
      \\
      \begin{multlined}
        d^{q-1}
        \leq (1+d^{r-2}) d
        = d + d^{r-1}
        \\
        \leq (2\delta_*)^{2-q} \bigl( \max\{d,2\delta_*\} \bigr)^{q-1}
           + (2\delta_*)^{r-q} \bigl( \max\{d,2\delta_*\} \bigr)^{q-1}
      \end{multlined}
    \end{gathered}
  \end{align}
  where $q=\max\{r,2\}$.
  Let us define
  \begin{align*}
    q \coloneqq
    \left\{\begin{array}{ll}
      r           & \text{case \eqref{A.2a}}, \\
      \max\{r,2\} & \text{case \eqref{A.2b}}.
    \end{array}\right.
  \end{align*}
  Due to~\eqref{A.4} we have, for the both cases in~\eqref{A.2},
  \begin{align} \label{eq:Sest}
    C_1(\delta_*,r) \mathcal{H}(|\bD|-2\delta_*) |\bD|^{q-1}
    \leq
    \mathcal{S}(|\bD|)|\bD|
    \leq
    C_2(\delta_*,r) \bigl( \max\{|\bD|,2\delta_*\} \bigr)^{q-1}
  \end{align}
  with certain $C_1(\delta_*,r)$, $C_2(\delta_*,r)>0$ independent of $\bD$.
  Notice also that
  \begin{align} \label{eq:plus}
    \tfrac12 \mathcal{H}(|\bD|-2\delta_*)
    \leq \frac{(|\bD|-\delta_*)^+}{|\bD|}
    \leq 1.
  \end{align}
  Now define
  $\bS=\frac{\left(|\bD|-\delta_*\right)^+}{|\bD|}\mathcal{S}(|\bD|)\bD$
  and observe that with the help of the right-wing inequalities
  of~\eqref{eq:Sest} and~\eqref{eq:plus} we obtain
  \begin{multline*}
    |\bD|^q + |\bS|^{q'}
    \leq |\bD|^q + \Bigl( C_2(\delta_*,r)
                          \bigr(\max\{|\bD|,2\delta_*\}\bigl)^{q-1}
                   \Bigr)^{q'}
    \leq C_3(\delta_*,r) \bigl( \max\{|\bD|,2\delta_*\} \bigr)^q
  \end{multline*}
  with certain $C_3(\delta_*,r)>0$ independent of $\bD$ and $\bS$.
  Hence, with the help of the left-wing inequalities
  of~\eqref{eq:Sest} and~\eqref{eq:plus},
  \begin{multline*}
    C_3(\delta_*,r)^{-1} \bigl( |\bD|^q + |\bS|^{q'} \bigr) - (2\delta_*)^q
    \leq \bigl( \max\{|\bD|,2\delta_*\} \bigr)^q - (2\delta_*)^q
    \\
    \leq \mathcal{H}(|\bD|-2\delta_*) |\bD|^q
    \leq 2 C_1(\delta_*,r)^{-1} \frac{(|\bD|-\delta_*)^+}{|\bD|} \mathcal{S}(|\bD|)|\bD|^2
    = 2 C_1(\delta_*,r)^{-1} \bS\mathbin{:}\bD
  \end{multline*}
  which is the last property~($\cG$\ref{Gcoercive}).
  We leave the case $\delta_*=0$ as an~exercise.
\end{proof}

\bibliographystyle{siamplain}

\end{document}